\documentclass[]{jloganal}
\usepackage{amssymb}
\usepackage{graphicx}
\usepackage{pgfplots}
\RequirePackage[utf8]{inputenc}

\arxivreference{}
\arxivpassword{}
\volumenumber{}
\issuenumber{}
\publicationyear{}
\papernumber{}
\startpage{}
\endpage{}
\doi{}
\MR{}
\Zbl{}
\received{}
\revised{}
\accepted{}
\published{}
\publishedonline{}
\proposed{}
\seconded{}
\corresponding{}
\editor{}
\version{}

\theoremstyle{definition}

\begin{document}

\title{A Generalization of the Cantor-Dedekind Continuum\\
with Nilpotent Infinitesimals}
\author{JOS{\'E} ROQUETTE}
\email{jroquet@math.tecnico.ulisboa.pt}
\email{joseroquette@iol.pt}

\begin{abstract}
We introduce a generalization of the Cantor-Dedekind continuum with explicit
infinitesimals. These infinitesimals are used as numbers obeying the same
basic rules as the other elements of the generalized continuum, in
accordance with Leibniz's original intuition, but with an important
difference: their product is null, as the Dutch theologian Bernard
Nieuwentijt sustained, against Leibniz's opinion. The starting-point is the
concept of shadow, and from it we define indiscernibility (the central
concept) and monad. Monads of points have a global-local nature, because in
spite of being infinite-dimensional real affine spaces with the same
cardinal as the whole generalized continuum, they are closed intervals with
length 0. Monads and shadows (initially defined for points) are then
extended to any subset of the new continuum, and their study reveals
interesting results of preservation in the areas of set theory and topology.
All these concepts do not depend on a definition of limit in the new
continuum; yet using them we obtain the basic results of the differential
calculus. Finally, we give two examples illustrating how the global-local
nature of the monad of a real number can be applied to the differential
treatment of certain singularities.
\end{abstract}

\maketitle

\givenname{JOS{\'E}} \surname{ROQUETTE}

\address{Department of Mathematics\\
IST-University of Lisbon\\
\newline Av. Rovisco Pais\\
1049-001 Lisboa\\
PORTUGAL}

\section{Introduction}

Up to 1960, when Abraham Robinson created \textit{Non-standard Analysis}, 
\textit{actual infinitesi-mals}, i.e. infinitesimals considered as numbers,
in the Leibniz's tradition $\left[ 6\right] $, were banished from
mathematical analysis by Weierstrass' $\varepsilon -\delta $ definition of 
\textit{limit} (in the 1850s), except for a minority of mathematicians and
at least one great philosopher (Charles S. Peirce). But physicists and
engineers (and differential geometers such as Sophus Lie, {\'{E}}lie Cartan,
and Hermann Weyl) refused to deprive themselves of the immense heuristic
power of that notion (and rightly so!).\newline
Today, there are two main rigorous theories of actual infinitesimals: 
\textit{Non-standard Analysis} (\textbf{NSA}) $\left[ 4\right] ,\left[ 5%
\right] ,\left[ 8\right] ,\left[ 9\right] ,\left[ 10\right] ,\left[ 11\right] $%
, using nonexplicit invertible infinitesimals, and \textit{Smooth
Infinitesimal Analysis} (\textbf{SIA}) (F.W. Lawvere, in the late 1960s) $%
\left[ 1\right], \left[ 2\right], \left[ 7\right] $, with nilpotent infinitesimals (i.e. infinitesimals $%
\varepsilon $ such that $\varepsilon ^{n}=0$, for some positive integer $n)$%
. But both theories are considered with suspicion by the immense majority of
the mathematical community, and physicists and engineers prefer their strong
intuitions.\newline
The generalization $(\widehat{\mathbb{R}})$ of the \textit{usual
Cantor-Dedekind continuum} $(\mathbb{R})$ we propose, and the ensuing
Calculus, have the following features:

\paragraph*{I}

-- The elements of $\widehat{\mathbb{R}}$, which we call \textit{generalized
real numbers}, are the convergent (in the usual sense) sequences in $\mathbb{%
R}$, and those sequences that converge to $0$ are called \textit{%
infinitesimals} (so \textit{infinitesimals are explicit}). The \textit{shadow%
} of a generalized real number is just its \textit{limit} as a convergent
sequence in $\mathbb{R}$, and from this concept we define a binary relation
on $\widehat{\mathbb{R}}$ that coincides with the identity of the shadows,
and which we call \textit{indiscernibility} ($\approx )$. The \textit{monad}
of a generalized real number $x_{0}$ ($m_{\approx }(x_{0})$) is the set of
all elements of $\widehat{\mathbb{R}}$ that are \textit{indiscernible} 
\textit{from} $x_{0}$. On the set $\widehat{\mathbb{R}}$ we define \textit{addition} term by term, but \textit{multiplication} and \textit{ordering} are introduced in a different manner, using the concept of shadow. We obtain
an ordered ring extension of $\mathbb{R}$ (though it is important to take into account $\mathbf{f_2})$ below); moreover, the quotient of $\widehat{\mathbb{R}}$ by $\approx $ is an ordered
field isomorphic to $\mathbb{R}$.\newline
Although we can embed $\mathbb{R}$ in $\widehat{\mathbb{R}}$ (through the
mapping $\xi \mapsto (\xi )$, where $(\xi )$ is the constant sequence
determined by the real number $\xi $), we must emphasize two features of $%
\widehat{\mathbb{R}}$ that are absent from $\mathbb{R}$:

$\mathbf{f_1})$ \textit{The product of two nonnull generalized} \textit{real
numbers or the square of a nonnull generalized real number may be null} (if and only all the factors are \textit{infinitesimal}).%
\newline
$\mathbf{f_2})$ \textit{Strict ordering is defined on} $\widehat{\mathbb{R}}$
\textit{except inside the monads} (as it should be expected, since
the elements of the monad of a generalized real number are indiscernible).
So we have this version of the usual \textit{trichotomy property}: 
\[
\left(\forall x,y\in \widehat{\mathbb{R}}) (x < y \vee x \approx y \vee y <
x\right). 
\]

\paragraph*{II}

-- We work in two \textit{modes}:

\textit{\underbar{The  mode  of  potentiality}}, i.e. the totality of
notions and concepts that can be defined within the structure $\mathbb{R}$.%
\newline
\textit{\underbar{The mode of actuality}}, i.e. the totality of notions and
concepts that can be defined within the structure $\widehat{\mathbb{R}},$
with the exception of any definition of \textit{limit}.

We use the \textit{mode of potentiality} emphasizing the usual definition of 
\textit{limit}, but in the \textit{ mode} \textit{of actuality}, in the
absence of such a definition, we must introduce the fundamental concepts of 
\textit{generalized real number}, and \textit{shadow}, in the \textit{mode
of potentiality}. Nevertheless, we must stress that this translation is only
made for the sake of defin\nolinebreak ition: once defined, the two
fundamental concepts are used in the \textit{mode of actuality}. Every
notion or concept in \textit{the mode of actuality} could be translated into 
\textit{the mode of potentiality}, but then we would renounce the intuitive
and computational power of \textit{actual methods}.\newline
Our work in these two \textit{modes}, sometimes simultaneously (as in the
definition of \textit{differentiability}), reflects our conviction that a
concept of \textit{actual infinitesimal} and a definition of \textit{limit}
are both necessary to a Calculus fit, not only for mathematicians, but also
for experimental scientists.

\paragraph*{III}

-- Each generalized real number \textit{x} is indiscernible from exactly one
real number: its \textit{shadow}, which we denote by $\sigma x$. In fact, 
\textit{each generalized real number x admits a unique decomposition as the
sum of a real number (its shadow) and an infinitesimal}. We denote this
infinitesimal by \textit{dx}, and we call it the \textit{differential }of%
\textit{\ x}. So we have, for each $x\in \widehat{\mathbb{R}}$, the unique
decomposition, which we call the $\sigma +d$ \textit{decomposition}: 
\[
x=\sigma x+dx. 
\]%
For each $x\in \widehat{\mathbb{R}}$, and $\xi \in \mathbb{R}$, we have, as
a direct consequence of the $\sigma +d$ \textit{decomposition} (and we
stress its uniqueness!): 
\[
\sigma \xi =\xi , 
\]%
\[
d\xi =0, 
\]%
\[
\sigma (dx)=0, 
\]%
\[
d(dx)=dx, 
\]%
\[
\sigma (\xi +dx)=\xi , 
\]%
\[
d(\xi +dx)=dx. 
\]%
Although we do not use a definition of \textit{limit} in $\widehat{\mathbb{R}%
}$, we can easily derive the \textit{basic algebraic rules of differentiation%
}, using the $\sigma +d$ \textit{decomposition}.

\paragraph*{IV}

-- For each subset $\widehat{A}$ of $\widehat{\mathbb{R}}$, we define its 
\textit{monad} ($m_{\approx }(\widehat{A})$) and \textit{shadow} ($\sigma (%
\widehat{A})$), and we obtain interesting \textit{set-theoretic} and \textit{%
topological} results of preservation.\newline
The \textit{intervals} in $\widehat{\mathbb{R}}$ are simply the monads of
the corresponding intervals in $\mathbb{R}$, and the \textit{length} of
those that are \textit{bounded} (i.e. those intervals in $\widehat{\mathbb{R}%
}$ that are monads of bounded intervals in $\mathbb{R})$ is the same as the length of
their originals in $\mathbb{R}$; for instance, the \textit{bounded open} and
the \textit{bounded closed intervals} in $\widehat{\mathbb{R}}$ are 
\[
\widehat{\left] \alpha ,\beta \right[ }:=m_{\approx }(\left] \alpha ,\beta %
\right[ ), 
\]%
\[
\widehat{\left[ \alpha ,\beta \right] }{{\ :=m_{\approx }(\left[ \alpha
,\beta \right] }}), 
\]%
respectively, where $\alpha ,\beta \in \mathbb{R}$, and $\alpha \leq \beta $
(their \textit{length} is $\beta -\alpha )$.

Intervals in $\widehat{\mathbb{R}}$ do not have \textit{pointlike}
extremities, and this feature is reminiscent of \textit{Stoic philosophical
view} about \textit{segments of Space} or \textit{Time} $\left[ 12\right] $;
for instance, if $\alpha ,\beta ,\gamma \in \mathbb{R}$, and $\alpha \leq
\beta \leq \gamma $, then 
\[
m_{\approx }(\alpha )=\widehat{\left[ \alpha ,\alpha \right] }
\]%
\[
m_{\approx }(\alpha ),m_{\approx }(\beta )\subseteq \widehat{\left[ \alpha
,\beta \right] },
\]%
\[
\widehat{\left[ \alpha ,\beta \right] }\cap \widehat{\left[ \beta ,\gamma %
\right] }\hspace{2pt}{{=m_{\approx }(\beta }}).
\]

\paragraph*{V}

-- The monad of each generalized real number \textit{x} has a \textit{%
global-local} nature since it is an infinite-dimensional real affine space
with the same cardinal as $\widehat{\mathbb{R}}$ (more precisely, $%
\left|m_{\approx} \left(x\right)\right|=\left|\widehat{\mathbb{R}}%
\right|=2^{\aleph_0}$), yet it is also a closed interval of length 0 (it is
easy to prove that $m_{\approx} (x) = m_{\approx} \left( \sigma x \right)$,
so $m_{\approx}(x)= \widehat{\left[\sigma x,\sigma x \right]}$ ).

We use this \textit{dual} nature in two examples of differential treatment
of \textit{singularities}.\newline

\paragraph*{VI}

-- For each function $\phi :I\rightarrow \mathbb{R}$, where \textit{I} is an
open interval in $\mathbb{R}$, its \textit{indiscernible extensions} are the
functions $f:m_{\approx }\left( I\right) \rightarrow \widehat{\mathbb{R}}$
such that 
\[
f(\sigma x)=\phi (\sigma x), 
\]%
\[
f(x)\approx \phi (\sigma x). 
\]%
If $\xi _{0}\in I$, and $f:m_{\approx }(I)\rightarrow \widehat{\mathbb{R}}~$
is an indiscernible extension of $\phi $, then $f$ is said to be \textit{%
differentiable }at $\xi _{0}$ iff there exists a real number $\alpha $ such
that 
\[
\left( \forall x\in m_{\approx }(\xi _{0})\right) f(x)=\phi \left( \xi
_{0}\right) +\alpha dx, 
\]%
with the proviso that $\alpha :=\lim_{\xi \rightarrow \xi _{0}}\frac{\phi
\left( \xi \right) -\phi (\xi _{0})}{\xi -\xi _{0}}$, when such limit exists
in $\mathbb{R}$.

$\alpha $ (which is \textit{unique}) is said to be the \textit{derivative }%
of $f$ at $x$, for each $x\in m_{\approx }(\xi _{0})$, and we denote it by $%
f^{\prime }(x)$, as usual.

So we have, when \textit{f} is differentiable at $\xi _{0}$:\newline
$\mathbf{d_{1}})$ If $x\in m_{\approx }(\xi _{0})$, then $f^{\prime
}(x)=f^{\prime }(\xi _{0})$.\newline
$\mathbf{d_{2}})$ For each $x\in m_{\approx }(\xi _{0})$, 
\[
f(x)=f(\xi _{0})+f^{\prime }(\xi _{0})dx.
\]

This is the expression, in analytical terms, of the \textit{geometric idea}
associated with the concept of \textit{differentiability,} according to
Leibniz primeval conception:

\textit{If }$f$\textit{\ is differentiable at }$\xi _{0}$\textit{, then the
graph of }$f$\textit{\ coincides locally }(\textit{i.e. for infinitesimal
increments of the argument around }$\xi _{0}$)\textit{\ with its tangent at
the point }$\left( \xi _{0},\hspace{1pt}f(\xi _{0})\right) $\textit{.}

Notice that if $\lim_{\xi \rightarrow \xi _{0}}\frac{\phi \left( \xi \right)
-\phi (\xi _{0})}{\xi -\xi _{0}}$ exists in $\mathbb{R}$ (i.e. $\phi $ is
differentiable at $\xi _{0}$, in the usual\smallskip \newline
sense) \hspace{1pt}and $\hspace{1pt}f$ \hspace{1pt}is \hspace{1pt}%
differentiable \hspace{1pt}at $\hspace{1pt}\xi _{0}$, \hspace{1pt}then $%
\hspace{1pt}f^{\prime }(\xi _{0})$ \hspace{1pt}is \hspace{1pt}identical 
\hspace{1pt}with \hspace{1pt}this \hspace{1pt}limit; \hspace{1pt}however,
\smallskip \vspace{-6pt}

$f^{\prime }(\xi _{0})$ may exist in the absence of $\lim_{\xi \rightarrow
\xi _{0}}\frac{\phi \left( \xi \right) -\phi (\xi _{0})}{\xi -\xi _{0}}$ ,
as it is the case for $\xi _{0}:=0,$ and $\hspace{2pt}\phi :\mathbb{R}%
\rightarrow \mathbb{R}$, $\hspace{1pt}f:\widehat{\mathbb{R}}\rightarrow 
\widehat{\mathbb{R}}$ \hspace{2pt}defined \hspace{2pt}by $\hspace{2pt}\phi
(\xi ):=\left\vert \xi \right\vert ,$ $f(x):=\left\{ 
\begin{array}{l}
x,\text{ if \ }x>0 \\ 
0,\text{ if \ }x\in m_{\approx }(0) \\ 
-x,\text{ if \ }x<0%
\end{array}%
\right. $ (clearly, $f^{\prime }(0)=0).$

Keeping in mind that the \textit{derivatives} are always associated with 
\textit{indiscernible extensions}, and using the definition, we obtain not
only the \textit{algebraic rules of derivation}, but also fundamental
theorems like the \textbf{Chain Rule}, the \textbf{Inverse Function Theorem}%
, the \textbf{Mean Value Theorem}, and \textbf{Taylor's Theorem}.

If $\lim_{\xi \rightarrow \xi _{0}}\frac{\phi \left( \xi \right) -\phi (\xi
_{0})}{\xi -\xi _{0}}$ exists, for each $\xi _{0}\in I$, then, among the
infinity of indiscernible \vspace{-3pt}

extensions of $\phi $, there exists exactly one that is differentiable at
each $\xi _{0}\in I$; we call

\vspace{-3pt}this function the \textit{natural indiscernible extension of} $%
\phi $, and we denote it by $\hat{\phi}$.\newline
So $\hat{\phi}:m_{\approx }(I)\rightarrow ~\widehat{\mathbb{R}}$ is the
function defined by 
\[
\hat{\phi}(x):=\phi (\xi _{0})+\lambda _{\phi }(\xi _{0})dx, 
\]%
where $\lambda _{\phi }(\xi _{0})$ denotes $\lim_{\xi \rightarrow \xi _{0}}%
\frac{\phi \left( \xi \right) -\phi (\xi _{0})}{\xi -\xi _{0}}$.

The concept of \textit{natural indiscernible extension} provides a rule for
the definition of the \textit{analogues} (and \textit{extensions}) of the
usual functions of Real Analysis. For instance, the natural indiscernible
extensions of exp, log, sin, cos, are the functions (where $\widehat{\mathbb{%
R}}^{+}$ is the set of positive generalized real numbers): 
\[
\widehat{\exp }:\widehat{\mathbb{R}}\rightarrow \widehat{\mathbb{R}},
\]%
\[
\widehat{\log }:\widehat{\mathbb{R}}^{+}\rightarrow \widehat{\mathbb{R}},
\]%
\[
\widehat{\sin }:\widehat{\mathbb{R}}\rightarrow \widehat{\mathbb{R}},
\]%
\[
\widehat{\cos }:\widehat{\mathbb{R}}\rightarrow \widehat{\mathbb{R}},
\]%
defined by 
\[
\widehat{\exp }(x):=\exp (\sigma x)+\exp (\sigma x)dx,
\]%
\[
\widehat{\log }(x):=\log (\sigma x)+\frac{1}{\sigma x}dx,
\]%
\[
\widehat{\sin }(x):=\sin (\sigma x)+\cos (\sigma x)dx,
\]%
\[
\widehat{\cos }(x):=\cos (\sigma x)-\sin (\sigma x)dx.
\]%
We show that these functions have the same basic properties as the usual
ones, and we obtain, \textit{rigorously}, some identities that physicists
and engineers often use \textit{intuitively}. For example (since $\sigma
(dx)=0$, and $\sigma ({1}+dx)=1$, as seen in \textbf{III}) : 
\[
\widehat{\exp }(dx)=\exp (\sigma (dx))+\exp (\sigma (dx))dx=\exp (0)+\exp
(0)dx=1+dx,
\]%
\[
\widehat{\log }\left( 1+dx\right) =\log (\sigma ({1}+dx))+\frac{1}{\sigma
(1+dx)}dx=\log (1)+dx=dx,
\]%
\[
\widehat{\sin }(dx)=\sin (\sigma (dx))+\cos (\sigma (dx))dx=\sin (0)+\cos
(0)dx=dx,
\]%
\[
\widehat{\cos }(dx)=\cos (\sigma (dx))-\sin (\sigma (dx))dx=\cos (0)-\sin
(0)dx=1.
\]%
\newline
\newline

\section{The Generalized Real Numbers}

Let ($\mathbb{R},<,+,\cdot ,0,1$) be a model of the usual real number system
axioms (in any of the equivalent formulations of most calculus textbooks),
and let $\widehat{\mathbb{R}}$ be the set of all sequences $x=(\xi _{n})$ in 
$\mathbb{R}$ that are convergent for the usual absolute value in ($\mathbb{R}%
,<$,+,$\cdot $,0,1). We refer to ($\mathbb{R},<,+,\cdot ,0,1$) as the 
\textit{Cantor-Dedekind continuum}.

\bigskip

\textbf{Definition 2.1 \ }Let \textit{x,y} $\in \mathbb{\widehat R}$.\newline
If $\lim x$ is the usual limit of \textit{x} in ($\mathbb{R},<,+,\cdot ,0,1$%
), then we call the constant sequence ($\lim x$), the \textit{shadow} of 
\textit{x}, and we denote it by $\sigma x$.\newline
$x$ is said to be \textit{indiscernible from} $y$, and we denote it by $%
x\approx y$, iff $x$ and $y$ have the same shadow.\newline
$x$ is said to be an \textit{infinitesimal} iff $x$ is indiscernible from
the constant sequence (0).\newline
The \textit{monad} of $x$, denoted by $m_{\approx }(x),$ is the set of all $%
y\in \mathbb{\widehat R}$ such that $y$ is indiscernible from $x$.\newline
So $m_{\approx }((0))$ is the set of all infinitesimals.

\bigskip

Clearly, the \textit{indiscernibility} relation, $\approx $ , is an
equivalence relation on $\widehat{\mathbb{R}}$, and if \textit{x} is an
element of $\widehat{\mathbb{R}}$, then its equivalence class for $\approx $
is $m_{\approx }(x)$. \textit{Indiscernibility} is the first and more
important binary relation defined on $\widehat{\mathbb{R}}$.

\bigskip

The next definition introduces a ring structure for $\widehat{\mathbb{R}}$
with a kind of linear ordering.

\bigskip

\textbf{Definition 2.2 \ }On the set $\widehat{\mathbb{R}}$, we consider two
binary operations, denoted by $\hat{{+}}$ and $\hat{\cdot}$, and called 
\textit{addition} and \textit{multiplication}, respectively. If $x=(\xi
_{n}) $ and $y=(\eta _{n})$ are elements of $\widehat{\mathbb{R}}$, then
these operations are defined by 
\[
x\text{ }\hat{+}\text{ }y:=(\xi _{n}+\eta _{n}), 
\]%
\[
x\text{ }\hat{\cdot}\text{ }y:=(\lim x\cdot \eta _{n}+\lim y\cdot \xi
_{n}-\lim x\cdot \lim y), 
\]%
where at the right-hand of the previous identities we consider the obvious
operations on $\mathbb{R}$ (clearly, $ x\text{ }\hat{+}\text{ }y, x\text{ }\hat{\cdot}\text{ }y \in \widehat{\mathbb{R}} $ and $ \lim (x\text{ }\hat{+}\text{ }y)=\lim x + \lim y, \lim (x\text{ }\hat{\cdot}\text{ }y)=\lim x \cdot \lim y)$.\newline
We say that \textit{x is less than y}, and we denote it by $x$ $\hat{<}$ $y$%
, iff $\lim x<\lim y$, and reciprocally, we say that \textit{x is greater
than y}, and we denote it by $x$ $\hat{>}$ $y$, iff $y$ $\hat{<}$ $x$, where
in $\lim x<\lim y$ we consider the usual linear ordering on $\mathbb{R}$.%
\newline
The elements of $\widehat{\mathbb{R}}^{+}:=\left\{ x\in \mathbb{R}|x\text{ }%
\hat{>}\text{ }(0)\right\} $ and $\widehat{\mathbb{R}}^{-}:=\left\{ x\in 
\mathbb{R}|x\text{ }\hat{<}\text{ }(0)\right\} $ will be called \textit{%
positive} and \textit{negative}, respectively.

\bigskip

\textbf{Proposition 2.3 \ a)} ($\widehat{\mathbb{R}},\widehat{<},\widehat{+},%
\hat{\cdot},(0),(1)$) is a commutative ring with the constant sequences (0)
and (1) as \textit{zero element} and \textit{identity element}, respectively.%
\newline
\textbf{b)} The \textit{shadow mapping} $\sigma :\widehat{\mathbb{R}}%
\rightarrow \widehat{\mathbb{R}}$, defined by $\sigma (x):=\sigma x,$ is an
idempotent ring endomorphism, i.e. 
\[
(\forall x\in \widehat{\mathbb{R}})\sigma (\sigma x)=\sigma x;
\]%
\[
(\forall x,y\in \widehat{\mathbb{R}})\sigma (x\text{ }\hat{+}\text{ }%
y)=\sigma x\text{ }\hat{+}\text{ }\sigma y,
\]%
\[
(\forall x,y\in \widehat{\mathbb{R}})\sigma (x\text{ }\hat{\cdot}\text{ }%
y)=\sigma x\text{ }\hat{\cdot}\text{ }\sigma y,
\]%
\[
\sigma (1)=(1).
\]%
Furthermore, 
\[
Ker(\sigma ):=\left\{ x\in \widehat{\mathbb{R}}|\sigma x=(0)\right\}
=m_{\approx }((0)),
\]%
\[
\sigma (\widehat{\mathbb{R}})=\left\{ y\in \widehat{\mathbb{R}}|\text{y is a
constant sequence}\right\} .
\]%
\textbf{c) }$m_{\approx }((0))$ is a nonnull ideal, so the sum of
infinitesimals is an infinitesimal, the additive inverse of an infinitesimal
is also an infinitesimal, (0) is an infinitesimal, the product of an element
of $\widehat{\mathbb{R}}$ and an infinitesimal is still an infinitesimal,
and there is a nonnull infinitesimal.\newline
\textbf{d)} The product of infinitesimals is always null, i.e. 
\[
(\forall x,y\in m_{\approx }((0)))\hspace{1pt}x\text{ }\hat{\cdot}\text{ }%
y=(0).
\]%
In particular, each infinitesimal is nilpotent, since $x$ $\hat{\cdot}$ $%
x=(0)$, for each $x\in m_{\approx }((0))$.\newline
\textbf{e)} An element of $\widehat{\mathbb{R}}$ has a \textit{%
multiplicative inverse} iff it is not an infinitesimal.\newline
\textbf{f)} If $x,y,z\in \widehat{\mathbb{R}}$, then 
\[
\lnot (x\text{ }\hat{<}\text{ }x),
\]%
\[
x\text{ }\hat{<}\text{ }y\wedge y\text{ }\hat{<}\text{ }z\Rightarrow x\text{ 
}\hat{<}\text{ }z,
\]%
\[
x\text{ }\hat{<}\text{ }y\vee x\approx y\vee y\text{ }\hat{<}\text{ }x,
\]%
\[
x\text{ }\hat{<}\text{ }y\Rightarrow x\text{ }\hat{+}\text{ }z\text{ }\hat{<}\text{ }y\text{ }\hat{+}\text{ }z,
\]%
\[
x\text{ }\hat{<}\text{ }y\wedge z\text{ }\hat{>}(0)\Rightarrow x\text{ }\hat{\cdot}\text{ }z\text{ }\hat{<}\text{ }y\text{ }\hat{\cdot}\text{ }z.
\]%
So, if we adopt the version of the usual \textit{trichotomy property} expressed by the third formula above, then ($\widehat{\mathbb{R}},%
\hat{<},\hat{+},\hat{\cdot},(0),(1)$) may be considered an ordered ring \nolinebreak .%
\newline
\textbf{g)} ($\widehat{\mathbb{R}},\hat{<},\hat{+},\hat{\cdot},(0),(1)$) is 
\textit{archimedean}, i.e. 
\[
(\forall x,y\in \mathbb{\widehat R})(x\text{ }\hat{>}\text{ }(0)\Rightarrow (\exists
m\in \mathbb{N})\widehat{m}x\text{ }\hat{>}\text{ }y),
\]%
where $\widehat{m}x$ abbreviates $x_{1}$ $\hat{+}$ $x_{2}$ $\hat{+}$ $\ldots 
$ $\hat{+}$ $x_{m}$, when $x_{1}=x_{2}=\ldots =x_{m}=x$ (assum\nolinebreak
ing $\widehat{{1}}x=x$).\newline
\textbf{h)} The mapping $\ast :\mathbb{R}\rightarrow \sigma (\widehat{%
\mathbb{R}}),$ defined by $\ast (\xi ):=(\xi )$, where ($\xi )$ is the usual
constant sequence determined by $\xi $, is a ring isomorphism of ($\mathbb{R}%
,<,+,\cdot ,0,1$) onto ($\sigma (\widehat{\mathbb{R}}),\hat{<},\hat{+},\hat{%
\cdot},(0),(1)$), and 
\[
(\forall \xi ,\eta \in \mathbb{R})(\xi <\eta \Leftrightarrow \ast (\xi )%
\text{ }\hat{<}\ast (\eta )).
\]%
So, using $\ast $, we can embed ($\mathbb{R},<,+,\cdot ,0,1$) in ($\widehat{%
\mathbb{R}},\hat{<},\hat{+},\hat{\cdot},(0),(1)$).

\bigskip

\textbf{Proof \ a)} \ Only the proofs of the associative property of multiplication and the distributive property of multiplication over addition offer some (slight) difficulty.\newpage If $x=(\xi_{n}),y=(\eta_{n}),z=(\zeta_{n})\in \widehat{\mathbb{R}}$, then\ \newline $ (x\ \widehat{\cdot}\ y)\ \widehat{\cdot}\ z =(\lim x\ \cdot\ \eta_{n}\ +\ \lim\ y\ \cdot\ \xi_{n}\ -\ \lim\ x\ \cdot\ \lim\ y)\ \widehat{\cdot}\ z=(\lim\ x\ \cdot\ \lim\ y\ \cdot\ \zeta_{n}\ +\ \newline +\ \lim\ z\ \cdot\ \lim\ x\ \cdot\ \eta_{n}\ +\ \lim\ z\ \cdot\ \lim\ y\ \cdot\ \xi_{n}\ -\ \lim\ z\ \cdot\ \lim\ x\ \cdot\ \lim\ y\ -\ \lim\ x\ \cdot\ \lim\ y\ \cdot\ \lim\ z)=\newline =(\lim\ x\ \cdot\ \lim\ y\ \cdot\ \zeta_{n}\ +\ \lim\ z\ \cdot\ \lim\ x\ \cdot\ \eta_{n}\ +\ \lim\ z\ \cdot\ \lim\ y\ \cdot\ \xi_{n}\ -\ 2\lim\ x\ \cdot\ \lim\ y\ \cdot\ \lim\ z),\ \newline x\ \widehat{\cdot}\ (y\ \widehat{\cdot}\ z)=x\ \widehat{\cdot}\ (\lim\ y\ \cdot\ \zeta_{n}\ +\ \lim\ z\ \cdot\ \eta_{n}\ -\lim\ y\ \cdot\ \lim\ z)=(\lim\ x\ \cdot\ \lim\ y\ \cdot\ \zeta_{n}\ +\newline +\ \lim\ x\ \cdot\ \lim\ z\ \cdot\ \eta_{n}\ -\ \lim\ x\ \cdot\ \lim\ y\ \cdot\ \lim\ z\ +\ \lim\ y\ \cdot\ \lim\ z\ \cdot\ \xi_{n}\ -\ \lim\ x\ \cdot\ \lim\ y\ \cdot\ \lim\ z)=\newline =(\lim\ x\ \cdot\ \lim\ y\ \cdot\ \zeta_{n}\ +\ \lim\ x\ \cdot\ \lim\ z\ \cdot\ \eta_{n}\ +\ \lim\ y\ \cdot\ \lim\ z\ \cdot\ \xi_{n}\ -\ 2 \lim\ x\ \cdot\ \lim\ y\ \cdot\ \lim\ z)=\newline =(x\ \widehat{\cdot}\ y)\ \widehat{\cdot}\ z;\ \newline x\ \widehat{\cdot}\ (y\ \widehat{+}\ z)=x\ \widehat{\cdot}\ (\eta_{n}\ +\ \zeta_{n})=(\lim\ x\ \cdot\ \eta_{n}\ +\ \lim\ x\ \cdot\ \zeta_{n}\ +\ \lim\ y\ \cdot\ \xi_{n}\ +\newline +\  \lim\ z\ \cdot\ \xi_{n}\ -\ \lim\ x\ \cdot\ \lim\ y\ -\ \lim\ x\ \cdot\ \lim\ z)=(\lim\ x\ \cdot\ \eta_{n}\ +\ \lim\ y\ \cdot\ \xi_{n}\ -\ \newline -\ \lim\ x\ \cdot\ \lim\ y)\ \widehat{+}\ ( \lim\ x\ \cdot\ \zeta_{n}\ +\ \lim\ z\ \cdot\ \xi_{n}\ -\ \lim\ x\ \cdot\ \lim\ z)=(x\ \widehat{\cdot}\ y)\ \widehat{+}\ (x\ \widehat{\cdot}\ z).$ \newline
\textbf{b)} is an immediate consequence of the usual algebraic properties of limits, and \textbf{c)}, \textbf{d)}  follow easily from \textbf{a)}, \textbf{b)}.\newline
\textbf{e)}\ If $x=\left( \xi _{n}\right) \in \widehat{%
%TCIMACRO{\U{211d} }%
%BeginExpansion
\mathbb{R}
%EndExpansion
}$ and $x$ is not infinitesimal, then a direct calculation shows that%
\[
x\text{ }\widehat{\cdot }\left( \frac{1}{\lim x}-\frac{\xi _{n}-\lim x}{%
\left( \lim x\right) ^{2}}\right) =\left( 1\right) ; 
\]%
\newline
so, since multiplication on $\widehat{%
%TCIMACRO{\U{211d} }%
%BeginExpansion
\mathbb{R}
%EndExpansion
}$ is associative, commutative, and $\left( 1\right) $ is its identity
element, $\left( \frac{1}{\lim x}-\frac{\xi _{n}-\lim x}{\left( \lim
x\right) ^{2}}\right) $ is the \textit{multiplicative} \textit{inverse} of $%
x=\left( \xi _{n}\right) .$\newline 
If $ x $ is infinitesimal, then we have (see \textbf{a)} and \textbf{b)}), for each $ y\in \widehat{\mathbb{R}}$:
\[ 
 \sigma(x\ \widehat{\cdot}\ y)=\sigma x\ \widehat{\cdot}\ \sigma y=\newline =(0)\ \widehat{\cdot}\ \sigma y=(0)\neq(1),  
\]
and we conclude that $ x $ is not invertible.\newline
Finally, \textbf{f)}, \textbf{g)}, \textbf{h)} admit a quite straightforward proof.  
$\blacksquare $

\bigskip

\textbf{Remark 2.4 \ }In accordance with \textbf{proposition 2.3 h)}, we
identify $\mathbb{R}$ with $\sigma (\widehat{\mathbb{R}})$ and $\xi $ with ($%
\xi $), for each $\xi \in \mathbb{R}$. For instance, we identify 0 with the
infinite sequence (0) and, for each $x\in \widehat{\mathbb{R}}$, $\xi \in 
\mathbb{R}$, we identify $\lim x$ with $\sigma x$ and $\xi $ with $\sigma
\xi $. Furthermore, from now on we shall use the symbols $+$, $\cdot $, $<$
not only for the usual addition, multiplication and linear ordering on $%
\mathbb{R}$, but also for the corresponding binary operations and relation $%
\hat{+},\hat{\cdot},\hat{<}$ on $\widehat{\mathbb{R}}$, and we shall even
drop the symbol $\cdot $ in most formulas. For example, revisiting part of 
\textbf{definition 2.2}, we have, for each $x,y\in \widehat{%
%TCIMACRO{\U{211d} }%
%BeginExpansion
\mathbb{R}
%EndExpansion
}:$%
\[
x<y:\Leftrightarrow \sigma x<\sigma y.
\]%
For the additive and multiplicative powers, we simply write $mx$ and $x^{m}$
instead of $\widehat{m}x$ and $x^{\widehat{m}}$ (where $x^{\widehat{m}}$
abbreviates $x_{1}$ $\hat{\cdot}$ $x_{2}$ $\hat{\cdot}$ $\ldots $ $\hat{\cdot%
}$ $x_{m}$, when $x_{1}=x_{2}=\ldots =x_{m}=x$ (assuming $x^{\widehat{1}}=x$%
)), respectively.

In the spirit of these identifications and notational simplifications,
notice that if $\xi \in \mathbb{R}$ and $x\in \widehat{\mathbb{R}}$, then $%
\xi x$ (previously denoted by $(\xi )$ $\hat{\cdot}$ $x)$ coincides with the
result of the scalar multiplication of the real number $\xi $ by the
sequence \textit{x}.

If $x,y\in \widehat{\mathbb{R}}$ and $x$ is not an infinitesimal, then we denote the \textit{multiplicative
inverse} of \textit{x} by $x^{-1}$ or $\frac{1}{x}$; so $\frac{1}{x}=\left( \frac{1}{\lim x}-\frac{\xi _{n}-\lim x}{(\lim
x)^{2}}\right)$. We also denote $yx^{-1}$
(\textit{the quotient }of\textit{\ y }by\textit{\ x}) by $\frac{y}{x}$, as
usual.

We maintain the general designation of \textit{real numbers} for the elements of $\mathbb{R}$ and call the elements of $\widehat{\mathbb{R}}$ \textit{generalized real numbers}.

\bigskip

Let us see some explicit \textit{generalized real numbers} (by \textit{explicit} we mean \textit{unambiguously defined  as a convergent sequence of real numbers}):

\bigskip

\textbf{Example 2.5 \ 1)} \ The \ eventually \ null \ sequences \ $%
(1,0,0,0,\ldots )$, \ $(0,1,0,0,0,\ldots )$,
($0,0,1,0,0,0,\ldots $), $\ldots $ are nonnull infinitesimal elements of $%
\mathbb{\widehat R}$. So \textit{we can exhibit nonnull infinitesimals}.\newline
\textbf{2)} Let $\xi _{0}$ be a nonnull real number. Then:\newline
The sequences $(0,\xi _{0},\xi _{0},\xi _{0},\ldots )$,$(0,0,\xi _{0},\xi
_{0},\xi _{0},\ldots )$, $(0,0,0,\xi _{0},\xi _{0},\xi _{0},\ldots )$,$%
\ldots $ are different elements of $m_{\approx }(\xi _{0})\backslash \{\xi
_{0}\}$.

\bigskip

In the next proposition, which admits a simple proof, \textbf{e)} and 
\textbf{f)} are particularly important.

\bigskip

\textbf{Proposition 2.6 \ a)} $\left( \forall x\in \widehat{\mathbb{R}}%
\right) \left( x=\sigma x\Leftrightarrow x\in \mathbb{R}\right) $.\newline
\textbf{b)} ($\forall \xi ,\eta \in \mathbb{R})\left( \xi \approx \eta
\Leftrightarrow \xi =\eta \right) $.\newline
\textbf{c)} $\mathbb{R}\cap m_{\approx }(0)=\{0\}$.\newline
\textbf{d)} Infinitesimals are not comparable with respect to the binary
relation $<$ on $\widehat{\mathbb{R}}$, i.e. if $\hat{\varepsilon}$ and $%
\hat{\delta }$ are infinitesimals, then 
\[
\lnot (\hat{\varepsilon}<\hat{\delta})\wedge \lnot (\hat{\delta}<\hat{\varepsilon}%
), 
\]%
\textbf{e)} \textit{An infinitesimal is less than any positive generalized
real number and greater than any negative generalized real number}, i.e. if $%
\hat{\varepsilon}$ is an infinitesimal, than 
\[
(\forall x\in \widehat{\mathbb{R}}^{+})\text{ }\hat{\varepsilon}<x, 
\]%
\[
(\forall y\in \widehat{\mathbb{R}}^{-})\text{ }\hat{\varepsilon}>y. 
\]%
In particular: 
\[
(\forall \xi \in {\mathbb{R}}^{+})\text{ }\hat{\varepsilon}<\xi , 
\]%
\[
(\forall \eta \in {\mathbb{R}}^{-})\text{ }\hat{\varepsilon}>\eta , 
\]%
where ${\mathbb{R}}^{+}$ and ${\mathbb{R}}^{-}$ are the usual sets of
(strictly) positive and (strictly) negative real numbers, respectively
(notice that ${\mathbb{R}}^{+}\hspace{-1pt}\subseteq \widehat{\mathbb{R}}%
^{+} $ and ${\mathbb{R}}^{-}\hspace{-1pt}\subseteq {\widehat{\mathbb{R}}}%
^{-} $, by \textbf{proposition 2.3 h)}).\newline
\textbf{f)} \textit{Each generalized real number is indiscernible from
exactly one real number: its shadow}, i.e. 
\[
(\forall x\in \mathbb{\widehat R})(x\approx \sigma x\wedge (\forall \xi \in \mathbb{R}%
)(x\approx \xi \Rightarrow \xi =\sigma x)). 
\]

\bigskip

\section{The $\protect\sigma +d$ Decomposition}

As a direct consequence of \textbf{proposition 2.3 a)}, \textbf{b)}, we have:

\bigskip

\textbf{Proposition 3.1} \ If \textit{x} is a generalized real number, then
there is a unique infinitesimal $\hat{\varepsilon}(x)$ such that 
\[
x=\sigma x+\hat{\varepsilon}(x). 
\]

\textbf{Definition 3.2} \ If \textit{x} is a generalized real number, then
we denote $\hat{\varepsilon}(x)$ by \textit{dx}, and we call it the \textit{%
differential} of \textit{x}.

\bigskip

\textbf{Proposition 3.3} \ If \textit{x} is a generalized real number then $%
x=\sigma x+dx$ is the unique decomposition of \textit{x} as the sum of a
real number and an infinitesimal.

\bigskip

\textbf{Proof.} We just have to use \textbf{proposition 2.3 a)}, \textbf{c)}%
, \textbf{proposition 2.6 c)}, \textbf{proposition 3.1}, and, of course, 
\textbf{definition 3.2}. $\blacksquare $

\bigskip

We call the decomposition stated by the previous proposition, the $\sigma +d$
\textit{decomposition}. Notice that the \textit{differential} of a
generalized real number \textit{x} is already inlaid in $\mathit{x}$, and
since $\sigma x$ and $dx$ are a constant sequence and a sequence converging
to 0, in $\mathbb{R}$, we are entitled to express the following intuition: a
generalized real number has a unique decomposition as the sum of a \textit{%
static part} (its \textit{shadow}) and a \textit{dynamic part} (its \textit{%
differential}).

\bigskip

\pagebreak Clearly:

\bigskip

\textbf{Corollary 3.4 \ a)} ($\forall x\in \widehat{\mathbb{R}}%
)(dx=0\Leftrightarrow x\in \mathbb{R})$.\newline
\textbf{b)} ($\forall x\in \widehat{\mathbb{R}})(x=dx\Leftrightarrow
x\approx 0)$.\newline
\textbf{c)} ($\forall x\in \widehat{\mathbb{R}})\hspace{1pt}d(dx)=dx$.

\bigskip

The following lemma is the key to obtain the basic algebraic rules of
differentiation.

\bigskip

\textbf{Lemma 3.5 \ a)} If $x,y\in \widehat{\mathbb{R}}$, then 
\[
x+y=\sigma x+\sigma y+dx+dy,\vspace{-3pt} 
\]%
\[
x-y=\sigma x-\sigma y+dx-dy. 
\]%
\begin{center}
\end{center}\textbf{b)} If $x,y\in \widehat{\mathbb{R}}$, then 
\[
xy=(\sigma x)(\sigma y)+(\sigma x)dy+(\sigma y)dx=(\sigma x)(\sigma
y)+xdy+ydx. 
\]%
In particular, for each $\xi \in \mathbb{R}$: 
\[
\xi x=\xi (\sigma x)+\xi dx. 
\]%
\textbf{c)} If $m\in \mathbb{N},$ and $x\in \widehat{\mathbb{R}}$, then
(with $x^{0}=1$) 
\[
x^{m}=(\sigma x)^{m}+m(\sigma x)^{m-1}dx=(\sigma x)^{m}+mx^{m-1}dx. 
\]%
\textbf{d)} If $x\in \widehat{\mathbb{R}},$ and $x$ is not an infinitesimal,
then 
\[
\frac{1}{x}=\frac{1}{\sigma x}-\frac{1}{(\sigma x)^{2}}dx=\frac{1}{\sigma x}-%
\frac{1}{x^{2}}dx. 
\]%
\textbf{e)} If $x,y\in \widehat{\mathbb{R}},$ and $x$ is not an
infinitesimal, then 
\[
\frac{y}{x}=\frac{\sigma y}{\sigma x}+\frac{\left( \sigma x\right) dy-\left(
\sigma y\right) dx}{(\sigma x)^{2}}=\frac{\sigma y}{\sigma x}+\frac{xdy-ydx}{%
x^{2}}. 
\]%
\textbf{f)} If $x\in \widehat{\mathbb{R}}^{+},$ $m\in \mathbb{N}$ and $m>1$,
then there is a unique $y\in \widehat{\mathbb{R}}^{+}$ such that 
\[
y^{m}=x. 
\]%
Such $y$ will be denoted by $\sqrt[m]{x}$, and we have: 
\[
\sqrt[m]{x}=\sqrt[m]{\sigma x}+\frac{1}{m\sqrt[m]{(\sigma x)^{m-1}}}dx=\sqrt[%
m]{\sigma x}+\frac{1}{m\sqrt[m]{x^{m-1}}}dx, 
\]%
where $\sqrt[m]{\sigma x}$ and $\sqrt[m]{(\sigma x)^{m-1}}$ are the usual 
\textit{positive mth roots }of $\sigma x$ and $(\sigma x)^{m-1}$,
re\nolinebreak spectively.

\bigskip

\textbf{Proof} \ Only the proof of \textbf{f)} has some difficulty.\newline
If $x,y\in \widehat{\mathbb{R}}^{+}$, then $\sigma x>$ 0 and $\sigma y>0$.%
\newline
So, using \textbf{c)} and \textbf{proposition 3.3}, we have: 
\[
\begin{array}{ll}
y^{m}=x & \Leftrightarrow (\sigma y+dy)^{m}=\sigma x+dx\Leftrightarrow
(\sigma y)^{m}+m(\sigma y)^{m-1}dy=\sigma x+dx\Leftrightarrow \\ 
& \Leftrightarrow \left\{ 
\begin{array}{l}
\sigma y=\sqrt[m]{\sigma x} \\ 
dy=\frac{1}{m\sqrt[m]{(\sigma x)^{m-1}}}dx%
\end{array}\Leftrightarrow y=\sqrt[m]{\sigma x}+\frac{1}{m\sqrt[m]{(\sigma x)^{m-1}}}dx%
\right. .%
\end{array}%
\]%
But $\sqrt[m]{\sigma x}>0$, since $\sigma x>$ $0$; so 
\[
\sqrt[m]{\sigma x}+\frac{1}{m\sqrt[m]{(\sigma x)^{m-1}}}dx>0. 
\]%
We have proven the existence (and uniqueness) of $ \sqrt[m]{x} $ and the identity 
\[
\sqrt[m]{x}=\sqrt[m]{\sigma x}+\frac{1}{m\sqrt[m]{(\sigma x)^{m-1}}}dx. 
\]%
In particular, if $x\in {\mathbb{R}}^{+}$, then 
\[
\sqrt[m]{x}=\sqrt[m]{\sigma x} 
\]%
Using \textbf{c)} and the result already proved (notice that $x^{m-1}>0$,
since $\sigma \left( x^{m-1}\right) =$\newline
$={(\sigma x)}^{m-1}>0$), we obtain: 
\[
\sqrt[m]{x^{m-1}}=\sqrt[m]{({\sigma x)}^{m-1}}+\hat{\varepsilon}, 
\]%
where $\hat{\varepsilon}$ is the infinitesimal defined by 
\[
\hat{\varepsilon}:=\frac{1}{m\sqrt[m]{(\sigma x)^{{(m-1)}^{2}}}}(m-1)({%
\sigma x)}^{m-2}dx. 
\]%
Then, using \textbf{d)}, 
\[
\frac{1}{\sqrt[m]{x^{m-1}}}=\frac{1}{\sqrt[m]{({\sigma x)}^{m-1}}}-\frac{1}{%
\left( \sqrt[m]{({\sigma x)}^{m-1}}\right) ^{2}}\hat{\varepsilon}. 
\]%
Since the product of infinitesimals is $0$, we have: 
\[
\frac{1}{\sqrt[m]{x^{m-1}}}dx=\frac{1}{\sqrt[m]{({\sigma x)}^{m-1}}}dx.\text{
}\blacksquare 
\]

\bigskip

As an immediate consequence of the previous lemma, we obtain, using \textbf{%
proposition 3.3}, \textit{the basic algebraic rules of differentiation},
without using any notion of \textit{limit} in $\widehat{\mathbb{R}}$:

\bigskip

\textbf{Proposition 3.6 \ a)} If $x,y\in \widehat{\mathbb{R}}$, then 
\[
d(x+y)=dx+dy,\vspace{-3pt} 
\]%
\[
d(x-y)=dx-dy. 
\]%
\textbf{b)} If $x,y\in \widehat{\mathbb{R}}$, then 
\[
d(xy)=(\sigma x)dy+(\sigma y)dx=xdy+ydx. 
\]%
In particular, for each $\xi \in \mathbb{R}$: 
\[
d(\xi x)=\xi dx. 
\]%
\textbf{c)} If $m\in \mathbb{N},$ and $x\in \widehat{\mathbb{R}}$, then 
\[
d(x^{m})=m(\sigma x)^{m-1}dx=mx^{m-1}dx. 
\]%
\textbf{d)} If $x\in \widehat{\mathbb{R}},$ and \textit{x} is not an
infinitesimal, then 
\[
d\left( \frac{1}{x}\right) =-\frac{1}{(\sigma x)^{2}}dx=-\frac{1}{x^{2}}dx. 
\]%
\textbf{e)} If $x,y\in \widehat{\mathbb{R}},$ and \textit{x} is not an
infinitesimal, then 
\[
d\left( \frac{y}{x}\right) =\frac{\left( \sigma x\right) dy-\left( \sigma
y\right) dx}{{(\sigma x)}^{2}}=\frac{xdy-ydx}{x^{2}}. 
\]%
\textbf{f)} If $x\in \widehat{\mathbb{R}}^{+},$ $m\in \mathbb{N}$ and $m>1$,
then 
\[
d(\sqrt[m]{x})=\frac{1}{m\sqrt[m]{{(\sigma x)}^{m-1}}}dx=\frac{1}{m\sqrt[m]{%
x^{m-1}}}dx. 
\]

\bigskip

We close this section with a density theorem, and a theorem relating the 
\textit{generalized real continuum},$\left( \widehat{\mathbb{R}},<,+,\cdot
,0,1\right)$, to the \textit{Cantor-Dedekind continuum}.

\bigskip

\textbf{Theorem 3.7 (The Density Theorem)\vspace{-3pt}}

\textbf{a)} If $x$ and \textit{y} are generalized real numbers such that $%
x<y $, then there exists $ \zeta \in \mathbb{R}$ such that $x<\zeta
<y$.\newline
\textbf{b)} If $\xi $ and $\eta $ are real numbers such that $\xi <\eta $,
then there exists $z\in\widehat{\mathbb{R}}\,\backslash\, \mathbb{R}$ such that $\xi <z<\eta $.

\bigskip

\textbf{Proof \ a)} We may choose $\varsigma =\frac{\sigma x\text{ }+\text{ }%
\sigma y}{2}$.

\textbf{b)} If $\hat{\varepsilon}$ is an infinitesimal and $\hat{\varepsilon}%
\not=0$, then we may choose $z=\frac{\xi \text{ }+\text{ }\eta }{2}+\hat{%
\varepsilon}.$ $\blacksquare $

\bigskip

We already mentioned the trivial facts that $\approx $ is an equivalence
relation on $\widehat{\mathbb{R}},$ and the equivalence class of each $x\in 
\widehat{\mathbb{R}}$ is $m_{\approx }(x)=x+m_{\approx }(0)$. On the \textit{%
quotient }of $\widehat{\mathbb{R}}$ by $\approx ,$ i.e. the set $\widehat{%
%TCIMACRO{\U{211d} }%
%BeginExpansion
\mathbb{R}
%EndExpansion
}$/$\approx $ $:=\left\{ m_{\approx }(x)|x\in \widehat{\mathbb{R}}\right\} ,$
we consider now two binary operations, denoted by $\boxplus $ and $\boxdot $%
, and called \textit{addition} and \textit{multiplication}, respectively,
and a binary relation denoted by $\sqsubset $. These operations and relation
are defined by: 
\[
m_{\approx }\hspace{-1pt}\left( x\right) \boxplus m_{\approx }\hspace{-1pt}%
\left( y\right) :=m_{\approx }\hspace{-1pt}\left( x+y\right) , 
\]%
\[
m_{\approx }\hspace{-1pt}\left( x\right) \boxdot m_{\approx }\hspace{-1pt}%
\left( y\right) :=m_{\approx }\hspace{-1pt}\left( xy\right) , 
\]%
\[
m_{\approx }\hspace{-1pt}\left( x\right) \sqsubset m_{\approx }\hspace{-1pt}%
\left( y\right) :\Leftrightarrow x<y, 
\]%
using, at the right-hand of the previous identities, the obvious binary
operations and relation on $\widehat{\mathbb{R}}$.

\bigskip

It is a simple task to show that $\boxplus $, $\boxdot $, $\sqsubset $ are
well-defined, and to prove the next theorem.

\bigskip

\textbf{Theorem 3.8 \ a)} $\left( \widehat{%
%TCIMACRO{\U{211d} }%
%BeginExpansion
\mathbb{R}
%EndExpansion
}\text{/\hspace{-2pt}}\approx ,\sqsubset ,\boxplus ,\boxdot ,m_{\approx }%
\hspace{-1pt}\hspace{-1pt}\left( 0\right) ,m_{\approx }(1)\right) $ is an
ordered field with $m_{\approx }\hspace{-1pt}\hspace{-1pt}\left( 0\right) $
and $m_{\approx }(1)$ as\textit{\ zero }and\textit{\ identity elements},
respectively.\newline
\textbf{b)} The mapping $\phi :\widehat{%
%TCIMACRO{\U{211d} }%
%BeginExpansion
\mathbb{R}
%EndExpansion
}$/$\approx $ $\rightarrow \mathbb{R}$, defined by $\phi (m_{\approx }%
\hspace{-1pt}\hspace{-1pt}\left( x\right) ):=$ $\sigma x$ , is an ordered
field isomorphism of $\left( \widehat{%
%TCIMACRO{\U{211d} }%
%BeginExpansion
\mathbb{R}
%EndExpansion
}\text{/\hspace{-2pt}}\approx ,\sqsubset ,\boxplus ,\boxdot ,m_{\approx }%
\hspace{-1pt}\hspace{-1pt}\left( 0\right) ,m_{\approx }(1)\right) $ onto the 
\textit{Cantor-Dedekind continuum}, $\left( \mathbb{R},<,+,\cdot ,0,1\right) 
$; so if we denote these fields simply by $\widehat{%
%TCIMACRO{\U{211d} }%
%BeginExpansion
\mathbb{R}
%EndExpansion
}$/$\approx $ and $\mathbb{R}$, we have: 
\[
\widehat{%
%TCIMACRO{\U{211d} }%
%BeginExpansion
\mathbb{R}
%EndExpansion
}\text{/\hspace{-2pt}}\approx \text{\ \hspace{-0.15cm}}\cong \text{%
\thinspace }\mathbb{R},
\]%
i.e. $\widehat{%
%TCIMACRO{\U{211d} }%
%BeginExpansion
\mathbb{R}
%EndExpansion
}$/$\approx $ is isomorphic to $\mathbb{R}$.

\bigskip

As we have just seen:

\bigskip

\textit{If we take the monads in the structure }$\widehat{\mathbb{R}}$%
\textit{\ for points, as we do in the structure }$\widehat{%
%TCIMACRO{\U{211d} }%
%BeginExpansion
\mathbb{R}
%EndExpansion
}$/$\approx $\textit{, then we obtain the Cantor-Dedekind continuum.
Otherwise, we have a richer continuum with indiscernibility and nilpotent
infinitesimals}.

\bigskip

\section{Monads and Shadows}

The next two propositions show that $\left\{ m_{\approx }(x)|x\in \widehat{%
\mathbb{R}}\right\} $ is a partition of $\widehat{\mathbb{R}}$ into
infinite- -dimensional real affine spaces, each one with the same cardinal
as $\widehat{\mathbb{R}}$, and this is also true for $\left\{ m_{\approx
}(\xi )|\xi \in {\mathbb{R}}\right\} $ (since $m_{\approx }\hspace{-2pt}%
\left( x\right) =m_{\approx }\hspace{-2pt}\left( \sigma x\right) $, for each 
$x\in \widehat{\mathbb{R}}$).

\bigskip

\textbf{Proposition 4.1} \ The monad of each generalized real number has the
same cardinal as $\widehat{\mathbb{R}}$.

\bigskip

\textbf{Proof} \ Since $m_{\approx }\hspace{-2pt}\left( x\right) =m_{\approx
}\hspace{-2pt}\left( \sigma x\right) $, for each $x\in \widehat{\mathbb{R}}$%
, we may prove the proposition only for the monads of real numbers.\smallskip

Let $\xi \in \mathbb{R}$, and let $\widehat{\mathbb{R}}_{\xi }$ be the set
of all generalized real numbers $x=(\xi _{n})$ such that $\xi _{n}=\xi $,
for $n>1$. Then (denoting by $\left\vert \widehat{A}\right\vert $ the
cardinal of each subset $\widehat{A}$ of $\widehat{\mathbb{R}})$: 
\[
\left\vert \widehat{\mathbb{R}}_{\xi }\right\vert \leq \left\vert m_{\approx
}\hspace{-2pt}\left( \xi \right) \right\vert \leq \left\vert \mathbb{R}^{%
\mathbb{N}}\right\vert , 
\]%
where ${\mathbb{R}}^{{\mathbb{N}}}$ denotes the set of all sequences in $%
\mathbb{R}$.

Obviously, 
\[
\left\vert \widehat{\mathbb{R}}_{\xi }\right\vert =\left\vert \mathbb{R}%
\right\vert =2^{\aleph _{0}}, 
\]%
and 
\[
\left\vert {\mathbb{R}}^{{\mathbb{N}}}\right\vert =(2^{\aleph _{0}})^{\aleph
_{0}}=2^{\aleph _{0}\aleph _{0}}=2^{\aleph _{0}}. 
\]%
So 
\[
\left\vert m_{\approx }\hspace{-2pt}\left( \xi \right) \right\vert
=2^{\aleph _{0}}. 
\]%
Finally, 
\[
\left\vert \widehat{\mathbb{R}}\right\vert =\left\vert \cup \left\{
m_{\approx }(\xi )|\xi \in {\mathbb{R}}\right\} \right\vert =2^{\aleph
_{0}}2^{\aleph _{0}}=2^{\aleph _{0}}.\text{ }\blacksquare 
\]%
\smallskip

\textbf{Proposition 4.2} \ \textbf{a)} $m_{\approx }(0)$ is an
infinite-dimensional real vector space, if we consider addition and
multiplication defined on $\widehat{\mathbb{R}}\times \widehat{\mathbb{R}}$,
as \textit{vector addition} and \textit{scalar} \textit{multiplication}
defined on $m_{\approx }(0)\times m_{\approx }(0)$ and $\mathbb{R}\times
m_{\approx }(0)$, respectively. Moreover, $ m_{\approx }(0) $ contains the {\itshape real spaces\/} $ l^{p} $, for each $ p \in [1,+ \infty[ $.\newline
\textbf{b)} If we consider $m_{\approx }(0)$ with the structure of real
vector space mentioned in \textbf{a)}, then 
\[
m_{\approx }(x)\text{ is an infinite-dimensional real affine space, for each 
}x\in \widehat{\mathbb{R}}. 
\]

\textbf{Proof \ a)} It is trivial to prove that $m_{\approx }(0)$ is a real
vector space, using \textbf{proposition 2.3 a)}, \textbf{c)}. Finally, if $ p\in[1,+ \infty[$ and $ x=(\xi_{n} )\in l^{p}$, then $ \sum_{n=1}^{+\infty}\mid\xi_{n}\mid ^{p}\   <+\infty $ and, consequently, $ x=( \xi_{n} )\in m_{\approx }(0)$.
\textbf{b)} follows from \textbf{a)}, since $m_{\approx
}(x)=x+m_{\approx }(0)$, for each $x\in \widehat{\mathbb{R}}.$ $\blacksquare 
$

\bigskip

The next definition generalizes the concepts of \textit{monad} and \textit{%
shadow} to any subset of $\widehat{\mathbb{R}}$.

\bigskip

\textbf{Definition 4.3} \ Let $\widehat{A}\subseteq \widehat{\mathbb{R}}$.%
\vspace{-2pt}

The \textit{monad} of $\widehat{A}$ and the \textit{shadow} of $\widehat{A}$%
, denoted by $m_{\approx }\hspace{-3pt}\left( \widehat{A}\right) $ and $%
\sigma \hspace{-3pt}\left( \widehat{A}\right) $, respectively, are defined
by: 
\[
m_{\approx }\hspace{-3pt}\left( \widehat{A}\right) :=\cup \left\{ m_{\approx
}(x)|x\in \widehat{A}\right\} , 
\]%
\[
{{\sigma \hspace{-3pt}\left( \widehat{A}\right) :=}}\cup {\left\{ \left\{
\sigma x\right\} |x\in \widehat{A}\right\} }. 
\]%
So 
\[
m_{\approx }\hspace{-3pt}\left( \widehat{A}\right) =\left\{ x\in \widehat{%
\mathbb{R}}|\left( \exists y\in \widehat{A}\right) x\approx y\right\} , 
\]%
\[
\sigma \hspace{-3pt}\left( \widehat{A}\right) =\left\{ \sigma x|x\in 
\widehat{A}\right\} . 
\]%
Clearly, we have, for each $x\in \widehat{\mathbb{R}}$ and $\widehat{A}%
\subseteq \widehat{\mathbb{R}}$, 
\[
m_{\approx }(\{x\})=m_{\approx }(x), 
\]%
\[
\sigma (\{x\})=\{\sigma x\}, 
\]%
\[
\widehat{A}\subseteq m_{\approx }\hspace{-3pt}\left( \widehat{A}\right) . 
\]

\bigskip

The next three propositions state some basic properties of \textit{monads}
and \textit{shadows}, and admit quite straightforward proofs.

\bigskip

\textbf{Proposition 4.4} \ Let $\widehat{A},\widehat{B}\subseteq \widehat{%
\mathbb{R}}$. Then:\vspace{-2pt}

\textbf{a)} $\widehat{A}\subseteq \widehat{B}\hspace{2pt}\Rightarrow
m_{\approx }\hspace{-3pt}\left( \widehat{A}\right) \subseteq m_{\approx }%
\hspace{-3pt}\left( \widehat{B}\right) \wedge \sigma \hspace{-3pt}\left( 
\widehat{A}\right) \subseteq \sigma \hspace{-3pt}\left( \widehat{B}\right) $.

\textbf{b)} $\widehat{A}\subseteq \mathbb{R}\Leftrightarrow \sigma \hspace{%
-3pt}\left( \widehat{A}\right) =\widehat{A}$.

\textbf{c)} $m_{\approx }\hspace{-3pt}\left( \sigma \hspace{-3pt}\left( 
\widehat{A}\right) \right) =m_{\approx }\hspace{-3pt}\left( \widehat{A}%
\right) \wedge \sigma \hspace{-3pt}\left( m_{\approx }\hspace{-3pt}\left( 
\widehat{A}\right) \right) =\sigma \hspace{-3pt}\left( \widehat{A}\right) $.

\textbf{d)} $\widehat{A},\widehat{B}\subseteq \mathbb{R}\hspace{2pt}%
\Rightarrow \left( m_{\approx }\hspace{-3pt}\left( \widehat{A}\right)
=m_{\approx }\hspace{-3pt}\left( \widehat{B}\right) \Leftrightarrow \widehat{%
A}=\widehat{B}\right) $.

\bigskip

The \textit{monad }and\textit{\ shadow operators} on subsets of $\widehat{%
\mathbb{R}}$ preserve the Boolean operations on sets, with some looseness in
the case of \textit{intersection} and \textit{complement} (this is the core
information expressed in the next two propositions).

\bigskip

\pagebreak \textbf{Proposition 4.5 \ a)} $m_{\approx }(\emptyset )=\emptyset $, $%
m_{\approx }\hspace{-3pt}\left( \widehat{\mathbb{R}}\right) =m_{\approx }(%
\mathbb{R})=\widehat{\mathbb{R}}$.\newline
Let $\widehat{A},\widehat{B}\subseteq \widehat{\mathbb{R}}$. Then:

\textbf{b)} $m_{\approx }\hspace{-3pt}\left( m_{\approx }\hspace{-3pt}\left( 
\widehat{A}\right) \right) =m_{\approx }\hspace{-3pt}\left( \widehat{A}%
\right) $,

\textbf{c)} $m_{\approx }\hspace{-3pt}\left( \widehat{A}\cup \widehat{B}%
\right) =m_{\approx }\hspace{-3pt}\left( \widehat{A}\right) \cup m_{\approx }%
\hspace{-3pt}\left( \widehat{B}\right) $,

\textbf{d)} $m_{\approx }\hspace{-3pt}\left( \widehat{A}\cap \widehat{B}%
\right) \subseteq m_{\approx }\hspace{-3pt}\left( \widehat{A}\right) \cap 
\hspace{1pt}m_{\approx }\hspace{-3pt}\left( \widehat{B}\right) ,$

$m_{\approx }\hspace{-3pt}\left( \widehat{A}\right) \hspace{-3pt}\backslash
m_{\approx }\hspace{-3pt}\left( \widehat{B}\right) \subseteq m_{\approx }%
\hspace{-3pt}\left( \widehat{A}\backslash \widehat{B}\right) $.\smallskip

If $\widehat{A},\widehat{B}\subseteq \mathbb{R}$, then 
\[
m_{\approx }\hspace{-3pt}\left( \widehat{A}\cap \widehat{B}\right)
=m_{\approx }\hspace{-3pt}\left( \widehat{A}\right) \cap \hspace{1pt}%
m_{\approx }\hspace{-3pt}\left( \widehat{B}\right) , 
\]%
\[
m_{\approx }\hspace{-3pt}\left( \widehat{A}\backslash \widehat{B}\right)
=m_{\approx }\hspace{-3pt}\left( \widehat{A}\right) \hspace{-3pt}\backslash
m_{\approx }\hspace{-3pt}\left( \widehat{B}\right) . 
\]%
Let $\widehat{{\mathbf{A}}}\subseteq P\hspace{-2pt}\left( \widehat{\mathbb{R}%
}\right) $ $\left( \text{i.e. }\widehat{{\mathbf{A}}}\text{ is a collection
of subsets of }\widehat{\mathbb{R}}\right) $. Then:

\textbf{e)} $m_{\approx }\hspace{-3pt}\left( \cup \left\{ \widehat{A}|%
\widehat{A}\in \widehat{{\mathbf{A}}}\right\} \right) =\cup \left\{
m_{\approx }\hspace{-3pt}\left( \widehat{A}\right) |\widehat{A}\in \widehat{{%
\mathbf{A}}}\right\} $,

\textbf{f)} $m_{\approx }\hspace{-3pt}\left( \cap \left\{ \widehat{A}|%
\widehat{A}\in \widehat{{\mathbf{A}}}\right\} \right) \subseteq \cap \left\{
m_{\approx }\hspace{-3pt}\left( \widehat{A}\right) |\widehat{A}\in \widehat{{%
\mathbf{A}}}\right\} $.

If $\widehat{{\mathbf{A}}}\subseteq P\hspace{1pt}(\mathbb{R})$ $\left( \text{%
i.e. }\widehat{{\mathbf{A}}}\text{ is a collection of subsets of }%
%TCIMACRO{\U{211d} }%
%BeginExpansion
\mathbb{R}
%EndExpansion
\right) $, then 
\[
m_{\approx }\hspace{-3pt}\left( \cap \left\{ \widehat{A}|\widehat{A}\in 
\widehat{{\mathbf{A}}}\right\} \right) =\cap \left\{ m_{\approx }\hspace{-3pt%
}\left( \widehat{A}\right) |\widehat{A}\in \widehat{{\mathbf{A}}}\right\} . 
\]

\bigskip \textbf{Proposition 4.6 \ a)} $\sigma (\emptyset )=\emptyset $, $%
\sigma \hspace{-2pt}\left( \widehat{\mathbb{R}}\right) =\sigma (\mathbb{R})=%
\mathbb{R}.$\newline
Let $\widehat{A},\widehat{B}\subseteq \widehat{\mathbb{R}}$. Then:\newline
\textbf{b)} $\sigma \hspace{-2pt}\left( \sigma \hspace{-2pt}\left( \widehat{A%
}\right) \right) =\sigma \hspace{-2pt}\left( \widehat{A}\right) $,\newline
\textbf{c)} $\sigma \hspace{-2pt}\left( \widehat{A}\cup \widehat{B}\right)
=\sigma \hspace{-2pt}\left( \widehat{A}\right) \cup \sigma \hspace{-2pt}%
\left( \widehat{B}\right) $,\newline
\textbf{d)} $\sigma \hspace{-2pt}\left( \widehat{A}\cap \widehat{B}\right)
\subseteq \sigma \hspace{-2pt}\left( \widehat{A}\right) \cap \sigma \hspace{%
-2pt}\left( \widehat{B}\right) $, \newline
$\vspace{-6pt}\sigma \hspace{-2pt}\left( \widehat{A}\right) \hspace{-2pt}%
\backslash \sigma \hspace{-2pt}\left( \widehat{B}\right) \subseteq \sigma 
\hspace{-2pt}\left( \widehat{A}\backslash \widehat{B}\right) .$\vspace{8pt}\newline
If $\widehat{A}$ and $\widehat{B}$ are monads of subsets of $\mathbb{R}$,
then 
\begin{eqnarray*}
\sigma \hspace{-2pt}\left( \widehat{A}\cap \widehat{B}\right)  &=&\sigma 
\hspace{-2pt}\left( \widehat{A}\right) \cap \sigma \hspace{-2pt}\left( 
\widehat{B}\right) , \\
\sigma \hspace{-2pt}\left( \widehat{A}\backslash \widehat{B}\right) 
&=&\sigma \hspace{-2pt}\left( \widehat{A}\right) \hspace{-2pt}\backslash
\sigma \hspace{-2pt}\left( \widehat{B}\right) .
\end{eqnarray*}%
Let $\widehat{{\mathbf{A}}}\subseteq P\left( \widehat{\mathbb{R}}\right) $ $%
\left( \text{i.e.}\ \widehat{{\mathbf{A}}}\text{ is a collection of subsets of 
}\widehat{\mathbb{R}}\right) $. Then:\vspace{-2pt}

\textbf{e)} $\sigma \hspace{-2pt}\left( \cup \left\{ \widehat{A}|\widehat{A}%
\in \widehat{{\mathbf{A}}}\right\} \right) =\cup \left\{ \sigma \hspace{-2pt}%
\left( \widehat{A}\right) |\widehat{A}\in \widehat{{\mathbf{A}}}\right\} $,%
\vspace{6pt}\newline
\textbf{f)} $\sigma \hspace{-2pt}\left( \cap \left\{ \widehat{A}|\widehat{A}%
\in \widehat{{\mathbf{A}}}\right\} \right) \subseteq \cap \left\{ \sigma 
\hspace{-2pt}\left( \widehat{A}\right) |\widehat{A}\in \widehat{{\mathbf{A}}}%
\right\} $.\vspace{4pt}\newline
If $\widehat{{\mathbf{A}}}$ is a collection of monads of subsets of $\mathbb{%
R}$, then 
\[
\sigma \hspace{-2pt}\left( \cap \left\{ \widehat{A}|\widehat{A}\in \widehat{{%
\mathbf{A}}}\right\} \right) =\cap \left\{ \sigma \hspace{-2pt}\left( 
\widehat{A}\right) |\widehat{A}\in \widehat{{\mathbf{A}}}\right\} .
\]

Using \textbf{proposition 4.4}, \textbf{proposition 4.5}, and \textbf{%
proposition 4.6}, we could prove that the \textit{monad }and \textit{shadow
operators} on subsets of $\widehat{\mathbb{R}}$ preserve the basic concepts
of topology, and the concept of $\sigma $-\textit{algebra}, which is
fundamental in \textit{Measure Theory}. This is clearly expressed in the
next two propositions.

\bigskip

\textbf{Proposition 4.7 \ a)} If $X\subseteq \mathbb{R},\mathbf{B}$ is a
base for a topology for $X$, and $\widehat{\mathbf{B}}:\mathbf{=}\left\{
m_{\approx }(A)|A\in \mathbf{B}\right\} $, then 
\[
\widehat{\mathbf{B}}\text{ is a base for a topology for }m_{\approx }(X). 
\]%
\textbf{b)} Let $\mathbf{T}$ be a topology for $\mathbb{R}$.\newline
If $\widehat{\mathbf{T}}:\mathbf{=}\left\{ m_{\approx }(A)|A\in {\mathbf{T}}%
\right\} $, then 
\[
\widehat{\mathbf{T}}\text{ is a topology for }\widehat{\mathbb{R}}. 
\]%
\textbf{c)} If $\mathbf{T}$ is a topology for $\mathbb{R}$, $\widehat{%
\mathbf{T}}:\mathbf{=}\left\{ m_{\approx }(A)|A\in {\mathbf{T}}\right\} ,$
and $X\subseteq \mathbb{R}$, then 
\[
{int}_{\widehat{\mathbf{T}}}m_{\approx }\hspace{-2pt}\left( X\right)
=m_{\approx }\hspace{-2pt}\left( {int}_{\mathbf{T}}X\right) , 
\]%
\[
{ext}_{\widehat{\mathbf{T}}}m_{\approx }\hspace{-2pt}\left( X\right)
=m_{\approx }\hspace{-2pt}\left( {ext}_{\mathbf{T}}X\right) , 
\]%
\[
{bd}_{\widehat{\mathbf{T}}}m_{\approx }\hspace{-2pt}\left( X\right)
=m_{\approx }\hspace{-2pt}\left( {bd}_{\mathbf{T}}X\right) , 
\]%
\[
{cl}_{\widehat{\mathbf{T}}}m_{\approx }\hspace{-2pt}\left( X\right)
=m_{\approx }\hspace{-2pt}\left( {cl}_{\mathbf{T}}X\right) , 
\]%
\[
m_{\approx }\hspace{-2pt}\left( X\right) \text{ is open for }\widehat{%
\mathbf{T}}\Leftrightarrow X\text{ is open for }\mathbf{T}, 
\]%
\[
m_{\approx }\hspace{-2pt}\left( X\right) \text{ is closed for }\widehat{%
\mathbf{T}}\Leftrightarrow X\text{ is closed for }\mathbf{T}, 
\]%
\[
m_{\approx }\hspace{-2pt}\left( X\right) \text{ is compact for }\widehat{%
\mathbf{T}}\Leftrightarrow X\text{ is compact for }\mathbf{T}\text{ }; 
\]%
where $int_{\widehat{\mathbf{T}}}$ , $int_{\mathbf{T}}$ , $ext_{\widehat{%
\mathbf{T}}}$ , $ext_{\mathbf{T}}$ , $bd_{\widehat{\mathbf{T}}}$ , $bd_{%
\mathbf{T}}$ , $cl_{\widehat{\mathbf{T}}}$ , $cl_{\mathbf{T}}$ are the 
\textit{interior, exterior, boundary} and \textit{closure} operators for the
topologies $\widehat{\mathbf{T}}$ and $\mathbf{T}$, respectively.\newline
\textbf{d)} If $\mathbf{T}$ is a topology for $\mathbb{R}$, $\widehat{%
\mathbf{T}}:\mathbf{=}\left\{ m_{\approx }(A)|A\in {\mathbf{T}}\right\} ,$ $%
X\subseteq \mathbb{R}$, and $Y\subseteq X,$ then 
\[
{\widehat{\mathbf{T}}}_{m_{\approx }\left( X\right) }{\mathbf{=}}\left\{
m_{\approx }\hspace{-2pt}\left( A\right) |A\in {\mathbf{T}}_{X}\right\} , 
\]%
\[
cl_{{\widehat{\mathbf{T}}}_{m_{\approx }\left( X\right) }}m_{\approx
}(Y)=m_{\approx }(cl_{{\mathbf{T}}_{X}}Y)\text{ }; 
\]%
where $\widehat{\mathbf{T}}_{m_{\approx }\left( X\right) }$, ${\mathbf{T}}%
_{X}$ are the relativizations of $\widehat{\mathbf{T}}$, $\mathbf{T}$ to $%
m_{\approx }\hspace{-2pt}\left( X\right) ,$ $X,$ respectively.\newline
\textbf{e) }If $\mathbf{T}$ is a topology for $\mathbb{R}$, $\widehat{%
\mathbf{T}}:\mathbf{=}\left\{ m_{\approx }(A)|A\in {\mathbf{T}}\right\} ,$
and $X\subseteq \mathbb{R}$, then%
\[
m_{\approx }\hspace{-2pt}\left( X\right) \text{ is connected for }\widehat{%
\mathbf{T}}\Leftrightarrow X\text{ is connected for }\mathbf{T}. 
\]%
\textbf{f)} Let $\mathbf{B}$ be a $\sigma $-\textit{algebra} of subsets of $%
\mathbb{R}$.\newline
If $\widehat{\mathbf{B}}:\mathbf{=}\left\{ m_{\approx }(A)|A\in {\mathbf{B}}%
\right\} ,$ then 
\[
\widehat{\mathbf{B}}\text{ is a }\sigma \text{-\textit{algebra} of subsets
of }\widehat{\mathbb{R}}. 
\]

\bigskip

\textbf{Proposition 4.8 \ a)} If $\widehat{X}\subseteq \widehat{\mathbb{R}}$%
, $\widehat{\mathbf{B}}$ is a base for a topology for $\widehat{X}$ and a
collection of monads of subsets of $\mathbb{R}$, and $\mathbf{B}:\mathbf{=}%
\left\{ \sigma \hspace{-3pt}\left( \widehat{A}\right) |\widehat{A}\in 
\widehat{\mathbf{B}}\right\} $, then 
\[
\mathbf{B}\text{ is a base for a topology for }\sigma \hspace{-3pt}\left( 
\widehat{X}\right) .
\]%
\textbf{b)} Let $\widehat{\mathbf{T}}$ be a topology for $\widehat{\mathbb{R}%
}$ and a collection of monads of subsets of $\mathbb{R}$.\newline
If $\mathbf{T}:\mathbf{=}\left\{ \sigma \hspace{-3pt}\left( \widehat{A}%
\right) |\widehat{A}\in \widehat{\mathbf{T}}\right\} ,$ then 
\[
\mathbf{T}\text{ is a topology for }\mathbb{R}.
\]%
\textbf{c)} If $\widehat{\mathbf{T}}$ is a topology for $\widehat{\mathbb{R}}
$ and a collection of monads of subsets of $\mathbb{R}$, $\widehat{X}$ is
the monad of a subset of $\mathbb{R}$, and $\mathbf{T}:\mathbf{=}\left\{
\sigma \hspace{-3pt}\left( \widehat{A}\right) |\widehat{A}\in \widehat{%
\mathbf{T}}\right\} $, then 
\[
{int}_{\mathbf{T}}\sigma \hspace{-3pt}\left( \widehat{X}\right) =\sigma 
\hspace{-3pt}\left( {int}_{\widehat{\mathbf{T}}}\widehat{X}\right) ,
\]%
\[
{ext}_{\mathbf{T}}\sigma \hspace{-3pt}\left( \widehat{X}\right) =\sigma 
\hspace{-3pt}\left( {ext}_{\widehat{\mathbf{T}}}\widehat{X}\right) ,
\]%
\[
{bd}_{\mathbf{T}}\sigma \hspace{-3pt}\left( \widehat{X}\right) =\sigma 
\hspace{-3pt}\left( {bd}_{\widehat{\mathbf{T}}}\widehat{X}\right) ,
\]%
\[
{cl}_{\mathbf{T}}\sigma \hspace{-3pt}\left( \widehat{X}\right) =\sigma 
\hspace{-3pt}\left( {cl}_{\widehat{\mathbf{T}}}\widehat{X}\right) ,
\]%
\[
\sigma \hspace{-3pt}\left( \widehat{X}\right) \text{ is open for }\mathbf{T}%
\Leftrightarrow \widehat{X}\text{ is open for }\widehat{\mathbf{T}},
\]%
\[
\sigma \hspace{-3pt}\left( \widehat{X}\right) \text{ is closed for }\mathbf{T%
}\Leftrightarrow \widehat{X}\text{ is closed for }\widehat{\mathbf{T}},
\]%
\[
\sigma \hspace{-3pt}\left( \widehat{X}\right) \text{ is compact for }\mathbf{%
T}\Leftrightarrow \widehat{X}\text{ is compact for }\widehat{\mathbf{T}}%
\text{ };
\]%
where $int_{\widehat{\mathbf{T}}}$ , $int_{\mathbf{T}}$ , $ext_{\widehat{%
\mathbf{T}}}$ , $ext_{\mathbf{T}}$ , $bd_{\widehat{\mathbf{T}}}$ , $bd_{%
\mathbf{T}}$ , $cl_{\widehat{\mathbf{T}}}$ , $cl_{\mathbf{T}}$ are the 
\textit{interior, exterior, boundary} and \textit{closure} operators for the
topologies $\widehat{\mathbf{T}}$ and ${\mathbf{T}}$, respectively.\newline
\textbf{d)} If $\hspace{1pt}\widehat{\mathbf{T}}$ \hspace{1pt}is\hspace{1pt}
a \hspace{1pt}topology \hspace{1pt}for $\hspace{1pt}\widehat{\mathbb{R}}$ 
\hspace{1pt}and \hspace{1pt}a \hspace{1pt}collection \hspace{1pt}of \hspace{%
1pt}monads \hspace{1pt}of \hspace{1pt}subsets \hspace{1pt}of $\hspace{1pt}%
\mathbb{R}$, $\widehat{Y}\subseteq \widehat{X}\subseteq \widehat{%
%TCIMACRO{\U{211d} }%
%BeginExpansion
\mathbb{R}
%EndExpansion
},\vspace{-0.2cm}$

$\widehat{X}$ and $\widehat{Y}$ are monads of subsets of $\widehat{\mathbb{R}}$, and $%
\mathbf{T}:\mathbf{=}\left\{ \sigma \hspace{-3pt}\left( \widehat{A}\right) |%
\widehat{A}\in \widehat{\mathbf{T}}\right\} $, then 
\[
{\mathbf{T}}_{\sigma \left( \widehat{X}\right) }{\mathbf{=}}\left\{ \sigma 
\hspace{-3pt}\left( \widehat{A}\right) |\widehat{A}\in {\widehat{\mathbf{T}}}%
_{\widehat{X}}\right\} ,
\]%
\[
{cl}_{\mathbf{T}_{\sigma \left( \widehat{X}\right) }}\sigma \hspace{-3pt}%
\left( \widehat{Y}\right) =\sigma \hspace{-3pt}\left( {cl}_{\widehat{\mathbf{%
T}}_{\widehat{X}}}\widehat{Y}\right) \text{ };
\]%
where ${\widehat{\mathbf{T}}}_{\widehat{X}}$, ${\mathbf{T}}_{\sigma \left( 
\widehat{X}\right) }$ are the relativizations of $\widehat{\mathbf{T}}$, $%
\mathbf{T}$ to $\widehat{X},$ $\sigma \hspace{-2pt}\left( \widehat{X}\right) 
$, respectively.\newline
\textbf{e) }If $\widehat{\mathbf{T}}$ is a topology for $\widehat{\mathbb{R}}
$ and a collection of monads of subsets of $\mathbb{R}$, $\widehat{X}$ is
the monad of a subset of $\mathbb{R}$, and $\mathbf{T}:\mathbf{=}\left\{
\sigma \hspace{-3pt}\left( \widehat{A}\right) |\widehat{A}\in \widehat{%
\mathbf{T}}\right\} $, then 
\[
\sigma \hspace{-3pt}\left( \widehat{X}\right) \text{ is connected for }%
\mathbf{T}\Leftrightarrow \widehat{X}\text{ is connected for }\widehat{%
\mathbf{T}}.\vspace{-7pt}
\]%
\newline
\textbf{f)} Let $\widehat{\mathbf{B}}$ be a $\sigma $-\textit{algebra} of
subsets of $\widehat{\mathbb{R}}$ and a collection of monads of subsets of $%
\mathbb{R}$.\vspace{3pt}\newline
If $\mathbf{B}:\mathbf{=}\left\{ \sigma \hspace{-3pt}\left( \widehat{A}%
\right) |\widehat{A}\in \widehat{\mathbf{B}}\right\} $, then 
\[
\mathbf{B}\text{ is a }\sigma \text{-\textit{algebra} of subsets of }%
%TCIMACRO{\U{211d} }%
%BeginExpansion
\mathbb{R}
%EndExpansion
.
\]

\section{The Derivative}

Throughout this section, we shall not use any concept of \textit{limit} in
the generalized real continuum $\left( \text{i.e. }\widehat{%
%TCIMACRO{\U{211d} }%
%BeginExpansion
\mathbb{R}
%EndExpansion
}\right) $, working instead, in an actual manner, with the concepts of 
\textit{indiscernibility, shadow, differential}, and \textit{monad}. The
concept of \textit{limit} is only used in the \textit{Cantor}-\textit{%
Dedekind continuum} $\left( \text{i.e. }\mathbb{R}\right) $.

\bigskip

The first important step is the introduction of the concept of \textit{%
indiscernible extension }of a function $\phi :X\rightarrow Y$, where $%
X,Y\subseteq \mathbb{R}$.

\bigskip

\textbf{Definition 5.1} \ Let $X,Y\subseteq \mathbb{R}$.\newline
If $\phi :X\rightarrow Y$ and $f:m_{\approx }(X)\rightarrow m_{\approx }(Y)$
are functions, then \textit{f} is said to be an \textit{indiscernible
extension }of $\phi $ iff 
\[
(\forall x\in m_{\approx }(X))(f(\sigma x)=\phi (\sigma x)\wedge f(x)\approx
\phi (\sigma x)). 
\]

\bigskip

\pagebreak Clearly:

\bigskip

\textbf{Proposition 5.2} \ Let $X,Y\subseteq \mathbb{R}$. \newline
If $\phi :X\rightarrow Y$, $\psi :X\rightarrow Y$, $f:m_{\approx
}(X)\rightarrow m_{\approx }(Y)$ are functions, and $f$ is an indiscernible
extension of $\phi $ and $\psi $, then 
\[
\phi =\psi . 
\]

\bigskip

Before introducing the concept of \textit{interval} in $\widehat{\mathbb{R}}$%
, we must define the analogue on $\widehat{\mathbb{R}}$ of the usual linear
ordering $\leq $ on $\mathbb{R}$.

\bigskip

\textbf{Definition 5.3} \ Let $x,y\in \widehat{\mathbb{R}}$.\newline
We say that \textit{x} \textit{is less than or indiscernible from} \textit{y}%
, and we denote it by $x\lesssim y$, iff $\sigma x\leq \sigma y$ (where in $%
\sigma x\leq \sigma y$ we consider the usual linear ordering $\leq $ on $%
\mathbb{R})$, and we say that \textit{x} \textit{is greater than or
indiscernible from} \textit{y}, and we denote it by $x\gtrsim y$, iff $%
y\lesssim x$ \nolinebreak .\vspace{3pt}\newline
$\widehat{\mathbb{R}}_{0}^{+}$ and $\widehat{\mathbb{R}}_{0}^{-}$ denote the
subsets of $\widehat{\mathbb{R}}$ defined by 
\[
{\widehat{\mathbb{R}}}_{0}^{+}:=\left\{ x\in \widehat{\mathbb{R}}|x\gtrsim
0\right\} ={\widehat{\mathbb{R}}}^{+}\cup m_{\approx }(0), 
\]%
\[
{\widehat{\mathbb{R}}}_{0}^{-}:=\left\{ x\in \widehat{\mathbb{R}}|x\lesssim
0\right\} ={\widehat{\mathbb{R}}}^{-}\cup m_{\approx }(0). 
\]

Clearly:

\bigskip

\textbf{Proposition 5.4 \ a)} If $x,y\in \widehat{\mathbb{R}}$, then\vspace{%
-3pt}\newline
\[
x\lesssim y\Leftrightarrow x<y\vee x\approx y. 
\]%
\textbf{b)} Let $x,y,z\in \widehat{\mathbb{R}}$. Then: 
\[
x\lesssim x, 
\]%
\[
x\lesssim y\wedge y\lesssim x\Rightarrow x\approx y, 
\]%
\[
x\lesssim y\wedge y\lesssim z\Rightarrow x\lesssim z, 
\]%
\[
x\lesssim y\vee y\lesssim x, 
\]%
\[
x\lesssim y\Rightarrow x+z\lesssim y+z, 
\]%
\[
x\lesssim y\wedge z\gtrsim 0\Rightarrow xz\lesssim yz. 
\]%
So if we adopt the version of the usual \textit{antisymmetry} expressed by the second formula above, then we may consider $\lesssim $ a linear
ordering on $\widehat{\mathbb{R}}$.\newline
\textbf{c)} If $\hat{\varepsilon}$ and $\hat{\delta}$ are infinitesimals,
then 
\[
\hat{\varepsilon}\lesssim \hat{\delta}\wedge \hat{\delta}\lesssim \hat{%
\varepsilon}. 
\]%
\textbf{d)} $\widehat{\mathbb{R}}_{0}^{+}$ and $\widehat{\mathbb{R}}_{0}^{-}$
are the sets of \textit{nonnegative} and \textit{nonpositive} generalized
real numbers, i.e. 
\[
{\widehat{\mathbb{R}}}_{0}^{+}=\widehat{\mathbb{R}}\backslash {\widehat{%
\mathbb{R}}}^{{-}}, 
\]%
and 
\[
{\widehat{\mathbb{R}}}_{0}^{-}=\widehat{\mathbb{R}}\backslash {\widehat{%
\mathbb{R}}}^{{+}}. 
\]%
Furthermore: 
\[
\widehat{\mathbb{R}}_{0}^{+}\cap {\widehat{\mathbb{R}}}_{0}^{-}=m_{\approx
}(0). 
\]

\smallskip The next definition introduce concepts that are adaptations to $%
\lesssim $ (and $\gtrsim )$, on $\widehat{\mathbb{R}}$, of the usual notions
for $\leq $ (and $\geq )$, on $\mathbb{R}$.

\bigskip

\textbf{Definition 5.5} \ Let $\widehat{A}\subseteq \widehat{\mathbb{R}},$
and $L,l\in \widehat{\mathbb{R}}$. Then:\smallskip \newline
$L$ is a $\lesssim $-\textit{upper bound} of $\widehat{A}$ iff 
\[
\left( \forall x\in \widehat{A}\right) x\lesssim L. 
\]%
$l$ is a $\lesssim $-\textit{lower bound} of $\widehat{A}$ iff 
\[
\left( \forall x\in \widehat{A}\right) x\gtrsim l 
\]%
$\widehat{A}$ is $\lesssim $-\textit{bounded above} iff $\widehat{A}$ has a $%
\lesssim $-upper bound, and $\widehat{A}$ is $\lesssim $-\textit{bounded
below} iff $\widehat{{A}}$ has a $\lesssim $-lower bound.

$\widehat{A}$ is $\lesssim $-\textit{bounded} iff $\widehat{A}$ is $\lesssim 
$-bounded above and $\lesssim $-bounded below.

$\widehat{A}$ is $\lesssim $-\textit{unbounded} iff $\widehat{A}$ is not $%
\lesssim $-bounded.

$L$ is a $\lesssim $-\textit{maximum} of $\widehat{A}$ iff $L$ $\in \widehat{%
A}$ and $L$ is a $\lesssim $-upper bound of $\widehat{A}$.

$l$ is a $\lesssim $-\textit{minimum} of $\widehat{A}$ iff $l$ $\in \widehat{%
A}$ and $l$ is a $\lesssim $-lower bound of $\widehat{A}$.

$L$ is a $\lesssim $-\textit{supremum} of $\widehat{A}$ iff $L$ is a $%
\lesssim $-minimum of $\lesssim $-\textit{Up}$\left( \widehat{A}\right) $,
where $\lesssim $-\textit{Up}$\left( \widehat{A}\right) $ is the set of all $%
\lesssim $-upper bounds of $\widehat{A}$.

$l$ is a $\lesssim $-\textit{infimum} of $\widehat{A}$ iff $l$ is a $\lesssim $-\textit{maximum} of $\lesssim $-\textit{Lo}$\left( \widehat{A}\right) $, where $%
\lesssim $-\textit{Lo}$\left( \widehat{A}\right) $ is the set of all $%
\lesssim $-lower bounds of $\widehat{A}$.

\bigskip

\pagebreak \textbf{Proposition 5.6} \ Let $\widehat{A}$ $\subseteq \widehat{\mathbb{R}}%
, $ and $L,L%
%TCIMACRO{\U{b4}}%
%BeginExpansion
{\acute{}}%
%EndExpansion
,l,l^{\prime }$ $\in \widehat{\mathbb{R}}$.\smallskip \newline
\textbf{a)} If $L$ $\approx $ $L^{\prime }$ and $l$ $\approx $ $l^{\prime }$%
, then

\begin{center}
$L$ is a $\lesssim $-upper bound of $\widehat{A}$ iff $L^{\prime }$ is a $%
\lesssim $-upper bound of $\widehat{A}$,
\end{center}

and

\begin{center}
$l$ is a $\lesssim $-lower bound of $\widehat{A}$ iff $l^{\prime }$ is a $%
\lesssim $-lower bound of $\widehat{A}$.
\end{center}

\textbf{b)} $\lesssim $-\textit{Up}$\left( \widehat{A}\right) $and $\lesssim 
$-\textit{Lo}$\left( \widehat{A}\right) $ are monads of subsets of $\mathbb{R%
}$.\smallskip \newline
\textbf{c)} If $L$ is a $\lesssim $-maximum of $\widehat{A}$, then\smallskip

\begin{center}
$L^{\prime }$ is a $\lesssim $-maximum of $\widehat{A}\Rightarrow L^{\prime
}\approx $ $L$.
\end{center}

Similarly, if $l$ is a $\lesssim $-minimum of $\widehat{A}$, then

\begin{center}
$l^{\prime }$ is a $\lesssim $-minimum of $\widehat{A}\Rightarrow l^{\prime
}\approx $ $l$.
\end{center}

If $L$ is a $\lesssim $-maximum of $\widehat{A},$ and $\widehat{A}$ is the
monad of a subset of $\mathbb{R}$, then

\begin{center}
$L^{\prime }$ $\approx L\Rightarrow L^{\prime }$ is a $\lesssim $-maximum of 
$\widehat{A}$.
\end{center}

Similarly, if $l$ is a $\lesssim $-minimum of $\widehat{A},$ and $\widehat{A}
$ is the monad of a subset of $\mathbb{R}$, then

\begin{center}
$l^{\prime }$ $\approx l\Rightarrow l^{\prime }$ is a $\lesssim $-minimum of 
$\widehat{A}$.
\end{center}
\textbf{d)} If $L$ is a $\lesssim $-supremum of $\widehat{A}$, then

\begin{center}
$L^{\prime }$ is a $\lesssim $-supremum of $\widehat{A}\Leftrightarrow $ $%
L^{\prime }$ $\approx $ $L$.
\end{center}

Similarly, if $l$ is a $\lesssim $-infimum of $\widehat{A}$, then

\begin{center}
$l^{\prime }$ is a $\lesssim $-infimum of $\widehat{A}\Leftrightarrow $ $%
l^{\prime }$ $\approx $ $l$.

\bigskip
\end{center}

\textbf{Proof} \ \textbf{a)} is trivial, since $\sigma L=\sigma L^{\prime }$
and $\sigma l=\sigma l^{\prime }$.\vspace{-3pt}

\textbf{b)} Using \textbf{a)}, we have:

\begin{center}
$m_{\approx }\hspace{-2pt}\left( \lesssim \text{-}Up\left( \widehat{A}%
\right) \right) =\hspace{1pt}\lesssim $-\textit{Up}$\left( \widehat{A}%
\right) $.
\end{center}

Then, using \textbf{proposition 4.4 c)}:

\begin{center}
$\lesssim $-\textit{Up}$\left( \widehat{A}\right) =m_{\approx }\hspace{-2pt}%
\left( \sigma \hspace{-2pt}\left( \lesssim \text{-}Up\left( \widehat{A}%
\right) \right) \right) $.
\end{center}

Similarly, for $\lesssim $-\textit{Lo}$\left( \widehat{A}\right) $.\newline
\textbf{c)} Let $L$ be a $\lesssim $-maximum of $\widehat{A}$.\newline
If $L^{\prime }$ is a $\lesssim $-maximum of $\widehat{A}$, then, since $%
L,L^{\prime }$ $\in \widehat{A}$, 
\[
L\lesssim L^{\prime }\wedge L^{\prime }\lesssim L. 
\]%
So, by \textbf{proposition 5.4 b)}, 
\[
L\approx L^{\prime }. 
\]%
Let $\widehat{A}$ be the monad of a subset of $\mathbb{R}$.\newline
If $L^{\prime }$ $\approx $ $L$, then, by \textbf{a)},

\begin{center}
$L^{\prime }$ is a $\lesssim $-upper bound of $\widehat{A}$.
\end{center}

On the other hand, since $L$ $\in \widehat{A}$, $L^{\prime }$ $\approx $ $L$%
, and $\widehat{A}$ is the monad of a subset of $\mathbb{R}$, we have: 
\[
L^{\prime }\in \hspace{1pt}\widehat{A}. 
\]%
So

\begin{center}
$L^{\prime }$ is a$\lesssim $-maximum of $\widehat{A}$.
\end{center}

Similarly, for the concept of $\lesssim $-minimum.\newline
\textbf{d)} follows directly from \textbf{b)} and \textbf{c)}. $\blacksquare 
$

\bigskip

We have just seen that the concepts of $\lesssim $-upper bound and $\lesssim 
$-lower bound are invariant under indiscernibility, and so are the concepts of $\lesssim $-supremum and $%
\lesssim $-infimum.

\bigskip

\textbf{Corollary 5.7} \ Let $\widehat{A}\subseteq \widehat{\mathbb{R}}$,
and $L,l$ $\in \widehat{\mathbb{R}}$.\newline
\textbf{a)} If $L$ is a $\lesssim $-supremum of $\widehat{A}$, then $\sigma
L $ is also a $\lesssim $-supremum of $\widehat{A}$, and each $\lesssim $-su-%
\newline
premum of $\widehat{A}$ has $\sigma L$ as its shadow.\newline
When $l$ is a $\lesssim $-infimum of $\widehat{A}$, $\sigma l$ is also a $%
\lesssim $-infimum of $\widehat{A}$, and each $\lesssim $-infimum of $%
\widehat{A}$ has $\sigma l$ as its shadow.\newline
\textbf{b)} If $L$ is a $\lesssim $-maximum of $\widehat{A}$ and $\sigma
L\in \widehat{A}$, then $\sigma L$ is a $\lesssim $-maximum of $\widehat{A}$%
, and each $\lesssim $-maximum of $\widehat{A}$ has $\sigma L$ as its shadow.%
\newline
When $l$ is a $\lesssim $-minimum of $\widehat{A}$ and $\sigma l\in \widehat{%
A}$, then $\sigma l$ is a $\lesssim $-minimum of $\widehat{A}$, and each $%
\lesssim $-minimum of $\widehat{A}$ has $\sigma l$ as its shadow.

\bigskip

\textbf{Proof} \ \textbf{a)} and \textbf{b)} follow immediately from \textbf{%
proposition 5.6 d)}, and \textbf{proposition 5.6} \textbf{a)}, \textbf{c)},
respectively. $\blacksquare $

\bigskip

\textbf{Definition 5.8} \ Let $\widehat{A}\subseteq \widehat{\mathbb{R}},$
and $L,l$ $\in \widehat{\mathbb{R}}$.\newline
If $L$ is a $\lesssim $-supremum of $\widehat{A}$, then $\sigma L$ is called
the \textit{real supremum} of $\widehat{A}$.\newline
Similarly, if $l$ is a $\lesssim $-infimum of $\widehat{A}$, then $\sigma l$
is called the \textit{real infimum} of $\widehat{A}$.\newline
If $L$ is a $\lesssim $-maximum of $\widehat{A}$ and $\sigma L\in \widehat{A}
$, then $\sigma L$ is said to be the \textit{real maximum} of $\widehat{A}$.%
\newline
In a similar manner, if $l$ is a $\lesssim $-minimum of $\widehat{A}$ and $%
\sigma l\in \widehat{A}$, then $\sigma l$ is said to be the \textit{real
minimum} of $\widehat{A}$.\newline
We denote the \textit{real supremum}, the \textit{real infimum}, the \textit{%
real maximum}, and the \textit{real minimum} of $\widehat{A}$ by $\sup_{r}%
\widehat{A}$, $\inf_{r}\widehat{A}$, $\max_{r}\widehat{A}$, and $\min_{r}%
\widehat{A}$, respectively.

\bigskip

Before presenting a \textit{Completeness Property} for $\widehat{\mathbb{R}}$%
, we need the following lemma:

\bigskip

\textbf{Lemma 5.9} \ Let $\widehat{A}\subseteq \widehat{\mathbb{R}},$ and $%
L,l$ $\in \widehat{\mathbb{R}}$.\newline
\textbf{a\textit{)}} $L$ is a $\lesssim $-upper bound of $\widehat{A}$ iff $%
\sigma L$ is an upper bound of $\sigma \hspace{-3pt}\left( \widehat{A}%
\right) $.\newline
$l$ is a $\lesssim $-lower bound of $\widehat{A}$ iff $\sigma l$ is a lower
bound of $\sigma \hspace{-3pt}\left( \widehat{A}\right) $.\newline
\textbf{b) }$\sigma \hspace{-1pt}\left( \lesssim \text{-\textit{Up}}\left( 
\widehat{A}\right) \hspace{-1pt}\right) \hspace{-1pt}\hspace{-2pt}=$\textit{%
Up}$\left( \sigma \hspace{-1pt}\hspace{-1pt}\left( \widehat{A}\right) 
\hspace{-1pt}\right) \hspace{-1pt},$ and $\sigma \hspace{-1pt}\left(
\lesssim \text{-\textit{Lo}}\left( \widehat{A}\right) \hspace{-1pt}\right) 
\hspace{-1pt}\hspace{-1pt}\hspace{-1pt}=$\textit{Lo}$\left( \sigma \hspace{%
-1pt}\hspace{-1pt}\left( \widehat{A}\right) \hspace{-1pt}\right) \hspace{-1pt%
};$ where \textit{Up}$\left( \sigma (\widehat{A})\hspace{-1pt}\right) $ and 
\textit{Lo}$\left( \sigma (\widehat{A})\hspace{-1pt}\right) $ are the sets
of all \textit{upper bounds} and \textit{lower bounds} of $\sigma (\widehat{A%
})$, \hspace{1pt}\hspace{1pt}respectively, for the usual linear ordering $%
\leq $ on $\mathbb{R}$.\newline
\textbf{c)} $L$ is a $\lesssim $-maximum of $\widehat{A}\Rightarrow \sigma
L=\max \sigma \hspace{-3pt}\left( \widehat{A}\right) .$\newline
$l$ is a $\lesssim $-minimum of $\widehat{A}\Rightarrow \sigma l=\min \sigma 
\hspace{-3pt}\left( \widehat{A}\right) .$\newline
If $\widehat{A}$ is the monad of a subset of $\mathbb{R}$, then 
\[
\sigma L=\max \sigma \hspace{-3pt}\left( \widehat{A}\right) \Rightarrow L%
\text{ is a }\lesssim \text{-maximum of }\widehat{A}. 
\]%
Similarly, if $\widehat{A}$ is the monad of a subset of $\mathbb{R}$, then 
\[
\sigma l=\min \sigma \hspace{-3pt}\left( \widehat{A}\right) \Rightarrow l%
\text{ is a }\lesssim \text{-minimum of }\widehat{A}. 
\]%
\textbf{d)} $L$ is a $\lesssim $-supremum of $\widehat{A}\Leftrightarrow
\sigma L=\sup \sigma \hspace{-3pt}\left( \widehat{A}\right) .$\newline
$l$ is a $\lesssim $-infimum of $\widehat{A}\Leftrightarrow \sigma l=\inf
\sigma \hspace{-3pt}\left( \widehat{A}\right) .$

\bigskip

\textbf{Proof} \ \textbf{a)} Clearly:%
\[
L\text{ is a}\lesssim \text{-upper bound of }\widehat{A}\Leftrightarrow
\left( \forall x\in \widehat{A}\right) \sigma x\leq \sigma L\Leftrightarrow 
\]%
\[
\Leftrightarrow \left( \forall \xi \in \sigma \hspace{-3pt}\left( \widehat{A}%
\right) \right) \xi \leq \sigma L\Leftrightarrow \sigma L\text{ is an upper
bound of }\sigma \hspace{-3pt}\left( \widehat{A}\right) . 
\]

We may use a similar proof for the notion of $\lesssim $-lower bound.\newline
\textbf{b)} For each $x\in \widehat{\mathbb{R}}$, we have, using \textbf{a)}%
, and \textbf{proposition} \textbf{4.4 b)}, \textbf{c)}:%
\[
x\in \text{\hspace{0.6pt}\hspace{-1pt}}\lesssim \text{-\textit{Up}}\left( 
\widehat{A}\right) \Leftrightarrow \sigma x\in \text{\textit{Up}}\left(
\sigma \hspace{-3pt}\left( \widehat{A}\right) \right) \Leftrightarrow x\in
m_{\approx }\hspace{-3pt}\left( \text{\textit{Up}}\left( \sigma \hspace{-3pt}%
\left( \widehat{A}\right) \right) \right) . 
\]%
So 
\[
\lesssim \text{-\textit{Up}}\left( \widehat{A}\right) =m_{\approx }\hspace{%
-3pt}\left( \text{\textit{Up}}\left( \sigma \hspace{-3pt}\left( \widehat{A}%
\right) \right) \right) . 
\]%
Then, using \textbf{proposition} \textbf{4.4 b)}, \textbf{c)}, 
\[
\sigma \hspace{-3pt}\left( \lesssim \text{-\textit{Up}}\left( \widehat{A}%
\right) \right) =\sigma \hspace{-3pt}\left( m_{\approx }\hspace{-3pt}\left( 
\text{\textit{Up}}\left( \sigma \hspace{-3pt}\left( \widehat{A}\right)
\right) \right) \right) =\text{\textit{Up}}\left( \sigma \hspace{-3pt}\left( 
\widehat{A}\right) \right) . 
\]%
Similarly, for $\lesssim $-\textit{Lo}$\left( \widehat{A}\right) $.\newline
\textbf{c)} If $L$ is $\lesssim $-maximum of $\widehat{A}$, then

\begin{center}
$L$ is a $\lesssim $-upper bound of $\widehat{A}$,
\end{center}

and so, by \textbf{a)},

\begin{center}
$\sigma L$ is an upper bound of $\sigma \hspace{-3pt}\left( \widehat{A}%
\right) $.
\end{center}

On the other hand, we have, since $L$ $\in \widehat{A}:$%
\[
\sigma L\in \sigma \hspace{-3pt}\left( \widehat{A}\right) . 
\]%
So 
\[
\sigma L=\max \sigma \hspace{-3pt}\left( \widehat{A}\right) . 
\]%
Let $\widehat{A}$ be the monad of a subset of $\mathbb{R}$.\vspace{3pt}%
\newline
If $\sigma L=\max \sigma \hspace{-3pt}\left( \widehat{A}\right) $, then $%
\sigma L$ is an upper bound of $\sigma \hspace{-3pt}\left( \widehat{A}%
\right) $, and so, by \textbf{a)}, $L$ is a $\lesssim $-upper bound of $%
\widehat{A}$.\newline
On the other hand, since $L\approx \sigma L$ and $\sigma L\in \sigma \hspace{%
-3pt}\left( \widehat{A}\right) $, 
\[
L\in m_{\approx }\hspace{-3pt}\left( \sigma \hspace{-3pt}\left( \widehat{A}%
\right) \right) . 
\]%
But $m_{\approx }\hspace{-3pt}\left( \sigma \hspace{-3pt}\left( \widehat{A}%
\right) \right) =m_{\approx }\hspace{-3pt}\left( \widehat{A}\right) $ (by 
\textbf{proposition 4.4 c)}), and $m_{\approx }\hspace{-3pt}\left( \widehat{A%
}\right) =\widehat{A}$ (by \textbf{proposition 4.5 b)}).\newline
So 
\[
L\in \widehat{A}. 
\]%
We have just proven that

\begin{center}
$L$ is a $\lesssim $-maximum of $\widehat{A}$.
\end{center}

Similarly, for the notion of $\lesssim $-minimum.\newline
\textbf{d)} Using \textbf{b)}, \textbf{c)}, and \textbf{proposition 5.6 b)},
we have:%
\[
L\text{ is a}\lesssim \text{-supremum of }\widehat{A}\Leftrightarrow L\text{
is a}\lesssim \text{-minimum of}\lesssim \text{-\textit{Up}}\left( \widehat{A%
}\right) \Leftrightarrow 
\]%
\[
\Leftrightarrow \sigma L=\min \sigma \hspace{-2pt}\left( \lesssim \text{-%
\textit{Up}}\left( \widehat{A}\right) \right) \Leftrightarrow \sigma L=\min 
\text{\textit{Up}}\left( \sigma \hspace{-3pt}\left( \widehat{A}\right)
\right) \Leftrightarrow \sigma L=\sup \sigma \hspace{-3pt}\left( \widehat{A}%
\right) .
\]%
Similarly, for the notion of $\lesssim $-infimum. $\blacksquare $

\bigskip

\textbf{Theorem 5.10 (The Completeness Property of} $\widehat{\mathbb{R}})$%
\newline
Let $\widehat{A}$ be a nonempty subset of $\widehat{\mathbb{R}}$.\newline
\textbf{a)} If $\widehat{A}$ is $\lesssim $-bounded above, then there exists 
$\sup_{r}\widehat{A}$.\newline
\textbf{b)} If $\widehat{A}$ is $\lesssim $-bounded below, then there exists 
$\inf_{r}\widehat{A}$.

\bigskip

\textbf{Proof \ a)} If $\widehat{A}$ is $\lesssim $-bounded above, then 
\[
\lesssim \text{-\textit{Up}}\left( \widehat{A}\right) \not=\emptyset . 
\]%
So 
\[
\sigma \hspace{-3pt}\left( \lesssim \text{-\textit{Up}}\left( \widehat{A}%
\right) \right) \not=\emptyset . 
\]%
Then, by \textbf{lemma 5.9} \textbf{b)}, 
\[
\text{\textit{Up}}\left( \sigma \hspace{-3pt}\left( \widehat{A}\right)
\right) \not=\emptyset . 
\]%
Since $\sigma \hspace{-3pt}\left( \widehat{A}\right) \not=\emptyset $
(because $\widehat{A}\not=\emptyset )$, we infer, using the \textbf{%
Completeness Property} of $\mathbb{R}$, that there exists $\sup \sigma 
\hspace{-3pt}\left( \widehat{A}\right) $.\newline
Denoting $\sup \sigma \hspace{-3pt}\left( \widehat{A}\right) $ by $L$, we
have, using \textbf{lemma 5.9 d), }and the fact that $L\in 
%TCIMACRO{\U{211d} }%
%BeginExpansion
\mathbb{R}
%EndExpansion
$:%
\[
L=\sup \sigma \hspace{-3pt}\left( \widehat{A}\right) \Leftrightarrow L\text{
is a }\lesssim \text{-supremum of }\widehat{A}\Rightarrow L=\text{sup}_{r}%
\text{ }\widehat{A}. 
\]

\textbf{b)} admits a similar proof. $\blacksquare $

\bigskip

\textbf{Definition 5.11} \ Let $\alpha ,\alpha _{1},\beta ,\beta _{1}\in 
%TCIMACRO{\U{211d} }%
%BeginExpansion
\mathbb{R}
%EndExpansion
$; with $\alpha \leq \beta $.\vspace{-5pt}

The \textit{closed, open}, and \textit{half-open intervals} determined by
the ordered pair $\left( \alpha ,\beta \right) $, de-\smallskip \newline
noted by $\widehat{[\alpha ,\beta ]},\widehat{{]}\alpha {,}\beta \lbrack },%
\widehat{]\alpha ,\beta ]}$ and $\widehat{[\alpha ,\beta \lbrack },$
respectively, are defined by: 
\[
\widehat{\lbrack \alpha ,\beta ]}:=\left\{ x\in \widehat{\mathbb{R}}|\alpha
\lesssim x\lesssim \beta \right\} , 
\]%
\[
\widehat{]\alpha ,\beta \lbrack }:=\left\{ x\in \widehat{\mathbb{R}}|\alpha
<x<\beta \right\} , 
\]%
\[
\widehat{]\alpha ,\beta ]}:=\left\{ x\in \widehat{\mathbb{R}}|\alpha
<x\lesssim \beta \right\} , 
\]%
\[
\widehat{\lbrack \alpha ,\beta \lbrack }:=\left\{ x\in \widehat{\mathbb{R}}%
|\alpha \lesssim x<\beta \right\} ~. 
\]%
The intervals just introduced are $\lesssim $-bounded sets.

We use the symbols $-\infty $ and $+\infty $ to introduce the intervals that
are $\lesssim $-unbounded sets: 
\[
\widehat{\lbrack \alpha _{1},+\infty \lbrack }:=\left\{ x\in \widehat{%
\mathbb{R}}|\alpha _{1}\lesssim x\right\} , 
\]%
\[
\widehat{]\alpha _{1},+\infty \lbrack }:=\left\{ x\in \widehat{\mathbb{R}}%
|\alpha _{1}<x\right\} , 
\]%
\[
\widehat{]-\infty ,\beta _{1}]}:=\left\{ x\in \widehat{\mathbb{R}}|x\lesssim
\beta _{1}\right\} , 
\]%
\[
\widehat{]-\infty ,\beta _{1}[}:=\left\{ x\in \widehat{\mathbb{R}}|x\text{ }{%
<}\text{ }\beta _{1}\right\} , 
\]%
\[
\widehat{]-\infty ,+\infty \lbrack }:=\widehat{\mathbb{R}}. 
\]

\vspace{2pt}The next proposition admits a quite straightforward proof (in
particular, \textbf{e)} follows easily from \textbf{proposition 4.4 b)}, 
\textbf{c)}, \textbf{proposition 4.7 e)}, \textbf{proposition 5.12 a)}, and
the well-known fact that \textit{the connected subsets of} $\mathbb{R}$, 
\textit{for the usual topology, are the intervals}).

\bigskip

\textbf{Proposition 5.12} \ \textbf{a)} The intervals in $\widehat{\mathbb{R}%
}$ are the monads of the correspondent intervals in $\mathbb{R}$, and the
intervals in $\mathbb{R}$ are the shadows of the correspondent intervals in $%
\widehat{\mathbb{R}}$; for example, if $\alpha ,\alpha _{1},\beta $ $\in 
\mathbb{R},$ and $\alpha \leq \beta $, then 
\[
\widehat{\lbrack \alpha ,\beta ]}=m_{\approx }([\alpha ,\beta ]), 
\]%
\[
\lbrack \alpha ,\beta ]=\sigma \hspace{-2pt}\left( \widehat{[\alpha ,\beta ]}%
\right) , 
\]%
\[
\widehat{\left[ \alpha _{1},+\infty \right[ }=m_{\approx }\hspace{-2pt}%
\left( \left[ \alpha _{1},+\infty \right[ \right) , 
\]%
\[
\left[ \alpha _{1},+\infty \right[ =\sigma \hspace{-2pt}\left( \widehat{%
\left[ \alpha _{1},+\infty \right[ }\right) . 
\]%
\textbf{b)} Let $\alpha ,\beta \in \mathbb{R},$ with $\alpha \leq \beta .$
Then:%
\[
\widehat{\lbrack \alpha ,\beta ]}\neq \emptyset , 
\]%
\[
\widehat{]\alpha ,\beta \lbrack }=\emptyset \Leftrightarrow \alpha =\beta , 
\]%
\[
\widehat{]\alpha ,\beta ]}=\emptyset \Leftrightarrow \alpha =\beta , 
\]%
\[
\widehat{\lbrack \alpha ,\beta \lbrack }=\emptyset \Leftrightarrow \alpha
=\beta . 
\]%
\textbf{c) }Let $\alpha ,\alpha ^{\prime },\alpha _{1},\alpha _{1}^{\prime
},\beta ,\beta ^{\prime },\beta _{1},\beta _{1}^{\prime }\in \mathbb{R},$
with $\alpha \leq \beta $ and $\alpha ^{\prime }\leq \beta ^{\prime }.$ Then:%
\[
\widehat{\lbrack \alpha ,\beta ]}=\widehat{[\alpha ^{\prime },\beta ^{\prime
}]}\Rightarrow \left( \alpha ,\beta \right) =\left( \alpha ^{\prime },\beta
^{\prime }\right) , 
\]%
\[
\alpha \text{\hspace{-2pt}}<\hspace{-2pt}\beta \text{\hspace{1pt}}\wedge
\left( \widehat{]\alpha ,\beta \lbrack }=\widehat{]\alpha ^{\prime },\beta
^{\prime }[}\vee \widehat{]\alpha ,\beta ]}=\widehat{]\alpha ^{\prime
},\beta ^{\prime }]}\vee \widehat{[\alpha ,\beta \lbrack }=\widehat{[\alpha
^{\prime },\beta ^{\prime }[}\right) \Rightarrow \left( \alpha ,\beta
\right) =\left( \alpha ^{\prime },\beta ^{\prime }\right) \text{\hspace{-1pt}%
}, 
\]%
\[
\alpha \text{\hspace{-2pt}}=\hspace{-2pt}\beta \text{\hspace{1pt}}\wedge
\left( \widehat{]\alpha ,\beta \lbrack }=\widehat{]\alpha ^{\prime },\beta
^{\prime }[}\vee \widehat{]\alpha ,\beta ]}=\widehat{]\alpha ^{\prime
},\beta ^{\prime }]}\vee \widehat{[\alpha ,\beta \lbrack }=\widehat{[\alpha
^{\prime },\beta ^{\prime }[}\right) \Rightarrow \alpha ^{\prime }=\beta
^{\prime }\text{\hspace{-1pt}}, 
\]%
\[
\widehat{\lbrack \alpha _{1},+\infty \lbrack }=\widehat{[\alpha _{1}^{\prime
},+\infty \lbrack }\vee \widehat{]\alpha _{1},+\infty \lbrack }=\widehat{%
]\alpha _{1}^{\prime },+\infty \lbrack }\Rightarrow \alpha _{1}=\alpha
_{1}^{\prime }, 
\]%
\[
\widehat{]-\infty ,\beta _{1}]}=\widehat{]-\infty ,\beta _{1}^{\prime }]}%
\vee \widehat{]-\infty ,\beta _{1}[}=\widehat{]-\infty ,\beta _{1}^{\prime }[%
}\Rightarrow \beta _{1}=\beta _{1}^{\prime }. 
\]%
Intervals of different kind are never identical, unless they are both the
empty set; for example (still with $\alpha ,\alpha _{1},\beta \in \mathbb{R}%
, $ and $\alpha \leq \beta $),%
\[
\widehat{\lbrack \alpha ,\beta ]}\neq \widehat{[\alpha ,\beta \lbrack }, 
\]%
\[
\widehat{]\alpha ,\beta \lbrack }\neq \widehat{[\alpha _{1},+\infty \lbrack }%
, 
\]%
\[
\widehat{]\alpha _{1},+\infty \lbrack }\neq \widehat{]-\infty ,\beta _{1}]}, 
\]%
\[
\widehat{\lbrack \alpha ,\beta \lbrack }=\widehat{]\alpha ,\beta \lbrack }%
\Leftrightarrow \widehat{[\alpha ,\beta \lbrack }=\widehat{]\alpha ,\beta
\lbrack }=\emptyset . 
\]

\textbf{d)} If $\widehat{I}$ is an interval in $\widehat{\mathbb{R}}$, then%
\[
(\exists ^{1}I)\left( I\text{ is an interval in }\mathbb{R}\wedge \widehat{I}%
=m_{\approx }(I)\right) , 
\]%
\[
m_{\approx }(\widehat{I})=\widehat{I}, 
\]%
\[
\sigma (\widehat{I})\subseteq \widehat{I}. 
\]%
\textbf{e)} Let $\mathbf{T}$ be the usual topology for $\mathbb{R}$, and let 
$\widehat{\mathbf{T}}:\mathbf{=}\left\{ m_{\approx }(A)|A\in {\mathbf{T}}%
\right\} $.

If $\widehat{X}$ is the monad of a subset of $\mathbb{R}$, then%
\[
\widehat{X}\text{ is connected for }\widehat{\mathbf{T}}\Leftrightarrow 
\widehat{X}\text{ is an interval in }\widehat{\mathbb{R}}. 
\]

\bigskip Now we may introduce the concept of \textit{length }\ of a $%
\lesssim $-bounded interval in $\widehat{%
%TCIMACRO{\U{211d} }%
%BeginExpansion
\mathbb{R}
%EndExpansion
}$ (notice how \textbf{proposition 5.12 b)}, \textbf{c) }is relevant to the
next definition)$.$

\bigskip \textbf{Definition 5.13 \ }Let $\alpha ,\beta \in 
%TCIMACRO{\U{211d} }%
%BeginExpansion
\mathbb{R}
%EndExpansion
,$ and $\alpha \leq \beta .$ If $\widehat{I}$ is one of the intervals $%
\widehat{[\alpha ,\beta ]},\widehat{]\alpha ,\beta \lbrack },\widehat{%
]\alpha ,\beta ]},\widehat{[\alpha ,\beta \lbrack },$ then the \textit{%
length }of $\widehat{I},$ denoted by $l\hspace{-2pt}\left( \widehat{I}%
\right) $, is defined by: 
\[
l\hspace{-2pt}\left( \widehat{I}\right) :=\beta -\alpha . 
\]

\bigskip Clearly:

\bigskip \textbf{Proposition 5.14 \ }If $\alpha \in \mathbb{R}$, then%
\[
l\hspace{-2pt}\left( \widehat{{]}\alpha ,\alpha \lbrack }\right) =l\hspace{%
-2pt}\left( \widehat{{]}\alpha ,\alpha {]}}\right) =l\hspace{-2pt}\left( 
\widehat{[\alpha ,\alpha \lbrack }\right) =l\hspace{-2pt}\left( \widehat{%
[\alpha ,\alpha {]}}\right) =0, 
\]%
but 
\[
\widehat{{]}\alpha ,\alpha \lbrack }=\widehat{{]}\alpha ,\alpha {]}}=%
\widehat{[\alpha ,\alpha \lbrack }=\emptyset , 
\]%
and 
\[
\widehat{\left[ \alpha ,\alpha \right] }=m_{\approx }(\alpha ). 
\]

\bigskip

\textbf{Remark 5.15} \ The intervals in $\widehat{\mathbb{R}}$ have no
clear-cut (i.e. pointlike) extremities.\vspace{-3pt}

For example, if $\alpha ,\beta ,\gamma \in \mathbb{R}$ and $\alpha
<\beta <\gamma $, then $\widehat{[\alpha ,\beta ]},\widehat{[\beta ,\gamma ]}
$ have $m_{\approx }(\alpha ),m_{\approx }(\beta )$ and $m_{\approx }(\beta
),m_{\approx }(\gamma )$ as extremities, respectively, and 
\[
\widehat{\lbrack \alpha ,\beta ]}\cap \widehat{[\beta ,\gamma ]}=m_{\approx
}(\beta ). 
\]%
The intervals in $\widehat{\mathbb{R}}$ are particularly fit to devise a
model for the \textit{flux of Time}:

A \textit{stretch of Time} is an interval $\widehat{[\alpha ,\beta ]}$ $%
(\alpha ,\beta \in \mathbb{R};$ $\alpha <\beta $ ) whose members will be
called \textit{instants}.

Each \textit{now} is the intersection of two adjacent stretches of Time,
such as%
\[
\widehat{\lbrack \alpha ,\beta ]},\widehat{[\beta ,\gamma ]}\text{ }(\alpha
,\beta ,\gamma \in \mathbb{R};\text{ }\alpha <\beta <\gamma ). 
\]%
So each \textit{now} is the monad of an instant, and consequently, a set of 
\textit{indiscernible in\nolinebreak stants} with the power of the continuum
and length $0$, since, for each $\beta \in \widehat{\mathbb{R}}$ ,%
\[
\left\vert m_{\approx }(\beta )\right\vert =2^{\aleph _{0}}\wedge
l(m_{\approx }(\beta ))=l(\widehat{\left[ \beta ,\beta \right] }=0. 
\]%
Also, being the intersection of two adjacent intervals, each \textit{now}
has a dual \textit{past}-\textit{future }nature.

This conception of Time is reminiscent of the ideas of the \textit{Stoic
philosophers} (especially Chrysippos) $\left[ 12\right] $.

\pagebreak 

We now present the concept of \textit{differentiability}.

\bigskip

\textbf{Definition 5.16} \ Let \textit{I} be an open interval in $\mathbb{R}$%
, let $\xi _{0}\in I$, and let $\phi :I\rightarrow \mathbb{R}$ be a function.%
\vspace{-0.25cm}

If $f:m_{\approx }(I)\rightarrow ~\widehat{\mathbb{R}}~$ is an indiscernible
extension of $\phi $, then \textit{f} is said to be \textit{differentiable }%
at $\xi _{0}$ iff there exists a real number $\alpha $ such that%
\[
\left( \forall x\in m_{\approx }(\xi _{0})\right) f(x)=\phi (\xi
_{0})+\alpha dx, 
\]

with the proviso that $\alpha :=\lim_{\xi \rightarrow \xi _{0}}\frac{\phi
\left( \xi \right) -\phi (\xi _{0})}{\xi -\xi _{0}}$, when such limit exists
in $\mathbb{R}$ (considering the usual definition of \textit{limit}).\vspace{%
-0.25cm}

If $J$ is an open subinterval (in $\mathbb{R})$ of $I$, then \textit{f} is
said to be \textit{differentiable }on $m_{\approx }(J)$ iff $f$ is
differentiable at each $\xi _{0}\in $ $J$.

\bigskip

\textbf{Proposition 5.17} \ Let $I$ be an open interval in $\mathbb{R}$, let 
$\xi _{0}\in $ \textit{I}, and let $\phi :I\rightarrow \mathbb{R}$ be a
function.\vspace{-0.25cm}

If $f:m_{\approx }(I)\rightarrow ~\widehat{\mathbb{R}}~$ is an indiscernible
extension of $\phi $, $\alpha $ and $\beta $ are real numbers, and%
\[
\left( \forall x\in m_{\approx }(\xi _{0})\right) \left( f(x)=\phi (\xi
_{0})+\alpha dx\wedge f(x)=\phi (\xi _{0})+\beta dx\right) , 
\]%
then 
\[
\alpha =\beta . 
\]

\bigskip

\textbf{Proof \ }If we choose $x\hspace{-2pt}\in \hspace{-2pt}m_{\approx
}(\xi _{0})$ such that \textit{dx} is the eventually null sequence $\left(
1,0,0,0,\dots \hspace{-1pt}\right) $\hspace{-1pt}\nolinebreak , then the
conclusion follows at once from $\alpha dx=\beta dx$, since%
\[
\alpha dx=\beta dx\Leftrightarrow (\alpha ,0,0,0,\dots )=(\beta ,0,0,0,\dots
).\text{ }\blacksquare 
\]

\textbf{Definition 5.18 \ }With the notation and the conditions of \textbf{%
definition 5.16}, if $f$ is differentiable at $\xi _{0}$, then $\alpha $ is
called the \textit{derivative }of\textit{\ }$f$\textit{\ }at\textit{\ }$x$,
for each $x\in m_{\approx }(\xi _{0})$, and we denote it by $f^{\prime }(x)$.

\bigskip

\textbf{Remark 5.19 \ }Let $I$ be an open interval in $\mathbb{R}$, let $\xi
_{0}\in $ $I$, and let $f:m_{\approx }(I)\rightarrow \widehat{\mathbb{R}}$
be an indiscernible extension of $\phi :I\rightarrow \mathbb{R}$ .\vspace{%
-0.25cm}

If $f$ is differentiable at $\xi _{0}$, then $f^{\prime }(x)$ exists (in $%
\mathbb{R})$, for each $x\hspace{-1pt}\in \hspace{-1pt}m_{\approx }(\xi
_{0}) $, and $f^{\prime }(x)\hspace{-1pt}=\hspace{-1pt}\nolinebreak
f^{\prime }(\xi _{0})$.\smallskip\ But the differentiability of $f$ at $\xi
_{0}$ does not entail the existence of $\lim_{\xi \rightarrow \xi _{0}}\frac{%
\phi \left( \xi \right) -\phi (\xi _{0})}{\xi -\xi _{0}}$ \nolinebreak ,
although if this is the case, then $f^{\prime }(\xi _{0})$ coincides with
this limit, by the proviso of \textbf{definition 5.16}.

As an example, let us consider the functions $\phi :\mathbb{R}\rightarrow 
\mathbb{R}$ and $f:\widehat{\mathbb{R}}\rightarrow \widehat{\mathbb{R}}$,
defined by $\phi (\xi ):=\left\vert \xi \right\vert $ and $f(x):=\left\{ 
\begin{array}{c}
x,\text{ if \ }x>0 \\ 
0,\text{ if \ }x\in m_{\approx }(0) \\ 
-x,\text{ if \ }x<0%
\end{array}%
\right. $ , where $|\hspace{5pt}|$ denote the usual absolute value in $%
\mathbb{R}$. Clearly, $f$ is an indiscernible extension of $\phi $,
differentiable at $\xi _{0}=0$ with $f^{\prime }(\xi _{0})=0$, but $%
\lim_{\xi \rightarrow \xi _{0}}\frac{\phi \left( \xi \right) -\phi (\xi _{0})%
}{\xi -\xi _{0}}$ does not exist in $\mathbb{R}$.

\bigskip

\textbf{Proposition 5.20 \ }Let $I$ be an open interval in $\mathbb{R}$, let 
$\xi _{0}\in $ $I$, and let $\phi :I\rightarrow \mathbb{R}$ be a function.%
\vspace{-0.25cm}

\textit{If }$f:m_{\approx }(I)\rightarrow \widehat{\mathbb{R}}$ is an
indiscernible extension of $\phi $, and $f$ is \textit{differentiable} at $%
\xi _{0}$, then%
\[
\left( \forall x\in m_{\approx }(\xi _{0})\right) f(x)=f(\xi _{0})+f^{\prime
}(\xi _{0})dx. 
\]

(Notice that we could have written%
\[
\left( \forall x\in m_{\approx }(\xi _{0})\right) f(x)=f(\sigma x)+f^{\prime
}(x)dx, 
\]

since $\sigma x=\xi _{0}$ and $f^{\prime }(x)=f^{\prime }(\xi _{0}),$ for
each $x\in m_{\approx }(\xi _{0})$).

\bigskip

\textbf{Proof \ }Just remember that $f(\xi _{0})=\phi (\xi _{0})$. $%
\blacksquare $

\bigskip

\textbf{Proposition 5.20} expresses, in analytic terms, the \textit{%
geometric idea} associated with the concept of \textit{differentiability}.
This idea was clearly expressed by G. W. Leibniz and G. de L'H\^{o}pital
(via Johann Bernoulli), and it is closely related to the use of nilpotent
infinitesimals, as the Dutch theologian and mathematician B. Nieuwentijt
first realized (around 1695):%
\[
\hspace{7pt}\text{\textit{If}}\mathit{\ }f\mathit{\ }\text{\textit{is
differentiable at }}\xi _{0}\text{, \textit{then the graph of} }f\mathit{\ }%
\text{\textit{coincides locally} }(\text{\textit{i.e. for infinitesimal}}
\]%
\textit{increments of the argument around} $\xi $$_{0})\mathit{\ }$\textit{%
with its tangent at the point} $(\xi _{0},$\thinspace $f(\xi _{0})).$

\bigskip

The next lemma is necessary to establish the basic algebraic rules of
derivation.

\bigskip

\textbf{Lemma 5.21 \ }Let $I$ be an open interval in $\mathbb{R}$, and let $%
f:m_{\approx }(I)\rightarrow \widehat{\mathbb{R}}$, $g:m_{\approx
}(I)\rightarrow \widehat{\mathbb{R}}$ be indiscernible extensions of $\phi
:I\rightarrow \mathbb{R}$, $\psi :I\rightarrow \mathbb{R}$, respectively.%
\vspace*{-0.25cm}

\textbf{a)} For fixed $\alpha ,\beta \in \mathbb{R}$, if $\phi (\xi
):=\alpha \xi +\beta $, then we may define $f$ by 
\[
f(x):=\alpha x+~\beta. 
\]
\textbf{b)} $f+g$\textit{, }$fg$ are indiscernible extensions of $\phi +\psi 
$, $\phi \psi $, respectively.

\textbf{c)} If $\psi \left( \xi \right) \neq 0$, for each $\xi \in I,$ then%
\[
\frac{f}{g}\text{ \ is an indiscernible extension of }\frac{\phi }{\psi }~. 
\]%
\textbf{d)} For fixed $m$ $\in \mathbb{N}$:%
\[
f^{m}\text{ is an indiscernible extension of }\phi ^{m}. 
\]

\vspace*{-0.25cm}If $\phi \left( \xi \right) >0,$ for each $\xi \in I,$ and $m>1$, then%
\[
\sqrt[m]{f}\text{ is an indiscernible extension of }\sqrt[m]{\phi }. 
\]

\vspace*{-0.25cm}\textbf{e)} Let $J$ be an open interval in $\mathbb{R}$
such that $\phi \left( I\right) \subseteq J$, and let $h:m_{\approx }\hspace{%
-2pt}\left( J\right) \rightarrow \widehat{\mathbb{R}}$ be an indiscernible
extension of $\theta :J\rightarrow {\mathbb{R}}$ \textit{.} Then:%
\[
h\circ f\text{ is an indiscernible extension of }\theta \circ \phi . 
\]

\vspace*{-0.25cm}\textbf{f)} If $f$ is injective and $m_{\approx }\hspace{%
-2pt}\left( \phi \left( I\right) \right) \subseteq f\hspace{-2pt}\left(
m_{\approx }\hspace{-2pt}\left( I\right) \right) $, then $\phi $ is also
injective and%
\[
f^{-1}\text{ is an indiscernible extension of }\phi ^{-1}. 
\]

\textbf{Proof \ }Only the proof of \textbf{e)} and \textbf{f)} has some
difficulty.\vspace{-3pt}

\textbf{e)} First, we shall prove that $h\circ f$ makes sense.\vspace{-3pt}

Let $z$ $\in m_{\approx }(I)$.\vspace{-3pt}

Then, since $\sigma z\in I$ (by \textbf{proposition 4.4 b)}, \textbf{c)})
and $\phi \left( I\right) \subseteq $ $J$, we have:%
\[
f\left( z\right) \approx \phi \left( \sigma {z}\right) \in J. 
\]%
So%
\[
f\left( z\right) \in m_{\approx }(J). 
\]%
We have proven that%
\[
f\left( m_{\approx }(I)\right) \subseteq m_{\approx }(J). 
\]%
Now let $x\in m_{\approx }(I)$.

Then 
\[
\left( h\circ f\right) \left( \sigma x\right) =h\left( f\left( \sigma
x\right) \right) =h\left( \phi \left( \sigma x\right) \right) =h(\sigma \phi
(\sigma x))=\theta \left( \sigma \phi \left( \sigma x\right) \right) =\theta
(\phi (\sigma x))=(\theta \circ \phi )(\sigma x).
\]%
On the other hand, since $\phi \hspace{-2pt}\left( \sigma x\right) =\sigma f%
\hspace{-2pt}\left( x\right) $, we have:%
\[
(h\circ f)(x)=h(f(x))\approx \theta (\phi (\sigma x))=(\theta \circ \phi
)(\sigma x).
\]%
We have proven that%
\[
h\circ f\text{ is an indiscernible extension of }\theta \circ \phi .
\]

\textbf{f)} If $f$ is injective, then so is $\phi $, since $\phi \left( \xi
\right) =f\left( \xi \right) $, for each $\xi \in I$ .\vspace{-5pt}

Let $z$ $\in m_{\approx }(I)$.\vspace{-5pt}

Then, since $\sigma f\hspace{-2pt}\left( z\right) \in f\hspace{-2pt}\left(
m_{\approx }(I)\right) $ (because $\sigma z\in I,$ by \textbf{proposition
4.4 b)}, \textbf{c)}, $I\subseteq m_{\approx }(I),$ and $\sigma f\hspace{-2pt%
}\left( z\right) =\phi (\sigma z)=f\hspace{-2pt}\left( \sigma z\right) $),
we have:%
\[
f^{-1}\hspace{-2pt}\left( \sigma f\left( z\right) \right) =f^{-1}\hspace{-2pt%
}\left( f(\sigma z)\right) =\sigma z=\phi ^{-1}\hspace{-2pt}\left( \phi
(\sigma z)\right) =\phi ^{-1}\hspace{-2pt}\left( \sigma f\left( z\right)
\right) ,
\]%
\[
f^{-1}\hspace{-2pt}\left( f\left( z\right) \right) =z\approx \sigma z=\phi
^{-1}\hspace{-2pt}\left( \phi \left( \sigma z\right) \right) ={\phi }^{-1}%
\hspace{-2pt}\left( \sigma f\left( z\right) \right) .
\]%
Since $f$ is an indiscernible extension of $\phi ,$ we have $f\hspace{-2pt}%
\left( m_{\approx }(I)\right) \subseteq \nolinebreak m_{\approx }(\phi
\left( I\right) ).$ So, from $m_{\approx }\left( \phi \left( I\right)
\right) \subseteq f\left( m_{\approx }\left( I\right) \right) ,$ we infer
that%
\[
f\left( m_{\approx }(I)\right) =m_{\approx }(\phi \left( I\right) ).
\]

\vspace*{-0.25cm}We have proven that\vspace*{-0.5cm}

\[
f^{-1}\text{ is an indiscernible extension of }\phi ^{-1}.\text{ }%
\blacksquare 
\]

\smallskip Let us state the \textit{basic algebraic properties of the
derivative}:

\bigskip

\textbf{Proposition 5.22} \ Let $I$ be an open interval in $\mathbb{R}$, let 
$f:m_{\approx }(I)\rightarrow \widehat{\mathbb{R}}$, $g:m_{\approx
}(I)\rightarrow \widehat{\mathbb{R}}$ be indiscernible extensions of $\phi
:I\rightarrow \mathbb{R},\psi :I\rightarrow \mathbb{R}$, respectively, and
let $\xi _{0}\in I$.\vspace{-6pt}

\textbf{a)} If $\alpha $ and $\beta $ are fixed real numbers, and $f$ is
defined by $f(x):=\alpha x+\beta $, then \textit{f} is differentiable at $%
\xi _{0}$, and%
\[
\left( \forall x\in m_{\approx }(\xi _{0})\right) f^{\prime }(x)=\alpha .
\]%
\textbf{b)} Let $f$ and $g$ be differentiable at $\xi _{0}$. Then:\vspace{%
-6pt}

If at least one of the limits $\lim_{\xi \rightarrow \xi _{0}}\frac{\phi
\left( \xi \right) -\phi (\xi _{0})}{\xi -\xi _{0}}$, $\lim_{\xi \rightarrow
\xi _{0}}\frac{\psi \left( \xi \right) -\psi (\xi _{0})}{\xi -\xi _{0}}$
exists in $\mathbb{R}$, then $f+g$ is differentiable at $\xi _{0}$, and for
each $x\in m_{\approx }(\xi _{0})$:%
\[
(f+g)^{\prime }(x)=f^{\prime }(x{)+}g^{\prime }(x). 
\]%
\textbf{c)} Let \textit{f} and \textit{g} be differentiable at $\xi _{0}$%
.\smallskip

\vspace{-6pt}If $\lim_{\xi \rightarrow \xi _{0}}\frac{\phi \left( \xi
\right) -\phi (\xi _{0})}{\xi -\xi _{0}}$ and $\lim_{\xi \rightarrow \xi
_{0}}\frac{\psi \left( \xi \right) -\psi (\xi _{0})}{\xi -\xi _{0}}$ exist
in $\mathbb{R}$, then \textit{fg} is differentiable at $\xi _{0}$, and we
have, for each $x\in m_{\approx }(\xi _{0})$:%
\[
(fg)^{\prime }(x)=f^{\prime }(x)g(\xi _{0})+g^{\prime }(x)f(\xi _{0}). 
\]

\vspace{-3pt}If $\phi \left( \xi _{0}\right) \neq 0,\lim_{\xi \rightarrow
\xi _{0}}\frac{\phi \left( \xi \right) -\phi (\xi _{0})}{\xi -\xi _{0}}$, $%
\lim_{\xi \rightarrow \xi _{0}}\psi \left( \xi \right) $ exist and $%
\lim_{\xi \rightarrow \xi _{0}}\frac{\psi \left( \xi \right) -\psi (\xi _{0})%
}{\xi -\xi _{0}}$ does not \smallskip \newline
exist in $\mathbb{R}$, then $fg$ is differentiable at $\xi _{0}$, and we
have, for each $x\in m_{\approx }(\xi _{0})$:%
\[
(fg)^{\prime }(x)=f^{\prime }(x)g(\xi _{0})+g^{\prime }(x)f(\xi _{0}). 
\]

If $\psi \left( \xi _{0}\right) \neq 0,\lim_{\xi \rightarrow \xi _{0}}\frac{%
\psi \left( \xi \right) -\psi (\xi _{0})}{\xi -\xi _{0}}$, $\lim_{\xi
\rightarrow \xi _{0}}\phi \left( \xi \right) $ exist and $\lim_{\xi
\rightarrow \xi _{0}}\frac{\phi \left( \xi \right) -\phi (\xi _{0})}{\xi
-\xi _{0}}$ does not \vspace{3pt}\newline
exist in $\mathbb{R}$, then $fg$ is differentiable at $\xi _{0}$, and we
have, for each $x\in m_{\approx }(\xi _{0})$:%
\[
(fg)^{\prime }(x)=f^{\prime }(x)g(\xi _{0})+g^{\prime }(x)f(\xi _{0}). 
\]%
\textbf{d)} Let $f$ and $g$ be differentiable at $\xi _{0}$, and let $\psi
\left( \xi \right) \neq 0$, for each $\xi \in I$.\vspace{-4pt}

If $\lim_{\xi \rightarrow \xi _{0}}\frac{\phi \left( \xi \right) -\phi (\xi
_{0})}{\xi -\xi _{0}}$ and $\lim_{\xi \rightarrow \xi _{0}}\frac{\psi \left(
\xi \right) -\psi (\xi _{0})}{\xi -\xi _{0}}$ exist in $\mathbb{R}$, then $\frac{f}{g}
$ is differentiable at $\xi _{0}$, \vspace{2pt}\newline
and we have, for each $x\in m_{\approx }(\xi _{0})$:%
\[
\left( \frac{f}{g}\right) ^{\prime }(x)=\frac{f^{\prime }(x)g(\xi _{0}{)-}%
g^{\prime }(x)f(\xi _{0})}{{g(\xi _{0})}^{{2}}}.
\]%
If $\phi \left( \xi _{0}\right) \neq 0$, $\lim_{\xi \rightarrow \xi _{0}}%
\frac{\phi \left( \xi \right) -\phi (\xi _{0})}{\xi -\xi _{0}}$, $\lim_{\xi
\rightarrow \xi _{0}}~\frac{1}{\psi (\xi )}$ exist and $\lim_{\xi
\rightarrow \xi _{0}}\frac{\frac{1}{\psi (\xi )}-\frac{1}{\psi (\xi _{0})}}{%
\xi -\xi _{0}}$ does not \vspace{4pt}\newline
exist in $\mathbb{R}$, then $\frac{f}{g}$ is differentiable at $\xi _{0}$, and we
have, for each $x\in m_{\approx }(\xi _{0})$:%
\[
\left( \frac{f}{g}\right) ^{\prime }(x)=\frac{f^{\prime }(x)g(\xi _{0}{)-}%
g^{\prime }(x)f(\xi _{0})}{{g(\xi _{0})}^{{2}}}.
\]%
If $\lim_{\xi \rightarrow \xi _{0}}\frac{\frac{1}{\psi (\xi )}-\frac{1}{\psi
(\xi _{0})}}{\xi -\xi _{0}}$, $\lim_{\xi \rightarrow \xi _{0}}\phi \left(
\xi \right) $ exist and $\lim_{\xi \rightarrow \xi _{0}}\frac{\phi \left(
\xi \right) -\phi (\xi _{0})}{\xi -\xi _{0}}$ does not \vspace{4pt}
exist in $\mathbb{R}$, then $\frac{f}{g}$ is differentiable at $\xi _{0}$, and for
each $x\in m_{\approx }(\xi _{0})$:%
\[
\left( \frac{f}{g}\right) ^{\prime }(x)=\frac{f^{\prime }(x)g(\xi _{0}{)-}%
g^{\prime }(x)f(\xi _{0})}{{g(\xi _{0})}^{{2}}}.
\]%
\textbf{e)} Let $m\in \mathbb{N}$, and let $f$ be differentiable at $\xi
_{0}.\vspace{-3pt}$

If $\lim_{\xi \rightarrow \xi _{0}}\frac{\phi \left( \xi \right) -\phi (\xi
_{0})}{\xi -\xi _{0}}$ exists in $\mathbb{R}$, then $f^{m}$ is
differentiable at $\xi _{0}$, and for each $x\in \nolinebreak m_{\approx
}(\xi _{0})$:%
\[
(f^{m})^{\prime }(x)=mf(\xi _{0})^{m-1}f^{\prime }(x). 
\]%
If $\phi $ is continuous at $\xi _{0}$ (considering the usual definition of 
\textit{continuity }at a point), $\phi \left( \xi _{0}\right) \neq 0,$ and $%
\lim_{\xi \rightarrow \xi _{0}}\frac{\phi \left( \xi \right) -\phi (\xi _{0})%
}{\xi -\xi _{0}}$ does not exist in $\mathbb{R}$, then $f^{m}$ is
differentiable at $\xi _{0}$, and for each $x\in m_{\approx }(\xi _{0})$:%
\[
(f^{m})^{\prime }(x)=mf(\xi _{0})^{m-1}f^{\prime }(x). 
\]

\textbf{f)} For fixed $m$ $\in \mathbb{N}$, let $f$ be differentiable at $%
\xi _{0}$, and let $\phi \left( \xi \right) >0$, for each $\xi \in I$ . Then 
$\sqrt[m]{f}$ is differentiable at $\xi _{0}$, and for each $x\in m_{\approx
}(\xi _{0})$:%
\[
(\sqrt[m]{f})^{\prime }(x)=\frac{f^{\prime }(x)}{m\sqrt[m]{{f\hspace{-2pt}%
\left( \xi _{0}\right) }^{m-1}}}.
\]%
\textbf{Proof \ }This \hspace{1pt}proposition \hspace{1pt}is \hspace{1pt}a 
\hspace{1pt}straightforward \hspace{1pt}consequence \hspace{1pt}of \hspace{%
1pt}\textbf{proposition 3.6} \hspace{1pt}and \hspace{1pt}\textbf{lemma 5.21}%
, except for the fact that we must be very careful with the proviso of 
\textbf{definition 5.16}. To illustrate the last point, we shall prove 
\textbf{c)}.\vspace{-0.25cm}

\textbf{c)} Let $f$ and $g$ be differentiable at $\xi _{0}$, and let $x\in
m_{\approx }(\xi _{0})$.

\vspace{-0.25cm}By \textbf{lemma 5.21 b)}, $fg$ is an indiscernible
extension of $\phi \psi $; so we have, using \textbf{proposition 3.6 b)}:%
\vspace{-0.5cm}

\[
\left( fg\right) \left( x\right) =\left( \phi \psi \right) \hspace{-2pt}%
\left( \xi _{0}\right) +d\left( fg\right) \left( x\right) =\left( \phi \psi
\right) \hspace{-2pt}\left( \xi _{0}\right) +\left( f(\xi _{0})g^{\prime
}(x)+g(\xi _{0})f^{\prime }(x)\right) dx.
\]

Before concluding that $fg$ is differentiable at $\xi _{0}$ and\smallskip 
\[
(fg)^{\prime }(x)=f(\xi _{0})g^{\prime }(x)+g(\xi _{0})f^{\prime }(x), 
\]%
we must be very careful with the proviso of \textbf{definition 5.16}.\vspace{%
-4pt}

If $\lim_{\xi \rightarrow \xi _{0}}\frac{\phi \left( \xi \right) -\phi (\xi
_{0})}{\xi -\xi _{0}},\lim_{\xi \rightarrow \xi _{0}}\frac{\psi \left( \xi
\right) -\psi (\xi _{0})}{\xi -\xi _{0}}$ exist in $\mathbb{R}$, then $%
\lim_{\xi \rightarrow \xi _{0}}\frac{\left( \phi \psi \right) \left( \xi
\right) -\left( \phi \psi \right) (\xi _{0})}{\xi -\xi _{0}}$ also exists in 
$\mathbb{R}$, and equals $f(\xi _{0})g^{\prime }(x)+g(\xi _{0})f^{\prime
}(x) $.\vspace{-4pt}

If $\phi \left( \xi _{0}\right) \neq 0$, $\lim_{\xi \rightarrow \xi _{0}}%
\frac{\phi \left( \xi \right) -\phi (\xi _{0})}{\xi -\xi _{0}}$, $\lim_{\xi
\rightarrow \xi _{0}}\psi \left( \xi \right) $ exist and $\lim_{\xi
\rightarrow \xi _{0}}\frac{\psi \left( \xi \right) -\psi (\xi _{0})}{\xi
-\xi _{0}}$ does not exist in $\mathbb{R}$, then it is easy to prove that $%
\lim_{\xi \rightarrow \xi _{0}}\frac{\left( \phi \psi \right) \left( \xi
\right) -\left( \phi \psi \right) (\xi _{0})}{\xi -\xi _{0}}$ does not exist
in $\mathbb{R}$, and therefore the proviso is not violated.\vspace{-4pt}

When $\psi \left( \xi _{0}\right) \neq $0, $\lim_{\xi \rightarrow \xi _{0}}%
\frac{\psi \left( \xi \right) -\psi (\xi _{0})}{\xi -\xi _{0}}$, $\lim_{\xi
\rightarrow \xi _{0}}~\phi \left( \xi \right) $ exist and $\lim_{\xi
\rightarrow \xi _{0}}\frac{\phi \left( \xi \right) -\phi (\xi _{0})}{\xi
-\xi _{0}}$ does not exist in $\mathbb{R}$, we may use the previous argument
to obtain the same conclusion. $\blacksquare $

\bigskip

\textbf{Theorem 5.23 (Chain Rule) \ }Let $f:m_{\approx }(I)\rightarrow 
\widehat{\mathbb{R}}$, $g:m_{\approx }(J)\rightarrow \widehat{\mathbb{R}}$
be indiscernible extensions of $\phi :I\rightarrow \mathbb{R},\psi
:J\rightarrow \mathbb{R}$, respectively, where $I,J$ are open intervals in $%
\mathbb{R}$ such that $\phi \left( I\right) \subseteq $ $J$, and let $\xi
_{0}\in I$.\vspace{-6pt}

If $f$ is differentiable at $\xi _{0}$, $g$ is differentiable at $\eta
_{0}:= $ $f$\hspace{-2pt}$\left( \xi _{0}\right) $, and both $\lim_{\xi
\rightarrow \xi _{0}}\frac{\phi \left( \xi \right) -\phi (\xi _{0})}{\xi
-\xi _{0}}$ and $\lim_{\eta \rightarrow \eta _{0}}\frac{\psi \left( \eta
\right) -\psi (\eta _{0})}{\eta -\eta _{0}}$ exist in $\mathbb{R}$, then $%
g\circ f$ is differentiable at $\xi _{0}$, and for each $x\in m_{\approx
}(\xi _{0})$:%
\[
(g\circ f)^{\prime }(x)=g^{\prime }(f\left( x\right) )f^{\prime }(x). 
\]

\textbf{Proof \ }Let $f$ be differentiable at $\xi _{0}$, and let $g$ be
differentiable at $\eta _{0}:=$ $f\left( \xi _{0}\right) $.\vspace{-6pt}

By \textbf{lemma 5.21 e)}, $g\circ f$ $\ $is an indiscernible extension of $%
\psi \circ \phi $; so we have, for each $x\in m_{\approx }(\xi _{0})$:%
\[
(g\circ f)(x)=\left( \psi \circ \phi \right) \hspace{-2pt}\left( \xi
_{0}\right) +d(g\circ f)(x)=\psi \hspace{-2pt}\left( \phi \hspace{-2pt}%
\left( \xi _{0}\right) \right) +d(g\circ f)(x)=\psi \hspace{-2pt}\left( f%
\hspace{-2pt}\left( \xi _{0}\right) \right) +d(g\circ f)(x)=
\]%
\[
=\psi \hspace{-2pt}\left( \eta _{0}\right) +d(g\circ f)(x)=g(\eta
_{0})+d(g\circ f)(x).
\]%
On the other hand, since $f$ is differentiable at $\xi _{0}$, and $g$ is
differentiable at $\eta _{0}=\nolinebreak f\left( \xi _{0}\right) ,$%
\[
(g\circ f)(x)=g(f\left( x\right) )=g(\eta _{0}{\ +}f^{\prime }(x)dx)=g%
\hspace{-2pt}\left( \eta _{0}\right) +g^{\prime }\hspace{-2pt}\left( \eta
_{0}\right) f^{\prime }(x)dx.
\]%
By comparison with the previous result for $(g\circ f)(x),$ we infer that%
\[
d(g\circ f)(x)=g^{\prime }\hspace{-3pt}\left( \eta _{0}\right) f^{\prime
}(x)dx.
\]%
Since $f\hspace{-2pt}\left( x\right) \in m_{\approx }(\eta _{0})$ (because $f
$ is differentiable at $\xi _{0}$), and $g$ is differentiable at $\eta _{0},$
we have:%
\[
d(g\circ f)(x)=g^{\prime }\hspace{-2pt}\left( f\left( x\right) \right)
f^{\prime }(x)dx..
\]%
And the proviso of \textbf{definition 5.16} is satisfied, since we obtain,
as an immediate consequence of the usual \textbf{Chain Rule} in $\mathbb{R}$
(and the differentiability of $f,g$ at $\xi _{0},\eta _{0},$ respectively) : 
\[
\lim_{\xi \rightarrow \xi _{0}}\frac{\left( \psi \circ \phi \right) \left(
\xi \right) -\left( \psi \circ \phi \right) (\xi _{0})}{\xi -\xi _{0}}%
=\left( \lim_{\eta \rightarrow \eta _{0}}\frac{\psi \left( \eta \right)
-\psi (\eta _{0})}{\eta -\eta _{0}}\right) \left( \lim_{\xi \rightarrow \xi
_{0}}\frac{\phi \left( \xi \right) -\phi (\xi _{0})}{\xi -\xi _{0}}\right) =
\]%
\[
=g^{\prime }\hspace{-2pt}\left( \eta _{0}\right) f^{\prime }(\xi
_{0})=g^{\prime }(f(x))f^{\prime }(x).
\]%
We have proven that $g\circ f$ is differentiable at $\xi _{0}$, and for each 
$x\in m_{\approx }(\xi _{0})$:%
\[
(g\circ f)^{\prime }(x)=g^{\prime }(f\left( x\right) )f^{\prime }(x).\text{ }%
\blacksquare 
\]

\bigskip

\textbf{Theorem 5.24} \textbf{(The Inverse Function Theorem) \ }Let $I$ be
an open interval in $\mathbb{R}$, let $f:m_{\approx }(I)\rightarrow \widehat{%
\mathbb{R}}$ be an injective indiscernible extension of a continuous
function $\phi :I\rightarrow \mathbb{R}$ (we consider the usual topology for 
$%
%TCIMACRO{\U{211d} }%
%BeginExpansion
\mathbb{R}
%EndExpansion
,$ and its relativization to $I$), and let $\xi _{0}\in I$.\vspace{-5pt}

If $m_{\approx }\hspace{-2pt}\left( \phi \hspace{-2pt}\left( I\right)
\right) \subseteq f\hspace{-2pt}\left( m_{\approx }\hspace{-2pt}\left(
I\right) \right) ,$ $f$ is differentiable at $\xi _{0}$, $f^{\prime }(\xi
_{0})\neq 0$, and $\lim_{\xi \rightarrow \xi _{0}}\frac{\phi \left( \xi
\right) -\phi (\xi _{0})}{\xi -\xi _{0}}$ exists in $\mathbb{R}$, then $%
f^{-1}$ (considered as a function with codomain $\widehat{%
%TCIMACRO{\U{211d} }%
%BeginExpansion
\mathbb{R}
%EndExpansion
}$) is differentiable at $\eta _{0}:=$ $f\hspace{-2pt}\left( \xi _{0}\right) 
$, and we have, for each $y\in m_{\approx }(\eta _{0})$: 
\[
\left( f^{-1}\right) ^{\prime }\hspace{-2pt}\left( y\right) =\frac{1}{%
f^{\prime }\hspace{-2pt}\left( f^{-1}\hspace{-2pt}\left( y\right) \right) }.
\]%
\vspace{4pt}

\textbf{Proof \ }Let $J:=\phi (I)$, $\alpha :=\lim_{\xi \rightarrow \xi _{0}}%
\frac{\phi \left( \xi \right) -\phi (\xi _{0})}{\xi -\xi _{0}}$, and $\eta
_{0}:=$ $f\hspace{-2pt}\left( \xi _{0}\right) $.\vspace{0.1cm}\newline
Since $\phi $ is continuous and injective (because $f$ is an injective
indiscernible extension of $\phi $), $\phi ^{-1}$ is also continuous
(considering the usual topology for $%
%TCIMACRO{\U{211d} }%
%BeginExpansion
\mathbb{R}
%EndExpansion
,$ and its relativization to $J$). So $J=\phi (I)=\left( \phi ^{-1}\right)
^{-1}\hspace{-2pt}\left( I\right) $ is an open interval in $\mathbb{R}$, and
the same is valid for $m_{\approx }(J)$ in $\widehat{\mathbb{R}}$ (see 
\textbf{proposition 5.12 a)}).\vspace{-6pt}

As $\alpha =f^{\prime }(\xi _{0})\neq 0,$ we know, by the usual \textbf{%
Inverse Function Theorem} in $\mathbb{R}$, that $\beta :=\lim_{\eta
\rightarrow \eta _{0}}\frac{\phi ^{-1}\left( \eta \right) -\phi ^{-1}(\eta
_{0})}{\eta -\eta _{0}}$ exists in $\mathbb{R}$, and \vspace{-4pt}%
\[
\beta =\frac{1}{\alpha }.\smallskip 
\]%
Since $m_{\approx }(J)=f\left( m_{\approx }\hspace{-2pt}\left( I\right)
\right) $ (see the proof of \textbf{lemma 5.21 f)}), we may consider the
function $g:m_{\approx }(J)\rightarrow \widehat{\mathbb{R}}$ defined by%
\[
g(y):=\left\{ 
\begin{array}{c}
f^{-1}\hspace{-2pt}\left( y\right) {,}\text{ {if \ }}y\in m_{\approx }(J){%
\backslash }m_{\approx }(\eta _{0}) \\ 
\phi ^{-1}\hspace{-3pt}\left( \eta _{0}\right) +\beta dy{,}\text{ {if \ }}%
y\in m_{\approx }(\eta _{0})%
\end{array}%
\right. .
\]%
Since, by \textbf{lemma 5.21 f)}, $f^{-1}$ is an indiscernible extension of $%
\phi ^{-1}$, to complete the prove we only need to show that $g(y)=f^{-1}(y)$%
, for each $y\in m_{\approx }(\eta _{0})$.\newline
If $y\in m_{\approx }(\eta _{0})$, then $g(y)\in m_{\approx }(\xi _{0})$
(because $\phi ^{-1}\hspace{-3pt}\left( \eta _{0}\right) =f^{-1}\hspace{-3pt}%
\left( \eta _{0}\right) =\xi _{0}$), and since $f$ is differentiable at $\xi
_{0},$ we have$:$%
\[
f\hspace{-2pt}\left( g(y)\right) =f(\phi ^{-1}\hspace{-3pt}\left( \eta
_{0}\right) +\beta dy)=f(\xi _{0}+\beta dy)=f(\xi _{0})+\alpha \beta dy=\eta
_{0}+dy=y.
\]%
So%
\[
g(y)=f^{-1}\hspace{-2pt}\left( f\hspace{-2pt}\left( g(y)\right) \right)
=f^{-1}\hspace{-2pt}\left( y\right) .\text{ }\blacksquare \medskip \vspace{%
4pt}
\]%
\textbf{Notation.} Let $I$ be a nonempty open interval in $\mathbb{R}$, let $%
\phi :I\rightarrow \mathbb{R}$ be a function, and let 
\[
\mathit{\Lambda }_{\phi }:=\left\{ \xi _{0}\in I|\lim_{\xi \rightarrow \xi
_{0}}\frac{\phi \left( \xi \right) -\phi (\xi _{0})}{\xi -\xi _{0}}\text{
exists in }{\mathbb{R}}\right\} .\vspace{2pt}
\]%
The \hspace{1pt}function $\hspace{1pt}\xi _{0}\in \mathit{\Lambda }_{\phi
}\mapsto \lim_{\xi \rightarrow \xi _{0}}\frac{\phi \left( \xi \right) -\phi
(\xi _{0})}{\xi -\xi _{0}}$, \hspace{1pt}from $\hspace{1pt}\mathit{\Lambda }%
_{\phi }$ to $\mathbb{R}$ , \hspace{1pt}will \hspace{1pt}be \hspace{1pt}%
denoted \hspace{1pt}by $\hspace{1pt}\lambda _{\phi }$ \vspace{2pt}\newline
(notice that we do not exclude, at least here, the case $\mathit{\Lambda }%
_{\phi }=\emptyset $).

\bigskip

\textbf{Theorem 5.25} \textbf{(The Mean Value Theorem)\vspace{-7pt}\vspace{%
1pt}}

Let $I$ be a nonempty open interval in $\mathbb{R}$, let $f:m_{\approx
}(I)\rightarrow \widehat{\mathbb{R}}$ be an indiscernible \vspace{1pt}%
\newline
extension of $\phi :I\rightarrow \mathbb{R}$, differentiable on $m_{\approx
}(I)$, and let $\mathit{\Lambda }_{\phi }=I$ \textbf{.\vspace{-7pt}\vspace{%
1pt}}

If $a,b$ $\in m_{\approx }(I)$ and $a<b$, then there exists $\gamma \in I$
such that $a<\gamma <b$, and%
\begin{equation}
f\hspace{-2pt}\left( b\right) -f\hspace{-2pt}\left( a\right) =f^{\prime }%
\hspace{-2pt}\left( \gamma \right) \left( b-a\right) +\left( f^{\prime }%
\hspace{-2pt}\left( b\right) -f^{\prime }(\gamma )\right) \hspace{-2pt}%
db+\left( f^{\prime }\hspace{-2pt}\left( \gamma \right) -f^{\prime
}(a)\right) \hspace{-2pt}da.  \label{GrindEQ__1_}
\end{equation}%
In particular, if $a,b\in I$\textit{,} then $(1)$ assumes the familiar form:%
\[
f\hspace{-2pt}\left( b\right) -f\hspace{-2pt}\left( a\right) =f^{\prime }%
\hspace{-2pt}\left( \gamma \right) \left( b-a\right) .
\]

\vspace{-0.25cm}The previous identities stay valid when we replace $\gamma $
by any $c$ $\in m_{\approx }(\gamma )$.

\bigskip

\textbf{Proof \ }Let $a,b$ $\in m_{\approx }(I)$, and $a<b$.\vspace{-6pt}

Then%
\[
\sigma a<\sigma b. 
\]%
So, by the usual \textbf{Mean Value Theorem}, there is $\gamma \in I$ such
that $\sigma a<\gamma <\sigma b$\textit{,} and%
\[
\phi \left( \sigma b\right) -\phi \left( \sigma a\right) =\lambda _{\phi
}(\gamma )\left( \sigma b-\sigma a\right) . 
\]

Then, since $f$ is an indiscernible extension of $\phi $, differentiable on $%
m_{\approx }(I)$, we have:

$f\hspace{-2pt}\left( b\right) -f\hspace{-2pt}\left( a\right) ={\phi \hspace{%
-1pt}\left( \sigma b\right) +f}^{\prime }\hspace{-2pt}\left( b\right) 
\hspace{-1pt}db-\phi \hspace{-1pt}\left( \sigma a\right) {-f}^{\prime }%
\hspace{-2pt}\left( a\right) \hspace{-1pt}da=f^{\prime }(\gamma )\left(
b-db-a+da\right) +$\newline
$+$ $f^{\prime }\hspace{-2pt}\left( b\right) \hspace{-1pt}db-f^{\prime }%
\hspace{-2pt}\left( a\right) \hspace{-1pt}da=f^{\prime }(\gamma )\left(
b-a\right) +\left( f^{\prime }\hspace{-2pt}\left( b\right) -f^{\prime
}(\gamma )\right) \hspace{-1pt}db+\left( f^{\prime }\hspace{-2pt}\left(
\gamma \right) -f^{\prime }(a)\right) \hspace{-1pt}da.$

Finally, by \textbf{definition 5.18}, $f^{\prime }\hspace{-2pt}\left( \gamma
\right) =f^{\prime }\hspace{-2pt}\left( c\right) $, for each $c$ $\in
m_{\approx }\hspace{-1pt}(\gamma )$. $\blacksquare $

\bigskip

\textbf{Corollary 5.26 \ }Let $I$ be a nonempty open interval in $\mathbb{R}$%
, let $f:m_{\approx }(I)\rightarrow \widehat{\mathbb{R}}$ be an
indiscernible extension of $\phi :I\rightarrow \mathbb{R}$, differentiable
on $m_{\approx }(I),$ and let $\mathit{\Lambda }_{\phi }=I$ .\vspace{-4pt}

\textbf{a)} If $f^{\prime }(x)=0$, for each $x$ $\in m_{\approx }(I)$, then $%
f$ is a constant function.\vspace{-4pt}

\textbf{b)} If $f^{\prime }\left( x\right) >0$, for each $x$ $\in m_{\approx
}(I)$, then $f$ is a strictly increasing function.\vspace{-4pt}

\textbf{b)} If $f^{\prime }\left( x\right) <0$, for each $x$ $\in m_{\approx
}(I)$, then $f$ is a strictly decreasing function.

\bigskip

\textbf{Proof \ a) }Let $a,b$ $\in m_{\approx }(I)$.

\vspace{-0.25cm}If $a\approx b,$ then, since $f$ is differentiable on $%
m_{\approx }(I)$ with null derivative,%
\[
f\hspace{-2pt}\left( a\right) =\phi \hspace{-2pt}\left( \sigma a\right)
+f^{\prime }\hspace{-3pt}\left( a\right) da=\phi \hspace{-2pt}\left( \sigma
a\right) =\phi \hspace{-2pt}\left( \sigma b\right) =\phi \hspace{-2pt}\left(
\sigma b\right) +f^{\prime }\hspace{-3pt}\left( b\right) db=f\hspace{-2pt}%
\left( b\right) . 
\]%
If $a<b$ or $b<a,$ then we obtain, as a direct consequence of identity $%
\left( 1\right) $ of \textbf{theorem 5.25},%
\[
f\left( a\right) =f\left( b\right) . 
\]%
\textbf{b) }and \textbf{c) }admit trivial proofs, since if $a,b$ $\in
m_{\approx }(I)$ and $a<b,$ then we easily obtain, using identity $\left(
1\right) $ of \textbf{theorem 5.25}: 
\[
\sigma f(b)-\sigma f(a)=f^{\prime }\hspace{-2pt}\left( \gamma \right) \left(
\sigma b-\sigma a\right) ,\text{ for a certain }\gamma \text{ such that }%
a<\gamma <b.\text{ }\blacksquare 
\]

\bigskip

We close this section with the introduction and elementary study of the
concept of \textit{nat\nolinebreak ural indiscernible extension} of a
function $\phi :I\rightarrow \mathbb{R}$, where $I$ is a nonempty open
interval in $\mathbb{R}$ . \textit{Natural indiscernible extensions} are the 
\og natural\fg \ versions, in $\widehat{\mathbb{R}}$,
of the usual differentiable functions, in $\mathbb{R}$.

\bigskip

The starting point is the next proposition, which follows immediately from 
\textbf{definition 5.16} and \textbf{remark 5.19}.

\bigskip

\textbf{Proposition 5.27 \ }Let $I$ be a nonempty open interval in $\mathbb{R%
}$, and let $\phi :I\rightarrow \mathbb{R}$ be a function such that $\mathit{%
\Lambda }_{\phi }=I$ .

\vspace{-0.25cm}Then the function $f:m_{\approx }(I)\rightarrow \widehat{%
\mathbb{R}}$, defined by $f(x):=\phi (\sigma x)+\lambda _{\phi }(\sigma x)dx$%
, is the unique indiscernible extension of $\phi $ differentiable on $%
m_{\approx }(I)$.

\bigskip

\textbf{Definition 5.28 \ }With the notation and the hypothesis of \textbf{%
proposition 5.27}, we call $f:m_{\approx }(I)\rightarrow \widehat{\mathbb{R}}
$, defined by $f(x):=\phi (\sigma x)+\lambda _{\phi }(\sigma x)dx$, the 
\textit{natural indiscernible extension }of $\phi $, and we denote it by $%
\hat{\phi}$.

\bigskip

\textit{Natural indiscernible extensions} preserve addition, scalar
multiplication by a real number, multiplication, division, composition, and
inversion, in a sense clearly expressed by \textbf{a)} to \textbf{e)}, and 
\textbf{g)}, in the next proposition.

\bigskip

\textbf{Proposition 5.29 \ }Let $I$ be a nonempty open interval in $\mathbb{R%
}$, and let $\phi :I\rightarrow \mathbb{R},\psi :\nolinebreak I\rightarrow
\nolinebreak \mathbb{R}$ be functions such that $\mathit{\Lambda }_{\phi }=%
\mathit{\Lambda }_{\psi }=I$\textit{\ .}

\textbf{a)} $\mathit{\Lambda }_{\phi +\psi }=I,$ and $\widehat{\phi +\psi }=%
\hat{\phi}+\hat{\psi}$ .

\textbf{b)} If $\alpha \in \mathbb{R},$ then $\mathit{\Lambda }_{\alpha \phi
}=I,$ and $\widehat{\alpha \phi }\mathbf{=}\alpha \hat{\phi}$.

\textbf{c)} $\mathit{\Lambda }_{\phi \psi }=I,$ and $\widehat{\phi \psi }=%
\hat{\phi}\hat{\psi}$ .

\textbf{d)} If $\psi \left( \xi \right) \neq 0$, for each $\xi \in I$ , then 
$\mathit{\Lambda }_{\frac{\phi }{\psi }}=I,$ and 
\[
\widehat{\left( \frac{\phi }{\psi }\right) }=\frac{\hat{\phi}}{\hat{\psi}}. 
\]%
\textbf{e)} If $J$ is a nonempty open interval in $\mathbb{R},$ $\theta
:J\rightarrow \mathbb{R}$ is a function such that $\phi (I)\subseteq J,$ and 
$\mathit{\Lambda }_{\theta }=J$, then $\mathit{\Lambda }_{\theta \circ \phi
}=I,$ and%
\[
\widehat{\theta \circ \phi }=\hat{\theta}\circ \hat{\phi}. 
\]%
\textbf{f)} If $A$ is a nonempty subset of $I$, then%
\[
\hat{\phi}(m_{\approx }(A))\subseteq m_{\approx }(\phi (A)). 
\]%
If $\lambda _{\phi }(\xi )\neq 0$, for each $\xi \in I$, then%
\[
\hat{\phi}(m_{\approx }(A))=m_{\approx }(\phi (A)). 
\]%
If $\alpha ,\beta \in I,\alpha <\beta $, $\lambda _{\phi }(\xi )\neq 0$, for
each $\xi \in \left] \alpha ,\beta \right[ $, and $\lambda _{\phi }(\alpha
)=\lambda _{\phi }(\beta )=0$, then%
\[
\hat{\phi}\left( \widehat{\left[ \alpha ,\beta \right] }\right) =m_{\approx
}(\phi (\left] \alpha ,\beta \right[ ))\cup \left\{ \phi \left( \alpha
\right) ,\phi \left( \beta \right) \right\} . 
\]%
\textbf{g)} If $\phi $ is continuous, injective, and $\lambda _{\phi }(\xi
)\neq 0$, for each $\xi \in I$, then $\hat{\phi}$ is injective, $\mathit{%
\Lambda }_{\phi ^{-1}}=\phi \left( I\right) ,$ and 
\[
\widehat{\phi ^{-1}}={\hat{\phi}}^{-1}, 
\]%
considering $\phi ^{-1},$ ${\hat{\phi}}^{-1}$ as functions with codomains $%
%TCIMACRO{\U{211d} }%
%BeginExpansion
\mathbb{R}
%EndExpansion
,$ $\widehat{%
%TCIMACRO{\U{211d} }%
%BeginExpansion
\mathbb{R}
%EndExpansion
},$ respectively.

\textbf{h)} If $I=\mathbb{R}$ and $\phi $ is an even function, then $\hat{%
\phi}$ is also an even function, i.e. $\hat{\phi}(-x)=\nolinebreak \hat{\phi}%
(x)$, for each $x\in \widehat{\mathbb{R}}$ .\vspace{-6pt}

Similarly, if $I=\mathbb{R}$ and $\phi $ is an odd function, then $\hat{\phi}
$ is also an odd function, i.e. $\hat{\phi}(-x)=-\hat{\phi}(x)$, for each $%
x\in \widehat{\mathbb{R}}$ .\vspace{-6pt}

\textbf{i)} If $I=\mathbb{R}$ , $\lambda _{0}\in {\mathbb{R}}^{+}$, and $%
\phi $ is a periodic function with period $\lambda _{0}$, then $\hat{\phi}$
is also periodic with the same real period, i.e. 
\[
\lambda _{0}=\text{min}_{r}\left\{ l\in {\widehat{\mathbb{R}}}^{+}|\left(
\forall x\in \widehat{\mathbb{R}}\right) \hat{\phi}(x+l)=\hat{\phi}%
(x)\right\} . 
\]%
\textbf{Proof a)} and \textbf{b)} admit trivial proofs, using the
well-known identities (with different notation\textbf{) }$\mathit{\Lambda }%
_{\phi +\psi }=\mathit{\Lambda }_{\alpha \phi }=I,$ and $\lambda _{\phi
+\psi }=\lambda _{\phi }+\lambda _{\psi },\lambda _{\alpha \phi }=\alpha
\lambda _{\phi }$ .

\vspace{-0.25cm}\textbf{c)} Clearly, $\mathit{\Lambda }_{\phi \psi }=I,$ and
for each $x\in m_{\approx }(I)$, we have, using the well--known identity
(with different notation\textbf{) }$\lambda _{\phi \psi }=\lambda _{\phi
}\psi +\lambda _{\psi }\phi $ :

$\widehat{\phi \psi }(x)=(\phi \psi )(\sigma x)+\lambda _{\phi \psi }(\sigma
x)dx=\phi (\sigma x)\psi (\sigma x)+(\lambda _{\phi }(\sigma x)\psi \left(
\sigma x\right) +\lambda _{\psi }(\sigma x)\phi (\sigma x))dx=$

$=\hspace{-2pt}\left( \phi \left( \sigma x\right) \hspace{-1pt}+\hspace{-2pt}%
{\lambda }_{\phi }(\sigma x)dx\right) \hspace{-2pt}\psi (\sigma x)+\lambda
_{\psi }(\sigma x)\phi (\sigma x)dx=\hat{\phi}(x)\psi (\sigma x)+(\psi
(\sigma x)+\lambda _{\psi }(\sigma x)dx)\phi (\sigma x)-\vspace{0.15cm}$%
\newline
$-\phi (\sigma x)\psi (\sigma x)=\hat{\phi}(x)\psi (\sigma x)+\hat{\psi}%
(x)\phi (\sigma x)-\phi (\sigma x)\psi (\sigma x)=\widehat{\phi }(x)\hspace{%
-3pt}\left( \sigma \hat{\psi}(x)\right) +\hat{\psi}(x)\hspace{-3pt}\left(
\sigma \hat{\phi}(x)\right) -\smallskip $\newline
$\smallskip -\left( \sigma \hat{\phi}(x)\right) \hspace{-4pt}\left( \sigma 
\hat{\psi}(x)\right) =\hat{\phi}(x)\hat{\psi}(x)=\left( \hat{\phi}\hat{\psi}%
\right) \hspace{-3pt}(x)$ .

\textbf{d)} Clearly $\mathit{\Lambda }_{\frac{1}{\psi }}=I$ $;$ and for each 
$x\in m_{\approx }(I)$, we have, using the well-known identity (with
different notation) $\lambda _{\frac{1}{\psi }}=-\frac{\lambda _{\psi }}{%
\psi ^{2}}$ :\vspace{-12pt}

\[
\widehat{\left( \frac{1}{\psi }\right) }(x)=\frac{1}{\psi }(\sigma
x)+\lambda _{\frac{1}{\psi }}(\sigma x)dx=\frac{1}{\psi (\sigma x)}-\frac{%
\lambda _{\psi }(\sigma x)}{{\psi \hspace{-2pt}\left( \sigma x\right) }^{2}}%
dx=\frac{\psi (\sigma x{)-}\lambda _{\psi }(\sigma x)dx}{{\psi \hspace{-2pt}%
\left( \sigma x\right) }^{2}}. 
\]%
But (since the product of infinitesimals is always null) 
\[
\frac{\psi (\sigma x{)-}\lambda _{\psi }(\sigma x)dx}{{\psi \hspace{-2pt}%
\left( \sigma x\right) }^{2}}(\psi (\sigma x)+\lambda _{\psi }(\sigma x)dx)=%
\frac{{\psi \left( \sigma x\right) }^{2}{-}\lambda _{\psi }(\sigma x)\psi
(\sigma x)dx\hspace{-1pt}+\hspace{-2pt}\lambda _{\psi }(\sigma x)\psi
(\sigma x)dx}{{\psi \hspace{-2pt}\left( \sigma x\right) }^{2}}=1. 
\]%
So 
\[
\frac{\psi (\sigma x{)-}\lambda _{\psi }(\sigma x)dx}{{\psi \left( \sigma
x\right) }^{2}}=\frac{1}{\psi (\sigma x)+\lambda _{\psi }(\sigma x)dx}=\frac{%
1}{\hat{\psi}(x)}=\frac{1}{\hat{\psi}}(x). 
\]%
We have proven that 
\[
\widehat{\left( \frac{1}{\psi }\right) }=\frac{1}{\hat{\psi}}. 
\]%
Finally, using \textbf{c)}, we have:%
\[
\mathit{\Lambda }_{\frac{\phi }{\psi }}=\mathit{\Lambda }_{\phi \frac{1}{%
\psi }}=I, 
\]%
\newline
and 
\[
\widehat{\left( \frac{\phi }{\psi }\right) }=\hat{\phi}\widehat{\left( \frac{%
1}{\psi }\right) }=\hat{\phi}\text{ }\frac{1}{\hat{\psi}}=\frac{\hat{\phi}}{%
\hat{\psi}}. 
\]%
\textbf{e)} Clearly, $\mathit{\Lambda }_{\theta \circ \phi }=I$ $;$ and for
each $x\in m_{\approx }(I)$, we have, using the well-known identity (with
different notation) $\lambda _{\theta \circ \phi }={(\lambda }_{\theta
}\circ \phi )\lambda _{\phi }$ :

$\widehat{\theta \circ \phi }(x)=(\theta \circ \phi )(\sigma x)+\lambda
_{\theta \circ \phi }(\sigma x)dx=\theta (\phi (\sigma x))+{\ (\lambda }%
_{\theta }\circ \phi )(\sigma x)\lambda _{\phi }(\sigma x)dx=\theta (\phi
(\sigma x))+\vspace{-4pt}\smallskip $

$+{\lambda }_{\theta }(\phi (\sigma x))\lambda _{\phi }(\sigma x)dx=\theta
(\phi (\sigma x))+{\lambda }_{\theta }(\phi (\sigma x))d\hat{\phi}(x)=\theta 
\hspace{-3pt}\left( \sigma \hat{\phi}(x)\right) +{\lambda }_{\theta }\hspace{%
-3pt}\left( \sigma \hat{\phi}(x)\right) \hspace{-3pt}d\hat{\phi}(x)=$\newline
$=\hat{\theta}\hspace{-3pt}\left( \hat{\phi}\hspace{1pt}(x)\right) =\left( 
\hat{\theta}\circ \hat{\phi}\right) \hspace{-3pt}(x)$.

\textbf{f)} If $x\in m_{\approx }(A)$, then, since $\sigma x\in A$ (by 
\textbf{proposition 4.4 b)}, \textbf{c)}),%
\[
\hat{\phi}(x)=\phi (\sigma x)+{\ \lambda }_{\phi }(\sigma x)dx\in m_{\approx
}(\phi (A)). 
\]%
We have proven that 
\[
\hat{\phi}(m_{\approx }(A))\subseteq m_{\approx }(\phi (A)). 
\]%
Let $\lambda _{\phi }(\xi )\neq 0$, for each $\xi \in I$ .\vspace{-4pt}

If $\xi \in A$ and $\hat{\varepsilon}\approx 0$, then 
\[
\phi (\xi )+\hat{\varepsilon}=\phi (\xi )+\lambda _{\phi }(\xi )\frac{1}{%
\lambda _{\phi }(\xi )}\hat{\varepsilon}=\hat{\phi}\hspace{-2pt}\left( \xi +%
\frac{1}{\lambda _{\phi }(\xi)}\hat{\varepsilon}\right) \in \hat{\phi}%
(m_{\approx }(A)).
\]%
We have proven that%
\[
m_{\approx }(\phi (A))\subseteq \hat{\phi}(m_{\approx }(A)).
\]%
Let $\alpha ,\beta \in I$, with $\alpha <\beta $, let $\lambda _{\phi }(\xi
)\neq 0$, for each $\xi \in \left] \alpha ,\beta \right[ $, and let $\lambda
_{\phi }(\alpha )=\nolinebreak \lambda _{\phi }(\beta )=\nolinebreak 0$.

We have: 
\[
\hat{\phi}\hspace{-2pt}\left( \widehat{\left[ \alpha ,\beta \right] }\right)
=\hat{\phi}\hspace{-2pt}\left( \widehat{\left] \alpha ,\beta \right[ }\cup
m_{\approx }(\alpha )\cup m_{\approx }(\beta )\right) =\hat{\phi}\hspace{-2pt%
}\left( \widehat{\left] \alpha ,\beta \right[ }\right) \cup {\hat{\phi}(m}%
_{\approx }(\alpha ))\cup \hat{\phi}(m_{\approx }(\beta )).
\]%
Since $\lambda _{\phi }(\xi )\neq 0$, for each $\xi \in \left] \alpha ,\beta %
\right[ $, and $\lambda _{\phi }(\alpha )=\lambda _{\phi }(\beta )=0$, we
obtain (using the result we have just proven, and \textbf{definition 5.28}): 
\[
\hat{\phi}\hspace{-2pt}\left( \widehat{\left] \alpha ,\beta \right[ }\right)
=\hat{\phi}\left( m_{\approx }(\left] \alpha ,\beta \right[ )\right)
=m_{\approx }(\phi (\left] \alpha ,\beta \right[ )),
\]%
\[
\hat{\phi}(m_{\approx }(\alpha ))=\left\{ \phi (\alpha )\right\} ,
\]%
\[
\hat{\phi}(m_{\approx }(\beta ))=\left\{ \phi (\beta )\right\} .
\]%
So%
\[
\hat{\phi}\hspace{-2pt}\left( \widehat{\left[ \alpha ,\beta \right] }\right)
=m_{\approx }(\phi (\left] \alpha ,\beta \right[ ))\cup \left\{ \phi (\alpha
),\phi (\beta )\right\} .
\]

\textbf{g)} Clearly, $\phi \left( I\right) $ is a nonempty open interval,
and by the usual \textbf{Inverse Function The\nolinebreak orem}, $\mathit{%
\Lambda }_{\phi ^{-1}}=\phi \left( I\right) .\vspace{-0.2cm}$

On the other hand, for each $x_{1},x_{2}\in m_{\approx }(I)$, we have (since 
$\phi $ is injective and $\lambda _{\phi }(\xi )\neq \nolinebreak 0$, for
each $\xi \in I$):%
\[
\hat{\phi}(x_{1})=\hat{\phi}(x_{2})\Rightarrow \phi (\sigma x_{1})+\lambda
_{\phi }(\sigma x_{1})dx_{1}=\phi (\sigma x_{2})+\lambda _{\phi }(\sigma
x_{2})dx_{2}\Rightarrow 
\]%
\[
\Rightarrow \left\{ 
\begin{array}{c}
\phi \left( \sigma x_{1}\right) =\phi (\sigma x_{2}) \\ 
\lambda _{\phi }(\sigma x_{1})dx_{1}=\lambda _{\phi }(\sigma x_{2})dx_{2}%
\end{array}%
\right. \Rightarrow \left\{ 
\begin{array}{c}
\sigma x_{1}=\sigma x_{2} \\ 
\lambda _{\phi }(\sigma x_{1}{)(}dx_{1}-dx_{2})=0%
\end{array}%
\right. \Rightarrow 
\]%
\[
\Rightarrow \left\{ 
\begin{array}{c}
\sigma x_{1}=\sigma x_{2} \\ 
dx_{1}=dx_{2}%
\end{array}%
\right. \Rightarrow x_{1}=x_{2}. 
\]%
So $\hat{\phi}$ is also injective.\vspace{-5pt}

If $x\in m_{\approx }(I),$ then we have, using \textbf{e) }and denoting by $%
\iota _{I}$ the inclusion function of $I$ into $%
%TCIMACRO{\U{211d} }%
%BeginExpansion
\mathbb{R}
%EndExpansion
:$%
\[
\left( \widehat{\phi ^{-1}}\circ \hat{\phi}\right) \hspace{-3pt}\left(
x\right) =\widehat{\phi ^{-1}\circ \phi }\left( x\right) =\widehat{\iota _{I}%
}\left( x\right) =\iota _{I}(\sigma x)+{\lambda _{\iota _{I}}}(\sigma
x)dx=\sigma x+dx=x. 
\]%
If $y\in m_{\approx }(\phi \left( I\right) )$, we have, using \textbf{e)}
and denoting by $\iota _{\phi \left( I\right) }$ the inclusion function of $%
\phi \left( I\right) $ into $\mathbb{R}$ :%
\[
(\hat{\phi}\circ \widehat{\phi ^{-1}})(y)=\widehat{\phi \circ \phi ^{-1}}(y)=%
\widehat{\iota _{\phi \left( I\right) }}(y)=\iota _{\phi \left( I\right)
}(\sigma y)+{\lambda }\iota _{\phi \left( I\right) }(\sigma y)dy=\sigma
y+dy=y. 
\]%
Finally, since the domains of $\widehat{\phi ^{-1}},$ ${\hat{\phi}}^{-1}$
are $m_{\approx }(\phi (I),$ $\hat{\phi}(m_{\approx }(I)),$ and these sets
are identical, by \textbf{f)}, we may consider proven that%
\[
\widehat{\phi ^{-1}}={\hat{\phi}}^{-1}, 
\]%
viewing $\phi ^{-1},$ ${\hat{\phi}}^{-1}$ as functions with codomains $%
%TCIMACRO{\U{211d} }%
%BeginExpansion
\mathbb{R}
%EndExpansion
,$ $\widehat{%
%TCIMACRO{\U{211d} }%
%BeginExpansion
\mathbb{R}
%EndExpansion
},$ respectively.\vspace{-2pt}

\textbf{h)} admits a trivial proof, since%
\[
\phi \text{ is an even function}\Rightarrow \lambda _{\phi }\text{ is an odd
function,} 
\]%
\[
\phi \text{ is an odd function}\Rightarrow \lambda _{\phi }\text{ is an even
function,} 
\]%
\vspace{-0.5cm}\newline
and $d(-x)=-dx$, for each $x\in \widehat{\mathbb{R}}$ .

\textbf{i)} Let $\lambda _{0}\in {\mathbb{R}}^{+}$, let $\phi :\mathbb{R}%
\rightarrow \mathbb{R}$ be a periodic function with period $\lambda _{0}$,
and let 
\[
\hat{L}:=\left\{ l\in {\widehat{\mathbb{R}}}^{+}|\left( \forall x\in 
\widehat{\mathbb{R}}\right) \hat{\phi}(x+l)=\hat{\phi}(x)\right\} . 
\]%
For each $x\in \widehat{\mathbb{R}}$ , we have, using the well-known fact
that $\lambda _{\phi }$ is also periodic with period $\lambda _{0}$:%
\[
\hat{\phi}(x+\lambda _{0})=\phi (\sigma x+\lambda _{0})+{\lambda }_{\phi
}(\sigma x+\lambda _{0})dx=\phi (\sigma x)+{\lambda }_{\phi }(\sigma x)dx=%
\hat{\phi}(x). 
\]

Then, since ${\mathbb{R}}^{+}\subseteq \widehat{\mathbb{R}}^{+}$, we infer
that 
\[
\lambda _{0}\in \hat{L}. 
\]%
On the other hand, if $l\in \hat{L}$, we have, for each $\xi \in \mathbb{R}:$%
\[
\phi (\xi +\sigma l)=\sigma \hat{\phi}(\xi +l)=\sigma \hat{\phi}(\xi )=\phi
(\xi ). 
\]%
Then, since $\sigma l\in {\mathbb{R}}^{+}$, 
\[
\lambda _{0}\leq \sigma l. 
\]%
So 
\[
\lambda _{0}\lesssim l. 
\]%
Since $\lambda _{0}\in {\mathbb{R}}$, $\lambda _{0}\in \hat{L}$, and $%
\lambda _{0}$ is an $\lesssim $-lower bound of $\hat{L}$, we conclude that 
\[
\lambda _{0}=\text{min}_{r}\hat{L}.\text{ }\blacksquare 
\]

\bigskip

Frequently, physicists and engineers use identities like%
\[
(1+dx)^{\alpha }=1+\alpha dx\text{ \ (for fixed }\alpha \in \mathbb{R}\text{)%
}, 
\]%
\[
\sin \left( dx\right) =dx, 
\]%
\[
\cos \left( dx\right) =1, 
\]%
\[
\exp \left( dx\right) =1+dx,\nolinebreak 
\]%
\[
\log (1+dx)=dx; 
\]

and they work with the functions involved in these identities as if they had
the same basic properties as the usual ones. These procedures rely on
powerful intuitions, but they are not rigorous and lead to contradictions in
the framework of ordinary calculus. And yet they must be valid in a
satisfactory calculus, based on an adequate (both for mathematics and the
experimental sciences) generalization of the Cantor-Dedekind continuum. In
the next example, we shall see how the \textit{natural indiscernible
extensions} give a positive answer to this aim, in the context of $\widehat{%
\mathbb{R}}$ .

\bigskip

\textbf{Example 5.30 \ }Let $I$ be a nonempty open interval in $\mathbb{R}$,
and let $\phi :I\rightarrow \mathbb{R}$ be a function such that $\mathit{%
\Lambda }_{\phi }=I$.\vspace{-0.2cm}

\textbf{1)} If $\phi $ is a constant function, i.e. $\phi (\xi ):=\alpha $,
for each $\xi \in I$, where $\alpha $ is a fixed real number, then, clearly,
its natural indiscernible extension is also a constant function assuming the
same value, i.e. $\hat{\phi}:m_{\approx }(I)\rightarrow \widehat{\mathbb{R}}$
is defined by%
\[
\hat{\phi}(x):=\alpha . 
\]%
\vspace{-0.35cm}$\vspace{-0.25cm}$

\textbf{2)} If $\phi $ is the inclusion function of $I$ into $\mathbb{R}$,
i.e. $\phi (\xi ):=\xi $, for each $\xi \in I$, then, since $\lambda _{\phi
}\left( \xi \right) =1$ and $\sigma x+dx=x$, for each $\xi \in \mathbb{R}$
and $x\in m_{\approx }(I)$, its natural indiscernible extension is the
inclusion function of $m_{\approx }(I)$ into $\widehat{\mathbb{R}}$, i.e. $%
\hat{\phi}:m_{\approx }(I)\rightarrow \widehat{\mathbb{R}}$ is defined by%
\[
\hat{\phi}(x):=x. 
\]

\textbf{3)} If $\phi $ is a polynomial function, i.e. $\phi (\xi ):=\alpha
_{0}{+\alpha }_{1}\xi +\ldots +\alpha _{m}\xi ^{m}$, where $\alpha
_{0},\alpha _{1}$,\dots , $\alpha _{m}$ are fixed real numbers, then, by the
previous examples, \textbf{proposition 5.29 a)}, \textbf{c)}, and
mathematical induction, its natural indiscernible extension is also a
polynomial function with the same coefficients, i.e. $\hat{\phi}:m_{\approx
}(I)\rightarrow \widehat{\mathbb{R}}$ is defined by%
\[
\hat{\phi}(x):=\alpha _{0}+\alpha _{1}x+\ldots +\alpha _{m}x^{m}. 
\]

\textbf{4)} If $\phi $ is an algebraic function, i.e. $\phi (\xi ){{:=\frac{%
\psi \left( \xi \right) }{\theta (\xi )}}}$, where $\psi :I\rightarrow 
\mathbb{R}$, $\theta :I\rightarrow \mathbb{R}$ are polynomial functions with
real coefficients, and $\theta (\xi )\neq 0$, for each $\xi \in I$, then, by
the last example and \textbf{proposition 5.29 d)}, its natural indiscernible
extension is also an algebraic function, more precisely, $\hat{\phi}%
:m_{\approx }(I)\rightarrow \widehat{\mathbb{R}}$ is defined by%
\[
\hat{\phi}(x):=\frac{\hat{\psi}\left( x\right) }{\hat{\theta}(x)}, 
\]%
where $\hat{\psi}$ and $\hat{\theta}$ are the natural indiscernible
extensions of $\psi $ and $\theta $, respectively.\vspace{-2pt}

\textbf{5)} Let $I:={\mathbb{R}}$, and let $\phi $ be the usual \textit{%
exponential function}, denoted by exp .

\vspace{-7pt}Since $\lambda _{\exp }\hspace{-2pt}\left( \xi \right) =\exp
(\xi )$, for each $\xi \in \mathbb{R}$, and $m_{\approx }(\mathbb{R})=%
\widehat{\mathbb{R}}$, the natural indiscernible extension of exp is the
function $\widehat{\exp }:\widehat{\mathbb{R}}\rightarrow \widehat{\mathbb{R}%
}$ defined by%
\[
\widehat{\exp }(x):=\exp (\sigma x)+\exp (\sigma x)dx. 
\]%
\vspace{-0.7cm}

$\widehat{\exp }{\ }$has the same basic properties as exp. For instance:

Using \textbf{proposition 5.29 f)}, we obtain:%
\[
\widehat{\exp }\hspace{-3pt}\left( \widehat{\mathbb{R}}\right) =\widehat{%
\exp }\hspace{-0.2509pc}\hspace{2pt}\left( m_{\approx }(\mathbb{R})\right)
=m_{\approx }\hspace{-3pt}\left( \exp \hspace{-2pt}\left( 
%TCIMACRO{\U{211d} }%
%BeginExpansion
\mathbb{R}
%EndExpansion
\right) \right) =m_{\approx }\hspace{-3pt}\left( 
%TCIMACRO{\U{211d} }%
%BeginExpansion
\mathbb{R}
%EndExpansion
^{+}\right) =\widehat{\mathbb{R}}^{+}. 
\]%
If $x\in \widehat{\mathbb{R}}$, then%
\[
\widehat{\exp }^{\prime }(x)=\exp (\sigma x)=\widehat{\exp }\left( \sigma
x\right) . 
\]

\vspace{-5pt}If $x_{1},x_{2}\in \widehat{\mathbb{R}}$, then (since $%
dx_{1}dx_{2}=dx_{2}dx_{1}=0$)%
\[
\widehat{\exp }(x_{1})\widehat{\exp }(x_{2})=(\exp (\sigma x_{1})+\exp
(\sigma x_{1})dx_{1})(\exp (\sigma x_{2})+\exp (\sigma x_{2})dx_{2})= 
\]%
\[
=\exp (\sigma (x_{1}+x_{2}))+\exp (\sigma (x_{1}+x_{2}))d(x_{1}+x_{2})=%
\widehat{\exp }(x_{1}+x_{2}). 
\]%
$\widehat{\exp }$ is a strictly increasing function, by \textbf{Corollary} 
\textbf{5.26 b)}, since $\widehat{\exp }^{\prime }(x)=\exp (\sigma
x)>\nolinebreak 0$,\smallskip\ for each $x\in \widehat{\mathbb{R}}$ .

And, of course,%
\[
\widehat{\exp }(0)=\exp (0)=1. 
\]

$\widehat{\exp }$ is the adequate function for the afore mentioned
considerations of physicists and engineers (as it is the case for the next
examples of natural indiscernible extensions), since it has the basic
properties of exp and is defined not only for real numbers (where it assumes
the same value as exp), but also for arguments involv\nolinebreak ing
infinitesimals. Moreover, $\widehat{\exp }(x)$ is always indiscernible from $%
\exp (\sigma x)$.

Now we may infer, rigorously, that%
\[
\widehat{\exp }(dx)=\exp (\sigma (dx))+\exp (\sigma (dx))dx=\exp (0)+\exp
(0)dx=1+dx, 
\]%
for each $x\in \widehat{\mathbb{R}}$ .

\textbf{6)} Let $I :={\mathbb{R}}$, and let $\phi$ be the usual \textit{%
natural} \textit{logarithm function}, which we denote by $\log$.

\vspace{-7pt}Since $\lambda _{\log }(\xi )=\frac{1}{\xi }$, for each $\xi
\in {\mathbb{R}}^{+}$, and $m_{\approx }({\mathbb{R}}^{+})=\widehat{\mathbb{R%
}}^{+}$, the natural indiscernible extension of $\log $ is the function $%
\widehat{\log }:\widehat{\mathbb{R}}^{+}\rightarrow \widehat{\mathbb{R}}$
defined by%
\[
\widehat{\log }(x):=\log (\sigma x)+\frac{1}{\sigma x}dx. 
\]%
By \textbf{proposition 5.29 g)}, we have: 
\[
\widehat{\log }=\widehat{\exp ^{-1}}={\widehat{\exp }}^{-1}. 
\]%
This result, in conjunction with the considerations of the previous example,
suffices to assure that $\widehat{\log }$ has the same basic properties as $%
\log $ .

And since $\widehat{\log }={\widehat{\exp }}^{-1},$ and $\widehat{\exp }\mspace{1mu}(%
\widehat{\mathbb{R}})={\widehat{\mathbb{R}}}^{+}$, we have: 
\[
\widehat{\log }\hspace{-2pt}\left( {\widehat{\mathbb{R}}}^{+}\right) =%
\widehat{\mathbb{R}}.
\]%
Clearly,\vspace{-6pt}%
\[
\widehat{\log }^{\prime }(x)=\frac{1}{\sigma x},
\]%
for each $x\in \widehat{\mathbb{R}}^{+}$.

Finally, we may infer, rigorously, that\vspace{-6pt}\vspace{1pt}%
\[
\widehat{\log }{(1+}dx)=\log (\sigma (1+dx))+\frac{1}{\sigma (1+dx)}dx=\log
(1)+dx=dx,\vspace{-2pt} 
\]%
for each $x\in \widehat{\mathbb{R}}^{+}$.

\textbf{7)} Let $I :={\mathbb{R}}$, and let $\phi$ be the usual \textit{sine
function}, denoted by sin .

\pagebreak \vspace{-7pt}Since $\lambda _{\sin }\hspace{-2pt}\left( \xi \right) =\cos
(\xi )$, for each $\xi \in \mathbb{R}$, the natural indiscernible extension
of sin is the function $\widehat{\sin }:\widehat{\mathbb{R}}\rightarrow 
\widehat{\mathbb{R}}$ defined by%
\[
\widehat{\sin }\hspace{-2pt}\left( x\right) :=\sin (\sigma x)+\cos (\sigma
x)dx. 
\]

\vspace{-0.25cm}Now let $I:={\mathbb{R}}$, and let $\phi $ be the usual 
\textit{cosine function}, denoted by cos .

\vspace{-7pt}Since $\lambda _{\cos }\hspace{-2pt}\left( \xi \right) =-\sin
(\xi )$, for each $\xi \in \mathbb{R}$, the natural indiscernible extension
of cos is the function $\widehat{\cos }:\widehat{\mathbb{R}}\rightarrow 
\widehat{\mathbb{R}}$ defined by%
\[
\widehat{\cos }\hspace{-2pt}\left( x\right) :=\cos (\sigma x)-\sin (\sigma
x)dx. 
\]%
$\widehat{\sin }$ and $\widehat{\cos }$ have the same basic properties as
sin and cos, respectively. For instance:

$\widehat{\sin }$ and $\widehat{\cos }$ have real period $2\pi $, as it is
clear from \textbf{proposition 5.29 i)}.

Using the last result and \textbf{proposition 5.29 f)}, we obtain: 
\[
\widehat{\sin }\hspace{-3pt}\left( \widehat{\mathbb{R}}\right) =\widehat{%
\sin }\hspace{-3pt}\left( \widehat{\left[ -\frac{\pi }{2},\frac{\pi }{2}%
\right] }\cup \widehat{\left[ \frac{\pi }{2},\frac{3\pi }{2}\right] }\right)
=\widehat{\sin }\hspace{-2pt}\left( \widehat{\left[ -\frac{\pi }{2},\frac{%
\pi }{2}\right] }\right) \cup \widehat{\sin }\hspace{-3pt}\left( \widehat{%
\left[ \frac{\pi }{2},\frac{3\pi }{2}\right] }\right) = 
\]%
\[
=m_{\approx }\hspace{-3pt}\left( \sin \hspace{-2pt}\left( \left] -\frac{\pi 
}{2},\frac{\pi }{2}\right[ \right) \right) \cup \mspace{1mu}\left\{ -1,1\right\} \cup\mspace{2mu}
m_{\approx }\hspace{-3pt}\left( \sin \hspace{-2pt}\left( \left] \frac{\pi }{2%
},\frac{3\pi }{2}\right[ \right) \right) \cup \left\{ -1,1\right\} =\widehat{%
\left] -1,1\right[ }\cup \mspace{1mu}\left\{ -1,1\right\} . 
\]%}
Similarly, 
\[
\widehat{\cos }\hspace{-3pt}\left( \widehat{\mathbb{R}}\right) =\widehat{%
\left] -1,1\right[ }\cup \left\{ -1,1\right\} . 
\]%
If $x\in \widehat{\mathbb{R}}$, then (since the square of an infinitesimal
is always null) 
\[
\widehat{\sin }^{2}(x)={(\sin (\sigma x)+\cos (\sigma x)dx)}^{2}=\sin
^{2}(\sigma x)+2\sin \hspace{-2pt}\left( \sigma x\right) \cos \hspace{-2pt}%
\left( \sigma x\right)dx , 
\]%
\[
\widehat{\cos }^{2}(x)={(\cos (\sigma x)-\sin (\sigma x)dx)}^{2}=\cos ^{2}%
\hspace{-2pt}\left( \sigma x\right) -2\cos \hspace{-2pt}\left( \sigma
x\right) \sin \hspace{-2pt}\left( \sigma x\right)dx . 
\]%
So 
\[
\widehat{\sin }^{2}(x)+\widehat{\cos }^{2}(x)=\sin ^{2}(\sigma x)+\cos ^{2}%
\hspace{-2pt}\left( \sigma x\right) =1. 
\]%
If $x_{1},x_{2}\in \widehat{\mathbb{R}}$, then%
\[
\widehat{\sin }(x_{1}\pm x_{2})=\sin (\sigma x_{1}\pm {\sigma x}_{2}))+\cos
(\sigma x_{1}\pm {\sigma x}_{2})d(x_{1}\pm x_{2})=\sin (\sigma x_{1})\cos
(\sigma x_{2})\pm 
\]%
\[
\pm \sin (\sigma x_{2})\cos (\sigma x_{1})+(\cos (\sigma x_{1})\cos (\sigma
x_{2})\mp \sin (\sigma x_{1})\sin (\sigma x_{2}))d\left( x_{1}\pm
x_{2}\right) . 
\]%
On the other hand (since the product of infinitesimals is always null),%
\[
\widehat{\sin }(x_{1})\widehat{\cos }(x_{2})=(\sin (\sigma x_{1})+\cos
(\sigma x_{1})dx_{1})(\cos (\sigma x_{2})-\sin (\sigma x_{2})dx_{2})= 
\]%
\[
=\sin (\sigma x_{1})\cos (\sigma x_{2})+\cos (\sigma x_{1})\cos (\sigma
x_{2})dx_{1}-\sin (\sigma x_{1})\sin (\sigma x_{2})dx_{2}, 
\]%
\[
\widehat{\sin }(x_{2})\widehat{\cos }(x_{1})=(\sin (\sigma x_{2})+\cos
(\sigma x_{2})dx_{2})(\cos (\sigma x_{1})-\sin (\sigma x_{1})dx_{1})= 
\]%
\[
=\sin (\sigma x_{2})\cos (\sigma x_{1})-\sin (\sigma x_{1})\sin (\sigma
x_{2})dx_{1}+\cos (\sigma x_{1})\cos (\sigma x_{2})dx_{2}. 
\]%
So 
\[
\widehat{\sin }(x_{1}\pm x_{2})=\widehat{\sin }(x_{1})\hspace{1pt}\widehat{%
\cos }(x_{2})\pm \widehat{\sin }\hspace{-2pt}\left( x_{2}\right) \widehat{%
\cos }\hspace{-2pt}\left( x_{1}\right) . 
\]%
In a similar manner, we could have proven that 
\[
\widehat{\cos }(x_{1}\pm x_{2})=\widehat{\cos }(x_{1})\hspace{1pt}\widehat{%
\cos }(x_{2})\mp \widehat{\sin }\hspace{-2pt}\left( x_{1}\right) \widehat{%
\sin }\hspace{-2pt}\left( x_{2}\right) . 
\]%
And we clearly have, for each $x\in \widehat{\mathbb{R}}:$%
\[
\widehat{\sin }^{\prime }(x)=\cos (\sigma x)=\widehat{\cos }\hspace{-2pt}%
\left( \sigma x\right) , 
\]%
\[
\widehat{\cos }^{\prime }(x)=-{\sin }(\sigma x)=-\widehat{\sin }\hspace{-2pt}%
\left( \sigma x\right) . 
\]

\vspace{-7pt}Finally, we may infer, rigorously, that%
\[
\widehat{\sin }(dx)=\sin (\sigma (dx))+\cos (\sigma (dx))dx=\sin (0)+\cos
(0)dx=dx, 
\]

\vspace{-7pt}for each $x\in \widehat{\mathbb{R}}.$

Similarly,%
\[
\widehat{\cos }(dx)=\cos (\sigma (dx))-\sin (\sigma (dx))dx=\cos (0)-\sin
(0)dx=1. 
\]%
\textbf{8)} Let $I:={\mathbb{R}}^{+}$, let $\alpha $ be a fixed real number,
and let $\phi $ be defined by $\phi (\xi ):=\xi ^{\alpha }$.\vspace{1pt}

\vspace{-7pt}Since $\lambda _{\phi }\hspace{-2pt}\left( \xi \right) =\alpha {%
\xi }^{\alpha -1}$, for each $\xi \in {\mathbb{R}}^{+}$, the natural
indiscernible extension of $\phi $ is \vspace{1pt}\newline
the function $\hat{\phi}:\widehat{\mathbb{R}}^{+}\rightarrow \widehat{%
\mathbb{R}}$ defined by%
\[
\hat{\phi}\hspace{-1pt}\left( x\right) :=({\sigma x)}^{\alpha }+\alpha {%
(\sigma x)}^{\alpha -1}dx. 
\]%
Clearly, for each $x\in \widehat{\mathbb{R}}^{+}$:%
\[
\hat{\phi}^{\prime }(x)=\alpha (\sigma x)^{\alpha -1}. 
\]%
If we denote $\hat{\phi}\hspace{-1pt}\left( x\right) $ by $x^{\alpha },{\ }${%
then}%
\[
x^{\alpha }:=({\sigma x)}^{\alpha }+\alpha {(\sigma x)}^{\alpha -1}dx, 
\]%
\[
\left( x^{\alpha }\right) ^{\prime }=\alpha {(\sigma x)}^{\alpha -1}, 
\]%
for each $x\in \widehat{\mathbb{R}}^{+}$.\vspace{-1pt}

Trivially, $\hat{\phi}\hspace{-3pt}\left( {\widehat{\mathbb{R}}}^{{+}%
}\right) =\left\{ 1\right\} ,$ when $\alpha =0.$ If $\alpha \neq 0,$ then we
obtain, using \textbf{proposition 5.29 f)}: 
\[
\hat{\phi}\hspace{-3pt}\left( {\widehat{\mathbb{R}}}^{{+}}\right) =\hat{\phi}%
\hspace{-3pt}\left( m_{\approx }\hspace{-3pt}\left( {\mathbb{R}}^{+}\right)
\right) =m_{\approx }\hspace{-3pt}\left( \phi \hspace{-3pt}\left( {\mathbb{R}%
}^{+}\right) \right) =m_{\approx }({\mathbb{R}}^{+})={\widehat{\mathbb{R}}}^{%
{+}}.\vspace{-4pt} 
\]%
As $\phi (\xi )=\exp (\alpha \log (\xi ))$, for each $\xi \in {\mathbb{R}}%
^{+}$, we obtain, using the examples \textbf{1)}, \textbf{5)}, \textbf{6)},
and \textbf{proposition 5.29 c)}, \textbf{e)} :%
\[
x^{\alpha }=\hat{\phi}\hspace{-1pt}\left( x\right) =\widehat{\exp }(\alpha 
\widehat{\log }(x)),\text{ for each }x\in \widehat{\mathbb{R}}^{+}. 
\]%
\vspace{-4pt}

Finally, we may infer, with complete rigour, that%
\[
(1+dx)^{\alpha }=(\sigma (1+dx))^{\alpha }+\alpha (\sigma (1+dx))^{\alpha
-1}dx=1^{\alpha }+\alpha .1^{\alpha -1}dx=1+\alpha dx. 
\]%
\textbf{9)} Let $I:=\mathbb{R}$, let $\alpha $ be a fixed positive real
number, and let $\phi $ be defined by $\phi (\xi ):=\nolinebreak \alpha
^{\xi }$.

\vspace{-7pt}Since $\lambda _{\phi }\hspace{-2pt}\left( \xi \right) =\alpha
^{\xi }\log (\alpha )$, for each $\xi \in {\mathbb{R}}$, the natural
indiscernible extension of $\phi $ is the function $\widehat{\phi }:\widehat{%
\mathbb{R}}\rightarrow \widehat{\mathbb{R}}$ defined by%
\[
\hat{\phi}\left( x\right) :=\alpha ^{\sigma x}+\alpha ^{\sigma x}\log
(\alpha )dx. 
\]%
Clearly, for each $x\in \widehat{\mathbb{R}}$:%
\[
\hat{\phi}^{\prime }(x)=\alpha ^{\sigma x}\log (\alpha ). 
\]%
If we denote $\hat{\phi}\left( x\right) $ by $\alpha ^{x},{\ }${then we h}%
ave, for each $x\in \widehat{\mathbb{R}}:$%
\[
\alpha ^{x}:=\alpha ^{\sigma x}+\alpha ^{\sigma x}\log (\alpha )dx, 
\]%
\[
\left( \alpha ^{x}\right) ^{\prime }=\alpha ^{\sigma x}\log (\alpha ). 
\]%
So, if $e$ is \textit{Euler's number}, then%
\[
e^{x}=e^{\sigma x}+e^{\sigma x}\log (e)dx=\exp (\sigma x)+\exp (\sigma x)dx=%
\widehat{\exp }(x), 
\]%
for each $x\in \widehat{\mathbb{R}}$.\vspace{-2pt}

Trivially, $\hat{\phi}\left( \widehat{\mathbb{R}}\right) =\left\{ 1\right\}
, $ when $\alpha =1.$ If $\alpha \neq $ $1$, then we obtain, using \textbf{%
proposition 5.29 f)}: 
\[
\hat{\phi}\hspace{-2pt}\left( \widehat{\mathbb{R}}\right) =\hat{\phi}\hspace{%
-1pt}\left( m_{\approx }\hspace{-2pt}\left( {\mathbb{R}}\right) \right)
=m_{\approx }(\phi \hspace{-1pt}\left( {\mathbb{R}}\right) )=m_{\approx }({%
\mathbb{R}}^{+})={\widehat{\mathbb{R}}}^{{+}}. 
\]

\bigskip

The next definition introduces the concepts of \textit{mth natural
indiscernible extension} and \textit{mth derivative function}, for $m\in 
\mathbb{N}$ .

\bigskip

\textbf{Definition 5.31 \ }Let $I$ be a nonempty open interval in $\mathbb{R}
$, and let $\phi :I\rightarrow \mathbb{R}$ be a function such that $\mathit{%
\Lambda }_{\phi }=I$.

\vspace{-6pt}The functions $\hat{\phi}:m_{\approx }(I)\rightarrow \widehat{%
\mathbb{R}},$ $\hat{\phi}^{\prime }:m_{\approx }(I)\rightarrow \widehat{%
\mathbb{R}}$ defined by%
\[
\hat{\phi}(x):=\phi (\sigma x)+\lambda _{\phi }(\sigma x)dx, 
\]%
\[
\hat{\phi}^{\prime }(x):=\lambda _{\phi }(\sigma x), 
\]

\vspace{-7pt}will be called the\textit{\ first natural indiscernible
extension }of $\phi ,$ and the\textit{\ first derivative function }of $\hat{%
\phi}$, respectively. So the \textit{first natural indiscernible extension }%
of $\phi $ is, in fact, its \textit{natural indiscernible extension}, and,
most conveniently, the value of the \textit{first derivative function }of $%
\hat{\phi}$ at $\xi _{0}\in I$ is its \textit{derivative} at this point (see 
\textbf{definition 5.28 }and \textbf{definition 5.18}, respectively).\vspace{%
2pt}

\vspace{-6pt}If $\mathit{\Lambda }_{\lambda _{\phi }}=I$, then the functions 
$\hat{\phi}^{\left[ 2\right] }:m_{\approx }(I)\rightarrow \widehat{\mathbb{R}%
},$ $\hat{\phi}^{\prime \prime }:m_{\approx }(I)\rightarrow \widehat{\mathbb{%
R}}$ defined by%
\[
\hat{\phi}^{[2]}(x):=\lambda _{\phi }(\sigma x)+\lambda _{\lambda _{\phi
}}(\sigma x)dx, 
\]%
\[
\hat{\phi}^{\prime \prime }(x):=\lambda _{\lambda _{\phi }}(\sigma x), 
\]

\vspace{-6pt}will be called the\textit{\ second natural indiscernible
extension }of $\phi ,$ and the\textit{\ second derivative function }of $\hat{%
\phi}$, respectively.

\vspace{-6pt}If $\mathit{\Lambda }_{\lambda _{\lambda _{\phi }}}=I$, then
the functions $\hat{\phi}^{\left[ 3\right] }:m_{\approx }(I)\rightarrow 
\widehat{\mathbb{R}},$ $\hat{\phi}^{\prime \prime \prime }:m_{\approx
}(I)\rightarrow \widehat{\mathbb{R}}$ defined by%
\[
\hat{\phi}^{[3]}(x):=\lambda _{\lambda _{\phi }}(\sigma x)+\lambda _{\lambda
_{\lambda _{\phi }}}(\sigma x)dx, 
\]%
\[
\hat{\phi}^{\prime \prime \prime }(x):=\lambda _{\lambda _{\lambda _{\phi
}}}(\sigma x), 
\]

\vspace{-6pt}will be called the\textit{\ third natural indiscernible
extension }of $\phi ,$ and the\textit{\ third derivative function }of $\hat{%
\phi}$, respectively.

\vspace{-6pt}For the sake of uniformity, we also denote $\hat{\phi},$ $\hat{%
\phi}^{\prime },$ $\hat{\phi}^{\prime \prime },$ $\hat{\phi}^{\prime \prime
\prime }$ by $\hat{\phi}^{[1]},$ $\hat{\phi}^{(1)},$ $\hat{\phi}^{(2)},$ $%
\hat{\phi}^{(3)}$, respectively.

We define in a similar manner the\textit{\ fourth natural indiscernible
extension }of $\phi $ and the\textit{\ fourth derivative function }of $\hat{%
\phi}$, denoted by $\hat{\phi}^{[4]}$ and $\hat{\phi}^{\left( 4\right) }$,
respectively,\dots ; and if $m\in \mathbb{N}$, then we denote by $\hat{\phi}%
^{[m]}$ and $\hat{\phi}^{(m)}$ the\textit{\ mth natural indiscernible
extension }of $\phi $ and the\textit{\ mth derivative function }of $\hat{\phi%
}$, when such functions exist .

\bigskip

\textbf{Notation } Let $m \in \mathbb{N}$.

\vspace{-6pt}Under the conditions and with the notation of \textbf{%
definition 5.31}, $\lambda _{\phi }^{(m)}$ will indicate that the symbol $%
\lambda $ appears $m$ times. For example: 
\[
\lambda _{\phi }^{(1)}:=\lambda _{\phi }, 
\]%
\[
\lambda _{\phi }^{(2)}:=\lambda _{\lambda _{\phi }}, 
\]%
\[
\lambda _{\phi }^{(3)}:=\lambda _{\lambda _{\lambda _{\phi }}}. 
\]%
And if we define $\lambda _{\phi }^{(0)}:=\phi $, then we have, for each $%
x\in m_{\approx }(I)$, and $m\in \mathbb{N}$:%
\[
\hat{\phi}^{[m]}(x)=\lambda _{\phi }^{(m-1)}(\sigma x)+\lambda _{\phi
}^{(m)}(\sigma x)dx. 
\]%
Since $\lambda _{\phi }^{(0)}:=\phi ,$ it is \og natural\fg ~ to introduce the function $\hat{\phi}^{(0)}:m_{\approx
}(I)\rightarrow \widehat{%
%TCIMACRO{\U{211d} }%
%BeginExpansion
\mathbb{R}
%EndExpansion
},$ defined by $\hat{\phi}^{(0)}\left( x\right) :=\phi \left( \sigma
x\right) =\hat{\phi}\left( \sigma x\right) .$

\bigskip

Clearly:

\bigskip

\textbf{Proposition 5.32 \ }Let $m\in \mathbb{N}$. Then:\vspace{-6pt}

\textbf{a)} $\hat{\phi}^{[m]}$ is the (first) natural indiscernible
extension of $\lambda _{\phi }^{(m-1)}$, i.e. $\hat{\phi}^{[m]}=\widehat{%
\lambda _{\phi }^{(m-1)}}$.\vspace{-3pt}

\textbf{b)} $\hat{\phi}^{(m)}=\sigma \circ \hat{\phi}^{[m+1]}$ (where $%
\sigma :\widehat{\mathbb{R}}\rightarrow \widehat{\mathbb{R}}$ is the \textit{%
shadow function}, i.e. $\sigma (x):=\sigma x,$ for each $x\in \widehat{%
%TCIMACRO{\U{211d} }%
%BeginExpansion
\mathbb{R}
%EndExpansion
}$) .

\bigskip

\textbf{Remark 5.33 \ }Let $m\in \mathbb{N}$.\vspace{-3pt}

If $\hat{\phi}^{[m+1]}$ and $\hat{\phi}^{(m+1)}$exist, it is important to
notice that $\hat{\phi}^{(m+1)}$ is the derivative function of $\hat{\phi}%
^{[m+1]}$, and not the derivative function of $\hat{\phi}^{(m)}$. This is
not surprising since $\hat{\phi}^{[m+1]}$ is the (first) natural
indiscernible extension of $\lambda _{\phi }^{(m)}$, and $\lambda _{\phi
}^{(m)}$ is, in fact, \textit{the usual} \textit{mth derivative function of} 
$\phi $.\vspace{-7pt}

In blunt terms, the rule (valid for the\textit{\ derivative at a point} or
the\textit{\ derivative function}) is%
\[
\text{\textit{The derivative is always associated with an indiscernible
extension}.} 
\]%
Finally, it is important to realize that the range of $\hat{\phi}^{(m)}$ is
always a subset of $\mathbb{R}$, although its codomain is $\widehat{\mathbb{R%
}}$ .\vspace{-2pt}

\bigskip

\textbf{Example 5.34 \ 1)} Let $\phi :\mathbb{R}\rightarrow \mathbb{R}$ be
the function defined by $\phi (\xi ):=\xi ^{2}$. Then $\mathit{\Lambda }%
_{\lambda _{\phi }^{(m)}}=\mathbb{R}$, for each $m\in \mathbb{N}_{0}$ (where 
$\mathbb{N}_{0}:=\mathbb{N}\cup \left\{ 0\right\} )$, and we have, for each $%
\xi \in {\mathbb{R}}$ : 
\[
\lambda _{\phi }^{(0)}(\xi {)=}\phi (\xi )=\xi ^{2}, 
\]%
\[
\lambda _{\phi }^{(1)}(\xi )=\lambda _{\phi }(\xi )=2\xi , 
\]%
\[
\lambda _{\phi }^{(2)}(\xi )=2, 
\]%
\[
\lambda _{\phi }^{(m)}(\xi )=0,\text{ for }m\geq 3. 
\]%
Then, for each $x\in \widehat{\mathbb{R}}$, and $m\in \mathbb{N}$:%
\[
\hat{\phi}^{[m]}(x)=\lambda _{\phi }^{(m-1)}(\sigma x)+\lambda _{\phi
}^{(m)}(\sigma x)dx=\left\{ 
\begin{array}{c}
(\sigma x)^{2}+2\left( \sigma x\right) dx=x^{2}{,}\text{ {if \ }}m=1\  \\ 
2\sigma x+{2}dx=2x{,}\text{ {if \ }}m=2\  \\ 
2{,}\text{ {if \ }}m=3\  \\ 
0{,}\text{ {if \ }}m\geq 4%
\end{array}%
\right. , 
\]%
as it should be, according to \textbf{example 5.30 3)}, and \textbf{%
proposition 5.32 a)}.

For each $x\in \widehat{\mathbb{R}}$, and $m\in \mathbb{N}$, we have:%
\[
\hat{\phi}^{(m)}(x)=\lambda _{\phi }^{(m)}(\sigma x)=\left\{ 
\begin{array}{c}
2\sigma x{,}\text{ {if \ }}m=1\  \\ 
2{,}\text{ {if \ }}m=2\  \\ 
0{,}\text{ {if \ }}m\geq 3%
\end{array}%
\right. \text{,} 
\]%
as it should be, according to the results we obtained for $\hat{\phi}^{[m]}$%
, and \textbf{proposition 5.32 b)}.

We could have written the last identities more synthetically as%
\[
(x^{2})^{\prime }=2\sigma x, 
\]%
\[
(x^{2})^{\prime \prime }=2, 
\]%
\[
(x^{2})^{(m)}=0,\text{ for }m\geq 3. 
\]%
\textbf{2)} Let $\phi :\mathbb{R}\rightarrow \mathbb{R}$ be the function
defined by $\phi (\xi ):=\exp (\xi )$. Then $\mathit{\Lambda }_{\lambda
_{\phi }^{(m)}}=\mathbb{R}$, for each $m\in {\mathbb{N}}_{0}$ , and we have:%
\[
\lambda _{\phi }^{(m)}(\xi )=\exp (\xi ),\text{ for each }\xi \in {\mathbb{R}%
},\text{ and }m\in {\mathbb{N}}_{0}. 
\]%
Then, for each $x\in \widehat{\mathbb{R}}$, and $m\in \mathbb{N}:\smallskip $%
\[
\hat{\phi}^{[m]}(x)=\lambda _{\phi }^{(m-1)}(\sigma x)+\lambda _{\phi
}^{(m)}(\sigma x)dx=\exp (\sigma x)+\exp (\sigma x)dx=\widehat{\exp }(x), 
\]%
\[
\hat{\phi}^{(m)}(x)=\lambda _{\phi }^{(m)}(\sigma x)=\exp (\sigma x)=%
\widehat{\exp }(\sigma x). 
\]%
{More }{synthetically:}%
\[
\widehat{\exp }^{[m]}(x)=\widehat{\exp }(x), 
\]%
\[
\widehat{\exp }^{(m)}(x)=\widehat{\exp }(\sigma x); 
\]

\vspace{-6pt}for each $x\in \widehat{\mathbb{R}}$, and $m\in \mathbb{N}$.

\textbf{3)} Let $\phi :\mathbb{R}\rightarrow \mathbb{R}$ be the function
defined by $\phi (\xi ):=\sin (\xi )$. Then $\mathit{\Lambda }_{\lambda
_{\phi }^{(m)}}=\mathbb{R}$, for each $m\in {\mathbb{N}}_{0}$, and we have: 
\[
\lambda _{\phi }^{\left( m\right) }\left( \xi \right) =\left\{ 
\begin{array}{c}
(-1)^{\frac{m-1}{2}}\cos \left( \xi \right) {,}\text{ {if \ }}m{\ }\text{{is
odd}}\  \\ 
(-1)^{\frac{m}{2}}\sin \left( \xi \right) {,}\text{ {if }}m{\ }\text{{is even%
}}%
\end{array}%
\right. . 
\]%
Then, for each $x\in \widehat{\mathbb{R}}$, and $m\in \mathbb{N}$:%
\[
\hat{\phi}^{[m]}(x)=\lambda _{\phi }^{(m-1)}(\sigma x)+\lambda _{\phi
}^{(m)}(\sigma x)dx= 
\]%
\[
=\left\{ 
\begin{array}{c}
(-1)^{\frac{m-1}{2}}\sin \left( \sigma x\right) +\left( -1\right) ^{\frac{m-1%
}{2}}\cos \left( \sigma x\right) dx{,}\text{ {if }}m\text{ {is odd}}\  \\ 
(-1)^{\frac{m-2}{2}}\cos \left( \sigma x\right) +\left( -1\right) ^{\frac{m}{%
2}}\sin \left( \sigma x\right) dx{,}\text{ {if }}m\text{ {is even}}%
\end{array}%
\right. = 
\]%
\[
=\left\{ 
\begin{array}{c}
(-1)^{\frac{m-1}{2}}\left(\sin \left( \sigma x\right) +\cos
\left( \sigma x\right)\right) dx{,}\text{ {if }}m\text{ {is odd}}\  \\ 
(-1)^{\frac{m-2}{2}}(\cos \left( \sigma x\right) -\sin \left( \sigma x\right)
dx){,}\text{ {if }}m\text{ {is even}}%
\end{array}%
\right. = 
\]%
\[
=\left\{ 
\begin{array}{c}
(-1)^{\frac{m-1}{2}}\widehat{\sin }(x){,}\text{ {if }}m\text{ {is odd}}\  \\ 
(-1)^{\frac{m-2}{2}}\widehat{\cos }(x){,}\text{ {if }}m\text{ {is even}}%
\end{array}%
\right. . 
\]%
For each $x\in \widehat{\mathbb{R}}$, and $m\in \mathbb{N}$, we have:%
\[
\hat{\phi}^{(m)}(x)=\lambda _{\phi }^{(m)}(\sigma x)= 
\]%
\[
=\left\{ 
\begin{array}{c}
(-1)^{\frac{m-1}{2}}\cos \left( \sigma x\right) {,}\text{ {if }}m{\ }\text{{%
is odd}}\  \\ 
(-1)^{\frac{m}{2}}\sin \left( \sigma x\right) {,}\text{ {if }}m{\ }\text{{is
even}}%
\end{array}%
\right. = 
\]%
\[
=\left\{ 
\begin{array}{c}
(-1)^{\frac{m-1}{2}}\widehat{\cos }(\sigma x){,}\text{ {if }}m\text{ {is odd}%
}\  \\ 
(-1)^{\frac{m}{2}}\widehat{\sin }(\sigma x){,}\text{ {if }}m\text{ {is even}}%
\end{array}%
\right. . 
\]%
More synthetically, we have, for each $x\in \widehat{\mathbb{R}}$, and $m\in 
\mathbb{N}:$%
\[
\widehat{\sin }^{[m]}(x)=\left\{ 
\begin{array}{c}
(-1)^{\frac{m-1}{2}}\widehat{\sin }(x){,}\text{ {if }}m\text{ {is odd}}\  \\ 
(-1)^{\frac{m-2}{2}}\widehat{\cos }(x){,}\text{ {if }}m\text{ {is even}}%
\end{array}%
\right. , 
\]%
\[
\widehat{\sin }^{(m)}(x)=\left\{ 
\begin{array}{c}
(-1)^{\frac{m-1}{2}}\widehat{\cos }(\sigma x){,}\text{ {if }}m\text{ {is odd}%
}\  \\ 
(-1)^{\frac{m}{2}}\widehat{\sin }(\sigma x){,}\text{ {if }}m\text{ {is even}}%
\}%
\end{array}%
\right. . 
\]%
For the \textit{cosine function}, we have $\mathit{\Lambda }_{\lambda _{\cos
}^{(m)}}=\mathbb{R}$, and $\lambda _{\cos }^{(m)}=\lambda _{\lambda _{\sin
}}^{(m)}=\lambda _{\sin }^{(m+1)},$ for each $m\in \mathbb{N}_{0}.$ Then,
for each $x\in \widehat{\mathbb{R}},$ and $m\in \mathbb{N}$,%
\[
\widehat{\cos }^{[m]}(x)=\lambda _{\cos }^{(m-1)}(\sigma x)+\lambda _{\cos
}^{(m)}(\sigma x)dx= 
\]%
\[
=\lambda _{\sin }^{(m)}(\sigma x)+\lambda _{\sin }^{(m+1)}(\sigma x)dx=%
\widehat{\sin }^{(m)}(\sigma x)+\widehat{\sin }^{(m+1)}(\sigma x)dx= 
\]%
\[
=\left\{ 
\begin{array}{c}
(-1)^{\frac{m-1}{2}}(\widehat{\cos }(\sigma x)-\widehat{\sin }(\sigma x)dx){,%
}\text{ {if }}m\text{ {is odd}}\  \\ 
(-1)^{\frac{m}{2}}(\widehat{\sin }(\sigma x)+\widehat{\cos }(\sigma x)dx){,}%
\text{ {if }}m\text{ {is even}}%
\end{array}%
\right. = 
\]%
\[
=\left\{ 
\begin{array}{c}
(-1)^{\frac{m-1}{2}}\widehat{\cos }(x){,}\text{ {if }}m\text{ {is odd}}\  \\ 
(-1)^{\frac{m}{2}}\widehat{\sin }(x){,}\text{ {if }}m\text{ {is even}}%
\end{array}%
\right. ; 
\]%
\[
\widehat{\cos }^{(m)}(x)=\lambda _{\cos }^{(m)}(\sigma x)={\lambda }_{\sin
}^{(m+1)}(\sigma x)={\widehat{\sin }}^{(m+1)}(\sigma x)=\left\{ 
\begin{array}{c}
(-1)^{\frac{m+1}{2}}\widehat{\sin }(\sigma x){,}\text{ {if }}m\text{ {is odd}%
}\  \\ 
(-1)^{\frac{m}{2}}\widehat{\cos }(\sigma x){,}\text{ {if }}m\text{ {is even}}%
\end{array}%
\right. . 
\]

\bigskip

We close this section with \textbf{Taylor's Theorem}.

\bigskip

\textbf{Theorem 5.35 (Taylor's Theorem) \ }Let \textit{I} be an open
interval in $\mathbb{R}$, let $\xi _{0}\in I$ and $m\in {\mathbb{N}}_{0}$,
and let $\phi :I\rightarrow \mathbb{R}$ be a function such that $\mathit{%
\Lambda }_{\lambda _{\phi }^{(k)}}=I$, for each $0\leq k\leq m$.

\vspace{-9pt}Then for each $x\in m_{\approx }\hspace{-3pt}\left( I\backslash 
\hspace{-2pt}\left\{ \xi _{0}\right\} \right) $there exists a real number $\theta \in \left] 0,1\right[ $
such that:%
\[
\hat{\phi}(x)\approx \sum_{k=0}^{m}\frac{\hat{\phi}^{\left( k\right) }%
\hspace{-3pt}\left( \xi _{0}\right) }{k!}\left( x-\xi _{0}\right)
^{k}+\frac{{(\ x-\xi _{0})}^{m+1}}{(m+1)!}\hat{\phi}^{(m+1)}\hspace{-1pt%
}(\xi _{0}{+}\theta (\ x-\xi _{0}{))}. 
\]

\bigskip

\textbf{Proof \ }By the usual \textbf{Taylor's Theorem} with the Lagrange
form of the remainder, for each $x\in m_{\approx }\left( I\backslash 
\hspace{-2pt}\left\{ \xi _{0}\right\} \right) $ there exists a real number $\theta \in \left] 0,1%
\right[ $ such that we have:%
\[
\hat{\phi}(x)\approx \phi (\sigma x)=\sum_{k=0}^{m}\frac{\lambda _{\phi
}^{(k)}\hspace{-3pt}\left( \xi _{0}\right) }{k!}\left( \sigma x-\xi
_{0}\right) ^{k}+\frac{{(\sigma x-\xi _{0})}^{m+1}}{(m+1)!}\lambda _{\phi
}^{(m+1)}\hspace{-1pt}(\xi _{0}{+}\theta (\sigma x-\xi _{0}{))}\approx 
\]

\[
\approx\sum_{k=0}^{m}{\ \frac{\hat{\phi}^{\left( k\right) }\hspace{-3pt}\left( \xi
_{0}\right) }{k!}}{\left( \ x-\xi _{0}\right) }^{k}+\frac{{(\
x-\xi _{0})}^{m+1}}{(m+1)!}\hat{\phi}^{(m+1)}\hspace{-1pt}(\xi _{0}{+}\theta
(\ x-\xi _{0}{))}.\text{ }\blacksquare 
\]

\section{The Differential Treatment of Singularities (two examples)}

For each $\xi _{0}\in \mathbb{R},$ $m_{\approx }(\xi _{0})$ has three
remarkable features:

\textbf{(i)} It has the same cardinality as $\widehat{\mathbb{R}}$, since
(see \textbf{proposition 4.1} and its proof)%
\[
\left\vert m_{\approx }\hspace{-3pt}\left( \xi _{0}\right) \right\vert
=\left\vert \widehat{\mathbb{R}}\right\vert =2^{\aleph _{0}}. 
\]

\vspace{-6pt}\textbf{(ii)} It is a closed interval in $\widehat{\mathbb{R}}$
with length $0$, since (see \textbf{proposition 5.14}) 
\[
m_{\approx }\hspace{-3pt}\left( \xi _{0}\right) =\widehat{\left[ \xi
_{0},\xi _{0}\right] }, 
\]%
\[
l(\widehat{\left[ \xi _{0},\xi _{0}\right] })=0. 
\]%
In this sense, $m_{\approx }\hspace{-3pt}\left( \xi _{0}\right) $ may be
viewed as a \textit{tiny} subset of $\widehat{\mathbb{R}}$.

\textbf{(iii)} It has a geometric structure, since (see \textbf{proposition
4.2 b)})%
\[
m_{\approx }\hspace{-3pt}\left( \xi _{0}\right) \text{ is an
infinite-dimensional real affine space.} 
\]

We may use \textbf{(ii}) to obtain immediately:

$($\textbf{ii}$\mathbf{^{\prime }})$ If $\xi _{0}\in \mathbb{R}$, then 
\[
m_{\approx }\hspace{-3pt}\left( \xi _{0}\right) \cap \mathbb{R}=\left\{ \xi
_{0}\right\} . 
\]%
\textbf{(i}) and \textbf{(iii)} express properties of $m_{\approx }\hspace{%
-3pt}\left( \xi _{0}\right) $ that are shared with the \textit{entire}
generalized real continuum (the fact that $\widehat{\mathbb{R}}$ is an
infinite-dimensional real affine space may be easily derived from \textbf{%
proposition 2.3 a)} and \textbf{(iii)}). Nevertheless $m_{\approx }\hspace{%
-3pt}\left( \xi _{0}\right) $ is a \textit{tiny} subset of $\widehat{\mathbb{%
R}}$, by \textbf{(ii)}. This \textit{global-local} nature of $m_{\approx }%
\hspace{-3pt}\left( \xi _{0}\right) $ is the source of its usefulness for
the \textit{differential calculus}. In the next two examples, we apply this
dual nature to the differential treatment of a singularity, using $($\textbf{%
ii}$^{\prime })$ and \textbf{(iii)}.

\bigskip

\textbf{Example 6.1 \ 1)} Consider, in $\mathbb{R}$, the differential
equation: 
\begin{equation}
\xi ^{\prime }(\tau )=\left\{ 
\begin{array}{c}
-1,{\ }\text{{if \ }}\tau <0 \\ 
1,{\ }\text{{if \ }}\tau =0 \\ 
1,{\ }\text{{if \ }}\tau >0%
\end{array}%
\right. .  \label{GrindEQ__2_}
\end{equation}%
Equation (2) has no solution on any open interval $I$ in $\mathbb{R}$ such
that $0\in I$, since if such a solution $\xi :I\rightarrow {\mathbb{R}}$
existed, then $\xi ^{\prime }$ would not satisfy the \textit{intermediate
value} \textit{property} on $I$ $\left[ \text{{see }{Fig}. 1}\right] $,
violating \textbf{Darboux's Theorem}.

% \begin{figure}[ht!]
% \begin{center}
% \includegraphics[clip, trim=10cm 3cm 16cm 10cm]{./image1.eps}
% \end{center}
% \end{figure}

 \begin{center}
 \begin{tikzpicture} 
 \begin{axis}[ 
 axis x line=middle,
 axis y line=middle,
 ylabel={$\xi'$},
 xlabel={$\tau$},
 axis line style={->},
 xmin=-1,xmax=1,
 ymin=-2,ymax=2,
 xtick=\empty,
 ytick=\empty,
 %ytick={1},ytick pos=left,
 %ytick={-1}, ytick pos=right,
 tickwidth=0.3cm,
 enlargelimits=false
 ]
 \addplot[black,mark=none,
 domain=0:5,samples=20,line width=2pt]
 {1} node [pos=0,left] {$1$};
 \addplot[fill=white,
 only marks,mark=*, mark size=3pt,
 every mark/.append style={rotate=-90}]
 coordinates{(0,-1)}  node [pos=0,right] {$ -1$};
 \addplot[black,mark=none,
 domain=-5:0,samples=20,line width=2pt]
 {-1};
 \end{axis}
 \end{tikzpicture}
 \end{center}
Fig. 1: $\xi ^{\prime }$ would not satisfy the \textit{intermediate value} 
\textit{property} on $I$, for any open interval $I$ in $\mathbb{R}$ such
that $0\in I$ .

Now consider the corresponding differential equation in $\widehat{\mathbb{R}}%
:$%
\begin{equation}
x^{\prime }(t)=\left\{ 
\begin{array}{c}
-1,{\ }\text{{if \ }}t<0 \\ 
1,{\ }\text{\ {if \ }}t\in m_{\approx }(0) \\ 
1,{\ }\text{{if \ }}t>0%
\end{array}%
\right. .  \label{GrindEQ__3_}
\end{equation}%
Equation (3) has an infinity of solutions on $\widehat{\mathbb{R}}$; for
instance, one solution is $\left[ \text{{see }{Fig}. 2}\right] $%
\[
x(t):=\left\{ 
\begin{array}{c}
-t,{\ }\text{{if \ }}t<0 \\ 
dt,{\ }\text{{if \ }}t\in m_{\approx }{(0)} \\ 
t,{\ }\text{{if \ }}t>0%
\end{array}%
\right. =\left\{ 
\begin{array}{c}
-t,{\ }\text{{if \ }}t<0 \\ 
t,{\ }\text{{if \ }}t\in m_{\approx }{(0)} \\ 
t,{\ }\text{{if \ }}t>0%
\end{array}%
\right. . 
\]

% \begin{figure}[ht!]
% \begin{center}
% \includegraphics*[clip,trim=5.60in 3.65in 8.2in 4.60in]{./image2.eps}
% \end{center}
% \end{figure}

 \begin{center}
 \begin{tikzpicture} 
 \begin{axis}[ 
 axis x line=middle,
 axis y line=middle,
 ylabel={$x$},
 xlabel={$t$},
 axis line style={->},
 xmin=-2,xmax=2,
 ymin=-1,ymax=3,
 xtick=\empty,
 ytick=\empty,
 tickwidth=0.3cm,
 enlargelimits=false
 ]
 \addplot[black,mark=none,
 domain=-0.5:3,samples=20,line width=2.5pt]
 {x};
 \addplot[black,mark=none,
 domain=-3:-0.5,samples=20,line width=2.5pt]
 {-x};
 \end{axis}
 \end{tikzpicture}
 \end{center}
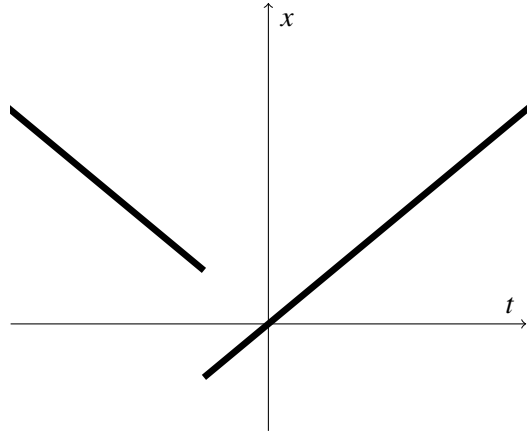
Fig. 2: A solution $x:\widehat{\mathbb{R}}\rightarrow \widehat{\mathbb{R}}$
of the differential equation (3)

\bigskip

\vspace{-6pt}Notice that $x:\widehat{\mathbb{R}}\rightarrow \widehat{\mathbb{%
R}}$ is an indiscernible extension of $\xi :\mathbb{R}\rightarrow \mathbb{R}$
, defined by $\xi (\tau ):=\left\vert \tau \right\vert $.

\textbf{2)} Consider, in $\mathbb{R}$, the differential equation: 
\begin{equation}
\xi ^{\prime }(\tau )={\delta }_{0}(\left\{ \tau \right\} )=\left\{ 
\begin{array}{c}
1,{\ }\text{{if \ }}\tau =0 \\ 
0,{\ }\text{{if \ }}\tau \neq 0%
\end{array}%
\right. .  \label{GrindEQ__4_}
\end{equation}%
By \textbf{Darboux's Theorem}, equation $\left( 4\right) $ has no solution
on any open interval $I$ in $\mathbb{R}$ such that $0\in I$ $\left[ \text{{%
see }{Fig}. 3}\right] .$

% \begin{figure}[ht!]
% \begin{center}
% \includegraphics*[clip,trim=4.83in 1.27in 7.59in 4.03in]{./image3.eps}
% \end{center}
% \end{figure}

 \begin{center}
 \begin{tikzpicture} 
 \begin{axis}[ 
 axis x line=middle,
 axis y line=middle,
 ylabel={$\xi'$},
 xlabel={$\tau$},
 axis line style={->},
 xmin=-2,xmax=2,
 ymin=-1,ymax=1.5,
 xtick=\empty,
 ytick={1},
 tickwidth=0.3cm,
 enlargelimits=false
 ]
 \addplot[black,mark=none,
 domain=-2:1.9,samples=20,line width=2.5pt]
 {0};
 \addplot[fill=white,
 only marks,mark=*, mark size=3pt,
 every mark/.append style={rotate=-90}]
 coordinates{(0,0)};
 \addplot[fill=black,
 only marks,mark=*, mark size=3pt,
 every mark/.append style={rotate=-90}]
 coordinates{(0,1)};
 \end{axis}
 \end{tikzpicture}
 \end{center}
Fig. 3: $\xi ^{\prime }$ would not satisfy the \textit{intermediate value} 
\textit{property} on $I$, for any open interval $I$ in $\mathbb{R}$ such
that $0\in I$ .

Now consider the corresponding differential equation in $\widehat{\mathbb{R}}%
:$ 
\begin{equation}
x^{\prime }(t)=\left\{ 
\begin{array}{c}
1,{\ }\text{{if \ }}t\in m_{\approx }(0) \\ 
0,{\ }\text{{if \ }}t\notin m_{\approx }(0)%
\end{array}%
\right. .  \label{GrindEQ__5_}
\end{equation}%
Equation $\left( 5\right) $ has an infinity of solutions on $\widehat{%
\mathbb{R}}$; for instance, one solution is $\left[ \text{{see }{Fig}. 4}%
\right] $%
\[
x(t):=\left\{ 
\begin{array}{c}
1+dt,{\ }\text{{if \ }}t\in m_{\approx }(0) \\ 
1,{\ }\text{{if \ }}t>0 \\ 
0,{\ }\text{{if \ }}t<0%
\end{array}%
\right. . 
\]

% \begin{figure}[ht!]
% \begin{center}
% \includegraphics*[clip, trim=5.05in 2.01in 7.31in 4.62in]{./image4.eps}
% \end{center}
% \end{figure}

 \begin{center}
 \begin{tikzpicture} 
 \begin{axis}[ 
 axis x line=middle,
 axis y line=middle,
 ylabel={$x$},
 xlabel={$t$},
 axis line style={->},
 xmin=-2,xmax=2,
 ymin=-1,ymax=2,
 xtick=\empty,
 ytick={1},
 tickwidth=0.3cm,
 enlargelimits=false
 ]
 \addplot[black,mark=none,
 domain=-2:-0.1,samples=20,line width=2.5pt]
 {0};
 \addplot[black,mark=none,
 domain=-0.1:0.1,samples=20,line width=2.5pt]
 {x+1};
 \addplot[black,mark=none,
 domain=0.1:2,samples=20,line width=2.5pt]
 {1};
 \end{axis}
 \end{tikzpicture}
 \end{center}

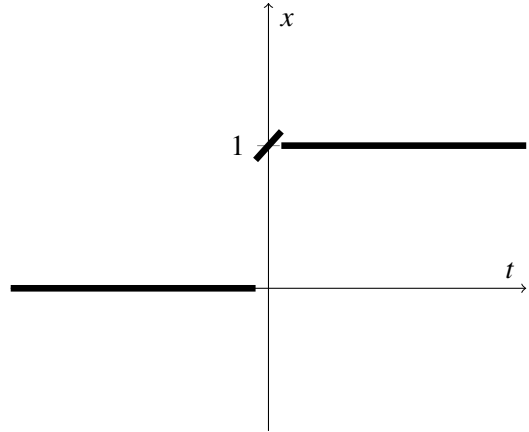
Fig. 4: A solution $x:\widehat{\mathbb{R}}\rightarrow \widehat{\mathbb{R}}$
of the differential equation $x^{\prime }(t)=\left\{ 
\begin{array}{c}
1,{\ }\text{{if \ }}t\in m_{\approx }(0) \\ 
0,{\ }\text{{if \ }}t\notin m_{\approx }(0)%
\end{array}%
\right. .\smallskip $

Notice \hspace{1pt}that $\hspace{1pt}x:\widehat{\mathbb{R}}\rightarrow 
\widehat{\mathbb{R}}$ \hspace{1pt}is \hspace{1pt}an \hspace{1pt}%
indiscernible \hspace{1pt}extension \hspace{1pt}of \hspace{1pt}the \hspace{%
1pt}well-known \hspace{1pt}\textit{Heaviside\vspace{0.2cm}}\newline
\textit{\ function }$H:\mathbb{R}\rightarrow \mathbb{R}$, defined by $%
H\left( \tau \right) :=\left\{ 
\begin{array}{c}
1,{\ }\text{{if \ }}\tau \geq 0 \\ 
0,{\ }\text{{if \ }}\tau <0%
\end{array}%
\right. $.

\section{Conclusion}

The purpose of this work was not to provide a tool to use the concept of 
\textit{actual infin\nolinebreak itesimal} as an alternative to the $%
\varepsilon\mspace{-3mu}\textsl{-}\delta $ definition of \textit{limit}. In fact, we use the
concept of \textit{infinitesimal} (and the concepts of \textit{shadow}, 
\textit{monad}, \textit{indiscernibility}) in the \textit{mode of actuality}
(in loose terms, the \textit{mode of} $\widehat{\mathbb{R}}$, without a
definition of \textit{limit}), and the usual definition of \textit{limit} in
the \textit{mode of potentiality} (in loose terms, the \textit{mode of} $%
\mathbb{R})$. It is our strong conviction that the \textit{modes of actuality%
} and \textit{potentiality} are both necessary (occasionally together, as in
the definition of \textit{differentiability}) to a Calculus suitable not
only for mathematicians, but also for experimental scientists. We must keep
in mind that physicists and engineers need the concept of \textit{limit},
and accept the usual $\varepsilon\mspace{-3mu}\textsl{-}\delta$ definition (though they use it
as little as possible, as most mathematicians), but they also want to use
the heuristic and computational power of \textit{actual infinitesimal methods%
}.

Five other features of this work are worth mentioning:\vspace{-6pt}\smallskip

$\mathbf{c}_{{\mathbf{1}}})$ \textit{The use of explicit actual
infinitesimals}.

${\mathbf{c}}_{{\mathbf{2}}})$ \textit{The local coincidence of the graph of
a function f}, \textit{differentiable at} $\xi _{0}\in \mathbb{R}$, \textit{%
with its tangent at} $\left( \xi _{0},\hspace{1pt}f(\xi _{0})\right) $.

${\mathbf{c}}_{{\mathbf{3}}})$ \textit{The global-local nature of monads of
points}.

${\mathbf{c}}_{{\mathbf{4}}})$ \textit{The set-theoretic and topological
properties of monads of subsets of} $\widehat{\mathbb{R}}$.

${\mathbf{c}}_{{\mathbf{5}}})$ \textit{The sets we use are those of} \textbf{%
ZFC} (\textit{Zermelo-Fraenkel Set Theory with the Axiom of Choice})\textit{%
, without any distinction between internal and external sets.}

${\mathbf{c}}_{{\mathbf{1}}})$ is a positive answer to the uneasiness caused
by the nonexplicit charac\nolinebreak ter of nonnull infinitesimals in 
\textit{Non-standard Analysis} (see, for example, Alain Connes' criticism in 
$\left[ 3\right] ,$ $\S 2,$ p. 211).

We believe that a generalization of ${\mathbf{c}}_{{\mathbf{2}}})$ is
instrumental in \textit{differential geom\nolinebreak etry}, especially for
the definition of the \textit{tangent space} to a manifold at a certain
point.

${\mathbf{c}}_{{\mathbf{3}}})$ was already used in the differential
treatment of some \textit{singularities}, but we are convinced of its
usefulness in the treatment of many others, in the area of \textit{%
differential equations}. Moreover, the fact that $ m_{\approx}(0) $ contains the \textit{real Hilbert space} $ l^{2} $ is very interesting since this space is isomorphic and isometric to any \textit{separable real Hilbert space}.

As to ${\mathbf{c}}_{{\mathbf{4}}})$, the \textit{set-theoretic} and \textit{%
topological} properties of monads of subsets of $\widehat{\mathbb{R}}$ seem
to reveal a pattern extensible to other areas of mathematics.

${\mathbf{c}}_{{\mathbf{5}}})$ is a positive answer to one major difficulty
encountered by non-standard analysts (especially those who work within the
framework of \textit{Internal Set Theory}): \textit{external} \textit{sets}.

Although this article concerns the \textit{differential calculus}, its
fundamental concepts can also be applied to the \textit{integral calculus}
(the work already done and its developments will be published in a future
article).

\bigskip

\pagebreak \textit{References:}

[1] J. Bell, \ \textit{A Primer of Infinitesimal Analysis}, Cambridge University Press, Cambridge, 1998.

[2] J. Bell, \ \textit{The Continuous and the Infinitesimal in mathematics and philosophy}, Polimetrica, Monza, 2006.

[3] A. Connes, Brisure de sym\'{e}trie spontan\'{e}e et g\'{e}om\'{e}trie du
point de vue spectral, \textit{Journal of Geometry and Physics}, \textbf{23}
(1997), 206--234.

[4] K. Hrba\v{c}ek, Axiomatic foundations for nonstandard analysis, \textit{%
Fund. \hspace{-2pt}Math}., \textbf{98} (1978), 1--19.

[5] K. Hrba\v{c}ek, Nonstandard set theory, \textit{Amer. Math. Monthly}, 
\textbf{86} (1979), 659--677.

[6] P. Mancuso, \textit{Philosophy of Mathematics \& Mathematical Practice
in the Seventeenth Century}, Oxford University Press, New York-Oxford, 1996.

[7] I. Moerdijk and G. E. Reyes, \textit{Models for} \textit{Smooth
Infinitesimal Analysis}, Springer-\newline
-Verlag, Berlin-Heidelberg-New York-Tokyo, 1991.

[8] E. Nelson, Internal Set Theory: a new approach to Nonstandard Analysis, 
\textit{Bulletin of the American Mathematical Society}, \textbf{83} (1977),
1165-1198.

[9] E. Nelson, \textit{Radically Elementary Probability Theory}, Annals of
Mathematics Stud\nolinebreak ies, PUP, Princeton, New Jersey, 1987.

[10] A. Robinson, Non-standard Analysis, \textit{Proc.\hspace{-4pt} of the
Royal Academy of Sciences,} \textit{Amsterdam}, ser. A, \textbf{64} (1961),
432-440.

[11] A. Robinson, \textit{Non-standard Analysis}, PUP, Princeton, New Jersey,
rev.ed. 1996.

[12] S. Sambursky, \textit{Physics of the Stoics}, PUP, Princeton, New
Jersey, first paperback ed. 1987.

\end{document}